\documentclass[titlepage]{amsart}
\usepackage{graphicx}
\usepackage[all]{xy}
\usepackage{amsthm}
\usepackage{amsmath}
\usepackage{amssymb}
\usepackage{epsf}
\usepackage{amscd}
\newcommand{\ts}{\tilde{\Sigma} }
\newcommand{\hs}{\hat{\Sigma}}
\newcommand{\mm}{$M(\ts , \Sigma )$}
\newcommand{\T}{$\mathbf{T_{z,y}^{c}}$}
\newcommand{\delS}{{\partial \Sigma}}

\newtheorem{lemma}{Lemma}
\newtheorem{theorem}{Theorem}

\newtheorem*{maintheorem}{Main Theorem}

\theoremstyle{remark}
\newtheorem{remark}{Remark}
\theoremstyle{definition}
\newtheorem{definition}{Definition}
\theoremstyle{proposition}
\newtheorem{proposition}{Proposition}
\theoremstyle{example}
\newtheorem{example}{Example}
\title{On the Lego-Teichm\"{u}ller game for finite $G$-cover}
\author{Tanvir Prince}
\begin{document}
\maketitle
\begin{abstract}
 Given a smooth, oriented, closed surface $\Sigma$ of genus zero, possibly with boundary, let $\ts  \longrightarrow \Sigma$ be a given $G$-cover of $\Sigma$, where $G$ is a given finite group. Let $S_{n}$ denote the standard sphere with $n$ holes. There are many ways of gluing together several $G$-cover of $S_{n}$ to construct the $G$-cover $\ts \longrightarrow \Sigma$, of $\Sigma$. We let M($\ts ,\Sigma$) be the set of all ways to construct the given $G$-cover, $\ts \longrightarrow \Sigma$, of $\Sigma$ from gluing of several $G$-covers of $S_{n}$,
 here $n$ may vary. In this paper, we define some simple moves and relation which will turn M($\ts ,\Sigma$) into a connected and simply-connected complex. This will be used in the future paper to construct $G$-equivariant modular functor. This $G$-equivariant modular functor will be an extension of the usual modular functor.
\end{abstract}
\tableofcontents

\section{Introduction}
This paper can be thought of as  an extension of the paper ``On the Lego-Teichm\"{u}ller Game'' by Bojko Bakalov and Alexander Kirillov Jr \cite{BK}. In the paper \cite{BK}, authors gave the description of how to represent a given smooth, oriented surface, possibly with boundary by gluing of several ``simple pieces''. Here ``simple pieces'' are just sphere with $n$ holes. One application of this is the construction of modular functor; it suffices to define the vector spaces associated with sphere with $n$ holes. Then since the behavior of modular functor under gluing is known, this defines a unique vector space associated to any surface. Of course, there are many different ways to represent a given surface as a result of gluing several spheres with $n$ holes, here $n$ can vary. So some natural questions arise; like how one can describe different ways of gluing ``standard pieces'' that give the same surface? Let $M(\Sigma)$ be the set of all such way of getting the surface $\Sigma$ from the ``simple pieces''. In the paper \cite{BK}, they described some simple moves or edges and some relations among them which turned $M(\Sigma)$ into a connected and simply-connected complex.\\

This definition of modular functor can be extended. Let $G$ be a finite group. And let $\ts  \longrightarrow \Sigma$ be a given $G$-cover of $\Sigma$. It is possible to extend the definition of modular functor from the surface, $\Sigma$, to the $G$-covers of surface, $\ts  \longrightarrow \Sigma$.This idea will be formalized in later paper. This extended modular functor will also satisfy a similar gluing axiom, just like the regular modular functor, but now we are gluing $G$-covers of surface and not just surface. Thus if we know the value of the extended modular functor on
 the ``simple pieces'', of course we need to know what are these ``simple pieces'' in this case, then this will be enough to define Extended modular functor to any $G$-cover. Since any $G$-cover can be constructed from the
 gluing of these ``simple pieces'' and the behavior of the extended modular functor under gluing is known, this will define a unique value for a given $G$-cover. \\

As in the paper of \cite{BK}, we are faced with similar questions. Let $M(\ts ,\Sigma)$ denote the set of all possible way of gluing together "simple pieces" to construct the given $G$-cover, $\ts  \longrightarrow \Sigma$. Is it possible to define some simple moves to go from one parameterization to the other? Is is possible to define all the relation between these moves? ie describe when a sequence of moves applied to a given parameterization yields the same parameterization? In other word, we are trying to define some simple moves and relations similar to the  paper \cite{BK}, which will turn $M(\ts ,\Sigma)$ into a connected and simply-connected complex. This paper runs side by side with the paper \cite{BK}, although there are some differences. It is recommended that readers first read the paper \cite{BK} before reading this one. This paper only deals with the case when the base surface, $\Sigma$, has genus zero. The case of positive genus will be considered in the subsequent papers.

\section{Some basic definition}
Let us start with some basic definitions.

\subsection{Extended Surface}\label{es}
\begin{definition}
An extended surface is a compact, smooth, oriented, closed surface, $\Sigma$, possibly with boundary and also comes with a choice  of distinguished or marked points on each of its boundary component.
\end{definition}
We denote by $A(\Sigma)$, the set of the boundary components. So an extended surface will be denoted by $(\Sigma,\{p_a\}_{a \in A(\Sigma)})$ where $p_a$ is the choice of marked point on the $a$ th boundary component.\\
Sometimes we will also denote a boundary circle by a Greek letter.
\subsection{$G$-cover of Extended surface \label{gcoverofes}}
$G$ will always denote a finite group, which is given and fixed throughout the whole paper. Let $(\Sigma,\{p_a\}_{a \in A(\Sigma)})$ be an extended surface.
\begin{definition}
By a $G$ cover of $(\Sigma,\{p_a\}_{a \in A(\Sigma)})$, we  mean $(\pi : \ts  \longrightarrow \Sigma, \{\tilde{p_a}\})$ where $(\pi : \ts  \longrightarrow \Sigma)$ is a principal $G$-cover and $\{\tilde{p_a}\}$ are choice of points on the fiber of $p_a$. In other word, $\tilde{p_a} \in \pi ^{-1} (p_a) \mbox{ for all } a \in A(\Sigma)$.
\end{definition}

\subsection{Morphism between $G$-cover of extended surface
\label{morphismbetweengcover}}
\begin{definition}
Given two $G$ covers of $(\Sigma,\{p_a\})$ say $(\tilde{\pi }: \ts  \longrightarrow \Sigma, \{\tilde{p_a}\})$ and
$(\hat{\pi }: \hat{\Sigma} \longrightarrow \Sigma, \{\hat{p_a}\})$, by a morphism between them we mean a homeomorphism $f:\ts \rightarrow \hat{\Sigma}$ so that the following conditions are satisfied:\\
\begin{flushleft}
 i) $f(\tilde{p_a}) = \hat{p_a}$ for all $a \in A(\Sigma)$\\
 ii) $\hat{\pi}f = f_{*} \tilde{\pi}$, here $f_{*}: \Sigma \to \Sigma$ is the homeomorphism we get by restriction of $f$ to $\Sigma$ \\
 iii) $f$ preserves the action of $G$ on each fiber. \\
\end{flushleft}
\end{definition}

See the diagram below:\\

\begin{displaymath}
\xymatrix{\ts  \ar[d]_{\tilde{\pi}} \ar[r]^{f} & \hs
\ar[d]^{\hat{\pi}} \\
\Sigma \ar[r]_{f_{*}} & \Sigma }
\end{displaymath}

\begin{remark}
Sometimes, we require $f_{*}: \Sigma \to \Sigma$  to be the identity.
\end{remark}

Although in this definition, we leave the base space fixed, there is no
need to do this. We can easily defined morphism between two $G$ covers
where the base space is not fixed.
\begin{definition}
 Let $(\tilde{\pi }: \ts
  \longrightarrow
\Sigma, \{\tilde{p_a}\})$ be a $G$ cover of $(\Sigma,\{p_a\})$ and let
$(\hat{\pi }: \hat{\Sigma} \longrightarrow {\Sigma'}, \{\hat{p_a}\})$
 be a
$G$ cover of $({\Sigma'},\{{p_a'\}})$. By a morphism between them, we
 mean
a homeomorphism $f:\ts \rightarrow \hat{\Sigma}$ so that it satisfies
 the
following condition:\\
\ i) $f(\tilde{p_a}) = \hat{p_a}$ for all $\tilde{p_a}$ \\
\ ii) $f$ commutes with the action of $G$.\\
\end{definition}
Because of the second condition, it is easily seen that $f$ descends to
 a
homeomorphism on the base surface by $f_*: \Sigma \rightarrow
 {\Sigma'}$
where we defined $f_*(b) = \hat{\pi}f(\tilde{b})$ where $\tilde{b} \in
{\tilde{\pi}}^{-1} (b)$ . It is easily seen that $f_*$ does not depend
 on
the choice of $\tilde{b}$. Also $\hat{\pi}f = f_*\tilde{\pi}$.\\
In other words the diagram below is commutative:\\
\begin{displaymath}
\xymatrix{\ts \ar[r]^{f} \ar[d]_{\tilde{\pi}} & \hs \ar[d]^{\hat{\pi}}
 \\
\Sigma \ar@{.>}[r]^{f_{*}} & \Sigma'}
\end{displaymath}
\begin{remark}
Here we introduce the notation, $\mbox{ Mor}_{\Sigma}(\ts , \ts)$ to denote all the morphisms between the $G$-cover $\ts$ and $\ts$ so that the induced map on $\Sigma$ is identity. Also we use the notation $\mbox{ Mor}(\ts , \hs)$ to denote all morphism between $G$-cover $\ts$ and $\hs$ where the induced map on the base surface can be anything.
\end{remark}

\subsection{Orientation of the boundary circle} \label{orientationoftheboundarycircle}
\begin{figure}
\includegraphics{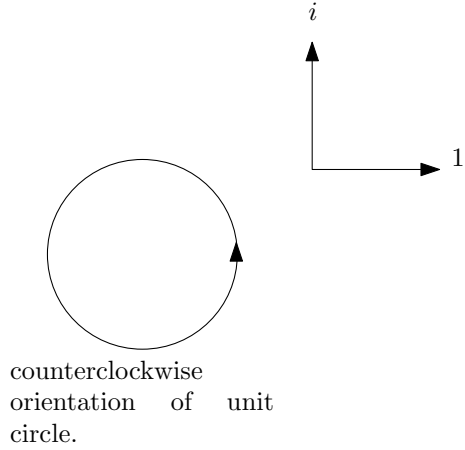}
\caption{orientation on the complement of the unit disk} \label{orientation}
\end{figure}

The orientation of the extended surface, $\Sigma$, naturally induces orientation on the boundary circle. We want to explain this in a little detail. Let $D = \{z \in \mathbb{C} | |z| < 1 \}$. First, we fix an orientation on the complement of $D$. This orientation is given by the choice of the basis, $\{1,i\}$,on $\mathbb{C}$, and the counterclockwise orientation on the unit circle. See the figure \ref{orientation}.
Now for each boundary circle of $\Sigma$, we choose a small neighborhood around the circle. This neighborhood is homeomorphic to $E = \{z \in \mathbb{C} | 1 \leq |z| \leq 1+\epsilon\}$ for some $\epsilon$. For any such homeomorphism which preserve the orientation of the surface, there is an unique choice of orientation of the boundary circle of $\Sigma$ so that the map also preserve the orientation on the boundary. This gives an orientation to each boundary component of $\Sigma$.\\
\subsection{Monodromy \label{monodromy}}
\subsubsection{Definition of Monodromy}
Monodromy will be an important tool to describe $G$-cover. In fact soon we will prove that two $G$-covers of a surface of genus 0 are isomorphic iff their monodromy is same.
\begin{definition}
\begin{description}
    \item{[Monodromy around a boundary circle:] Let $(\pi : \ts
    \longrightarrow \Sigma, \{\tilde{p_a}\})$ be a $G$-cover of
    $(\Sigma,\{p_a\}_{a \in A(\Sigma)})$. Consider the $a$ th boundary
    circle, $S$, where the base point on $S$ is $p_a$ and the point on
 the
    fiber above is $\tilde{p_a}$. This $S$ has an orientation which comes from the orientation of the surface. Let $\alpha: [0,1] \rightarrow S$ be
 a
    parameterization of the boundary circle $S$ which also preserves the
    orientation of $S$. We also assume that
 $\alpha(0)=\alpha(1)=p_a$.
    Then there is a unique lifting of $\alpha$ to the $G$-cover, say
    $\tilde{\alpha}:[0,1]\rightarrow \ts $ such that
    $\tilde{\alpha}(0)=\tilde{p_a}$. We define the monodromy, $m \in
 G$, of
    this $a$ th boundary circle to be that unique element of $G$ such
 that
    \begin{equation}\label{monodromyeq}
        m \tilde{\alpha}(0) = \tilde{\alpha}(1)
    \end{equation}}
\end{description}
\end{definition}
\begin{lemma}\label{monodromylemma}
    Monodromy does not depend on the choice of parameterization
\end{lemma}
\begin{proof}
    This is not hard to see. From the definition, we see that the
 monodromy
    only depends on the starting point of the lifting; if we know the
    starting point of the lifting, everything else is determined by the
    cover, including the end point of the lifting. And monodromy only depends on
 the
    starting point and the ending point of the lifting. The details are
    left for the reader.
\end{proof}

\subsubsection{Relation Between Monodromy and Lifting of a Map}
Given two $G$-covers and a map between their base surfaces, we want to
 know
under what condition this map can be lifted to the covers. Monodromy
 helps
us to partially answer this question. We have the following lemma.
\begin{lemma} \label{liftinglemme}
    Let $(\tilde{\pi }: \ts  \longrightarrow \Sigma, \{\tilde{p_a}\})$
 be a
    $G$ cover of $(\Sigma,\{p_a\})$ and let $(\hat{\pi }: \hat{\Sigma}
    \longrightarrow {\Sigma'}, \{\hat{p_a}\})$ be a $G$ cover of
    $({\Sigma'},\{{p_a'}\})$. Here we also assume the base surfaces to
 be
    connected. Let $f:\Sigma \rightarrow {\Sigma'}$ be a homeomorphism
 of
    the base surfaces which maps marked points to marked points, that is
    $f(p_a)={p_a'}$.
    See the diagram below:\\
    \begin{displaymath}
        \xymatrix{\ts \ar@{.>}[r]^{\tilde{f}} \ar[d]_{\tilde{\pi}} &
 \hs
        \ar[d]^{\hat{\pi}} \\
        \Sigma \ar[r]^{f} & \Sigma'}
    \end{displaymath}
Then
\begin{enumerate}
    \item{The lifting of $f$ to $G$-covers $\ts  \longrightarrow \hs$
 is
        unique, if it exists at all.}
    \item{If in addition, we assume that the genus of both base surfaces are 0, then $f$ can be lifted iff monodromy of the two $G$-covers match. That is $$m((\partial \Sigma)_{i}) = m(f(\partial \Sigma)_{i})$$
         where $m((\partial \Sigma)_{i})$ = monodromy of $(\partial \Sigma)_{i}$ boundary circle of $\Sigma$.\\
         and $m(f(\partial \Sigma)_{i})$ = monodromy of $f(\partial \Sigma)_{i}$ boundary circle of $\Sigma'$.\\
         Note that the homeomorphism, $f:\Sigma \rightarrow {\Sigma'}$, maps the boundary components of $\Sigma$ to the boundary components of $\Sigma'$.}

\end{enumerate}
\end{lemma}
\begin{proof}
    (1) is obvious. To be more specific, let $f_1$ and $f_2$ be two
 lifting
    of $f$. From the definition of lifting and morphism of $G$ covers
 (see
    \ref{morphismbetweengcover}), $f_1$ and $f_2$ must satisfy\\
    \begin{center}
        $\hat{\pi}f_1 = \hat{\pi}f_2 = f\tilde{\pi}$
    \end{center}
    So if $x \in \Sigma$ then both $f_1$ and $f_2$ maps the fiber above
 the
    $x$ to the fiber above the $f(x)$. Moreover from the definition of
    morphism between $G$-cover we must have
        $f_1(\tilde{p_a})=f_2(\tilde{p_a})=\hat{p_a}$\\
    Also both $f_1$ and $f_2$ commute with the action of $G$. This
    information guarantees that they must agree on all the fibers of same
    connected component. But since our base surfaces are connected they
    must agree everywhere. So $f_1$ = $f_2$.\\

For (2), we will use the following standard proposition of $G$ cover:\\
\begin{proposition}
Let $\tilde{\pi}: \ts \to \Sigma$ and $\hat{\pi}: \hs \to \Sigma'$ be two $G$ covers and $f: \Sigma \to \Sigma'$ be a map between the base surface. Then this map, $f$ can be lifted to the $G$ cover if and only if the monodromy of every loop, $\alpha$,  in $\tilde{\pi}: \ts \to \Sigma$ is equal to the monodromy of the corresponding loop, $f(\alpha)$, in  $\hat{\pi}: \hs \to \Sigma'$. And in the case that both of the base surfaces, $\Sigma$ and $\Sigma'$, have genus zero, then $f$ can be lifted to the $G$ covers if and only if the monodromy of the corresponding boundary circles match.
\end{proposition}
\begin{proof}
A proof of this or some equivalent statements can be found in many standard books on topology. For example see section 1.3 of \cite{AH}. 
\end{proof}
From the above proposition, the statement 2 of our lemma easily followed.

\end{proof}

\subsection{Gluing of two $G$-covers of extended surfaces}
\begin{definition}
Let $(\tilde{\pi }: \ts  \longrightarrow \Sigma, \{\tilde{p_a}\})$
 be a
    $G$ cover of $(\Sigma,\{p_a\})$ and let $(\hat{\pi }: \hat{\Sigma}
    \longrightarrow {\Sigma'}, \{\hat{p_a}\})$ be a $G$ cover of
    $({\Sigma'},\{{p_a'\}})$. Let $s_{i,j}: (\partial \Sigma)_{i} \to (\partial \Sigma')_{j}$ be an orientation reversing map of the $i$th boundary circle of $\Sigma$ to the $j$ th boundary circle of $\Sigma'$ so that $s_{i,j}(p_{i}) = p_{j}'$. Then we define the gluing of these two $G$ covers under $s_{i,j}$ to be the following $G$ cover:\\
$$\tilde{\pi} \bigsqcup_{s_{i,j}} \hat{\pi } : \ts \bigsqcup_{s_{i,j}} \hs \longrightarrow \Sigma \bigsqcup_{s_{i,j}} \Sigma'$$
where $\Sigma \bigsqcup_{s_{i,j}} \Sigma'$ is the surface obtained by identifying points on $(\partial \Sigma)_{i}$ to the points on $(\partial \Sigma')_{j}$ through the map $s_{i,j}$. And $\ts \bigsqcup_{s_{i,j}} \hs$ is the $G$ cover obtained by identifing a point $(t,g)$, where $t \in (\partial \Sigma)_{i}$ and $g \in G$, to the point $(s_{i,j}(t), g)$.
\end{definition}
\subsubsection{How to glue two $G$-covers}\label{gluingofgcover}
Such a gluing not always exists, but if it exists, it is unique. So we have the following lemma.
\begin{lemma}\label{gluinglemma}
    Let $(\tilde{\pi }: \ts  \longrightarrow \Sigma, \{\tilde{p_a}\})$
 be a
    $G$ cover of $(\Sigma,\{p_a\})$ and let $(\hat{\pi }: \hat{\Sigma}
    \longrightarrow {\Sigma'}, \{\hat{p_a}\})$ be a $G$ cover of
    $({\Sigma'},\{{p_a'\}})$. Take $b \in A(\Sigma)$ and $c \in
    A(\Sigma')$. We want to glue ${(\partial {\Sigma})}_b $ and
    ${(\partial {{\Sigma'}})}_c $. Then :\\
    \begin{enumerate}
        \item{If the gluing exists, there is a unique way, up to isomorphism of $G$ covers, to glue these
            $G$-covers}
        \item{Gluing is possible iff the monodromy, $m_b$, of
 $(\delS)_b$
        and the monodromy, $m_c$ of $(\delS')_c $
        are inverse of each other. That is $m_b m_c = 1$.}
        \end{enumerate}
\end{lemma}
\begin{proof}
    Not only we need to glue $(\delS)_b $ and $(\delS')_c $ but we also
    need to glue the cover above it. See the diagram on figure \ref{glue1}.
    \begin{figure}
        \includegraphics{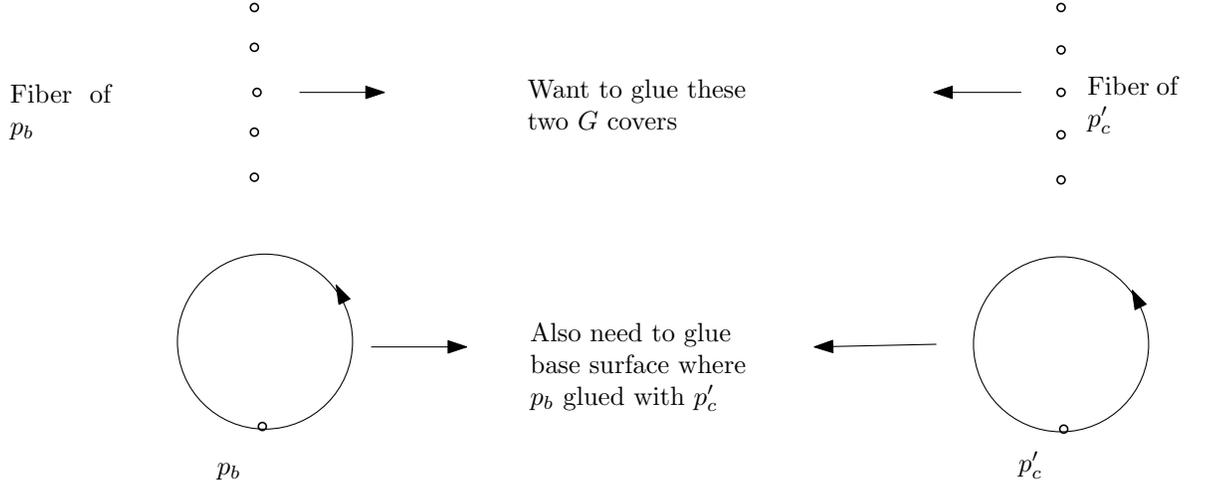}
        \caption{Fiber of the two boundary circles we want to glue}\label{glue1}
    \end{figure}
    We identify $p_b \in (\delS)_b $ with ${p_c'}\in {(\delS')}_c $ and
    this will basically tells us how to glue ${(\delS)}_b $
    with $(\delS')_c $ since, first of all the orientation of
    the boundary circle comes from the orientation of the surface, and
 the
    set of all orientation preserving homeomorphisms from ${(\delS)}_b $
    to $(\delS')_c $ which map $p_b$ to ${p_c'}$ is
    homotopic to each other. This is how we glue ${(\delS)}_b $ with
    ${(\delS')}_c $ . Now what about the cover? Let us take a
    small neighborhood of $p_b$ and ${p_c'}$ on the circle
 ${(\delS)}_b $
    and ${(\delS')}_c $ so that the fiber above these neighborhoods of
    circles break up into disjoint pieces and each piece maps
    homeomorphically    by $\tilde{\pi}$ and $\hat{\pi}$ to these
    neighborhoods of circle. see the diagram on figure \ref{gluing2}\\
    \begin{figure}
        \includegraphics{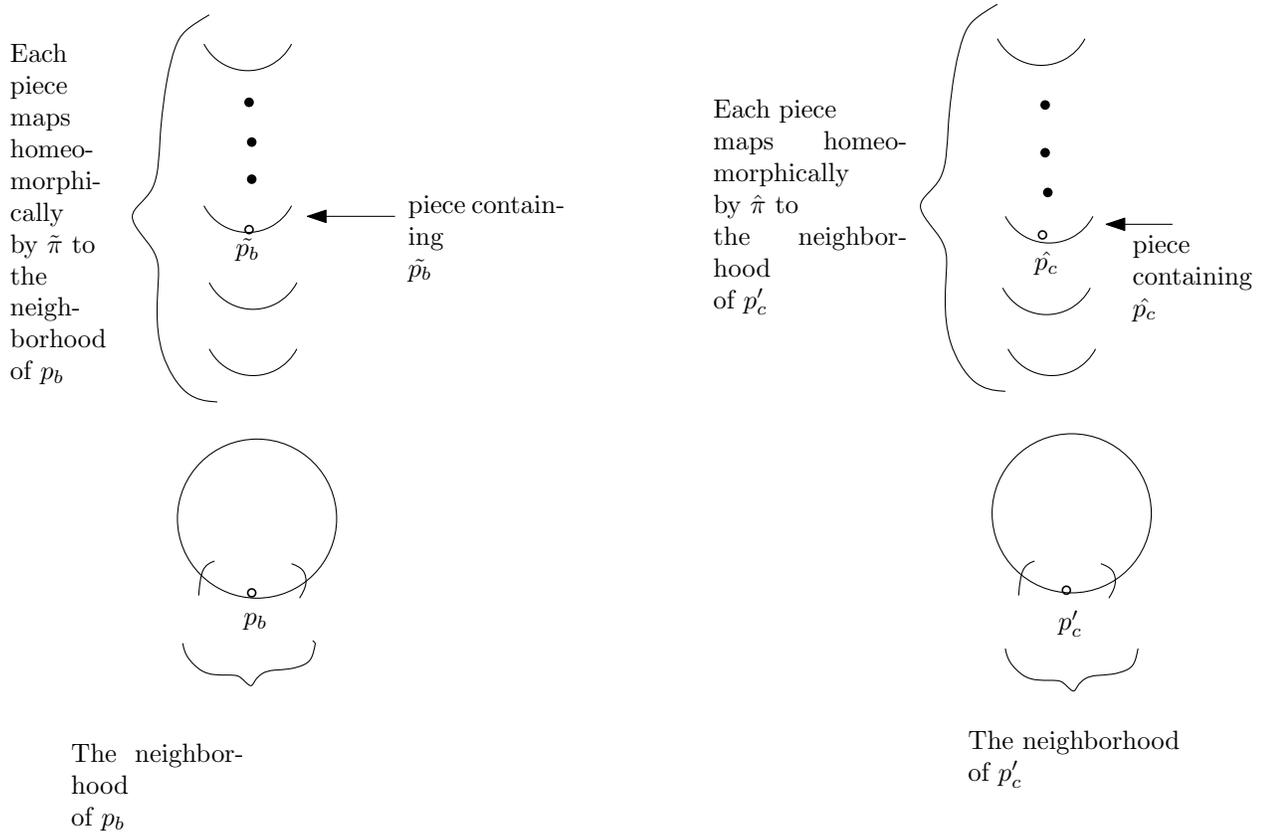}
        \caption{Neighborhood of distinguished points and the fiber above it}\label{gluing2}
    \end{figure}
    One of the pieces on the left contains $\tilde{p_b}$ and one of the
    piece on the right contains $\hat{p_c}$. Also $G$ acts simply
    transitively on these pieces. By definition, we require the piece
    containing $\tilde{p_b}$ to be glued with the piece containing
    $\hat{p_c}$ so that the action of $G$ commutes with the gluing.
 This
    fixes how the pieces of $G$-cover containing $\tilde{p_b}$ will be
    glued to the pieces of $G$-cover containing $\hat{p_c}$. Now we
 move
    around the circle and repeat the same process until we cover the
 whole
    circle, $(\delS)_b $ and $(\delS')_c$.

    This tells us that there is at most one way to glue. Note that not
    always we can glue two $G$- covers. For example take $|G| = 2 $.
 Take
    one $G$-cover to be the trivial $G$-cover of the circle and the
 other
    $G$-cover to be the double cover of the circle. These two
 $G$-covers
    can not be glued.

    The second part of the lemma is left to the reader. Basically one
    sees that as we moved around the whole circle and then glue or we
 first
    glue and then move around the whole circle, in either case we
 arrived
    at the same point since the monodromy are inverse of each other. So
 we
    don't have any problem to glue. Readers can supply the detail.
\end{proof}

\section{Standard block}
\subsection{Standard sphere, $S_{n}$}
\begin{figure}
\includegraphics{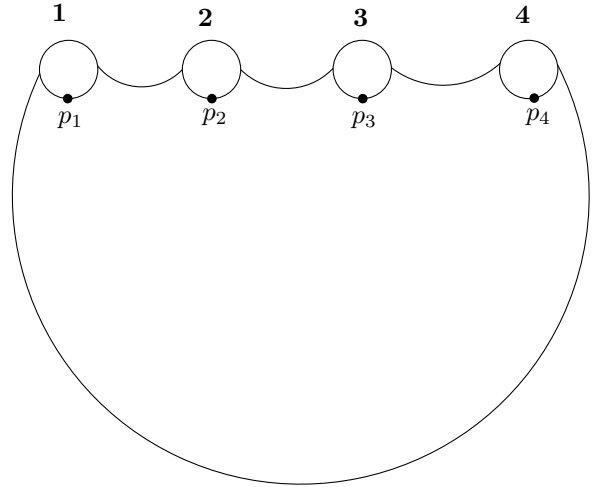}
\caption{A standard sphere with four holes where the marked points are denoted by $p_{1},...p_{4}$}\label{standardsphere}
\end{figure}

For every $n\geq 0$, we define the standard sphere, $S_{n}$, to be the Reimann sphere $\overline{\mathbb{C}}$ with $n$ disks $|z-k|<\frac{1}{3}$ removed and with the marked points being $k-\frac{i}{3}$, here $k = 1,2,3,...n$. Of course, we could replace these $n$ disks with any other $n$ non-overlapping disks with centers on the real line and with marked points in the lower half plane. Any two such spheres are homeomorphic and the homeomorphism can be chosen canonically up to homotopy. Note that the set of boundary components of the standard sphere is naturally indexed by numbers $1,2,...n$. The standard sphere, $S_{4}$, with four holes is shown on figure \ref{standardsphere}
\subsection{Standard blocks, $S_{n}(g_{1},...,g_{n};h_{1},...,h_{n})$}
\begin{figure}
\includegraphics[scale=.6]{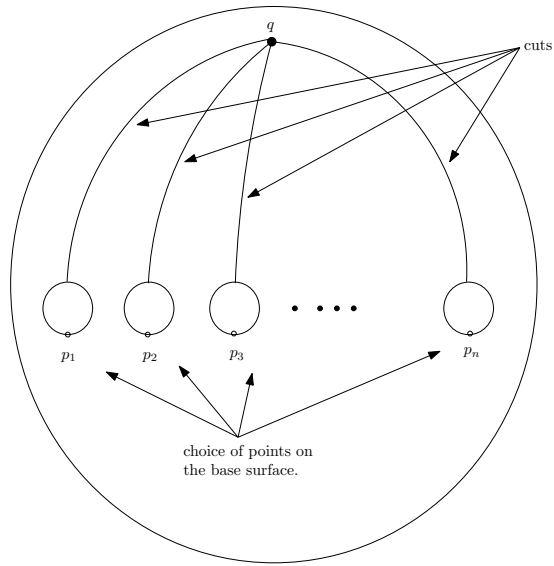}
\caption{Cuts on $S_{n}$}\label{cutsons11}
\end{figure}

\begin{figure}
\includegraphics{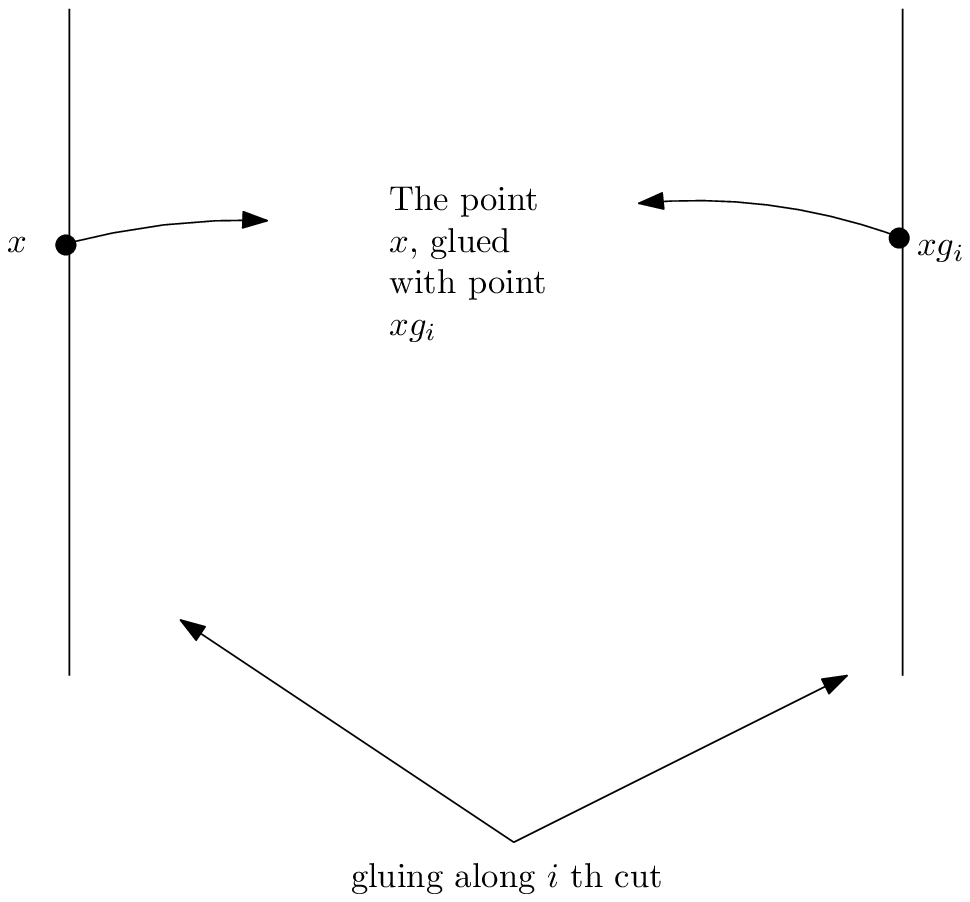}
\caption{A point $x$ on the left hand side of the fiber of the $i$th cut is glued to $xg_{i}$ on the right hand side }\label{gluingithcut}
\end{figure}

In the paper \cite{BK}, where $|G|=1$, standard blocks are just these standard spheres,$S_{n},n=1.2...$. So we need to extend this definition to the general case where $G$ is a finite group. To do this we start with the following construction. Let us start with a standard sphere with $n$ holes, $S_{n}$ and $2n$ elements from $G$ where we denote these $2n$ elements as $\{g_{1},...,g_{n}\}$ and $\{h_{1},...,h_{n}\}$ and we also required that $g_{1}...g_{n} = 1$.  First we make the cuts on $S_{n}$ as in figure \ref{cutsons11}.
Here the point $q \in S_{n}$ in figure \ref{cutsons11} is the point at $\infty$. In fact $q$ can be chosen to be any point on the upper hemisphere as long as it does not belong to the boundary circles. Then one can easily sees that $S_n \backslash$cuts is simply connected. So $G$-cover of $S_n \backslash$cuts is trivial. In other word, $G$-cover of $S_n \backslash$cuts looks like \\
\begin{displaymath}
\xymatrix{S_n \backslash \mbox{ cuts } \times G \ar[d] \\
S_n \backslash \mbox{ cuts }}
\end{displaymath}
Although there is only one $G$ cover of $S_n \backslash$cuts up to isomorphism, there are a total of $|G|$ many way to identify $G$ cover of $S_n \backslash$cuts with $G$, but the important thing here to notice is that any such identification is isomorphic. So we chose one such identification. Now consider the $i$th cut. We want to glue the fiber on the left hand side of the $i$th cut with the fiber on the right hand side of the $i$th cut. This identification must preserve the action of the group $G$. Thus the identification can be given by a multiplication on the right by some element of the group $G$. We chose this element to be $g_{i} \in G$. That is we glued along the $i$th cut by multiplication of the right by $g_i$. See figure \ref{gluingithcut}.
Now we choose $\tilde{p_i} = (p_i,h_i)$ as the choice of point in the fiber of $p_i$. Thus we end up with a $G$ cover of $S_{n}$. The reason we require that $g_1...g_n=1$ is easy to see. Consider a point $(t,x), \mbox{ where } x \in G$, on the $G$-cover. As we moved along each cut, we multiply on the right by $g_i$. After moving along all the cut, we end up with $(\mbox{ some point on of the base surface }, xg_1g_2...g_n)$, but then we must have $x= xg_1g_2...g_n$ or equivalently $g_1...g_n=1$. Thus from the data $g_1,...,g_n$ and $h_1,...,h_n$ where $g_1...g_n=1$, we have constructed a $G$ cover of $S_{n}$.
\begin{definition}
The $G$ cover of $S_{n}$ constructed above is called standard block and will be denoted by $S_n(g_{1},g_{2}....,g_{n} ; h_{1},h_{2}....,h_{n})$ where we always assume that $g_{1}g_{2}...g_{n} = 1$ although we do not always write this explicitly.
\end{definition}
\subsection{When two standard blocks are isomorphic}\label{isostand}
Given two $G$-covers of $S_n$, we want to know when they are isomorphic, hence the following lemma:\\
\begin{lemma}\label{isolemma}
Let $S_n(g_1,....g_n;h_1,...h_n)$ and $S_n(g_1',....g_n';h_1'...h_n')$ be two standard blocks. Then they are isomorphic iff $\exists x \in G$ so that $xg_ix^{-1}=g_i'$ and $h_ix^{-1}=h_i'$ for $i=1...n$. We denote the isomorphism, $S_{n}(\mathbf{g,h}) \to S_{n}(\mathbf{g',h'})$ by $\phi_{x}$.
\end{lemma}
\begin{proof}
First assume $S_n(g_1,....g_n;h_1,...h_n)$ and $S_n(g_1',....g_n';h_1'...h_n')$ are isomorphic. This means there exist an isomorphism of $G$-cover which by definition (see section \ref{morphismbetweengcover} on page \pageref{morphismbetweengcover}) maps $(p_i,h_i) \longrightarrow (p_i,h_i')$ and preserve the $G$-action. So on the fiber, this isomorphism is just the right multiplication by some element $x^{-1}\in G$. Also the following condition must satisfy along the cut:\\
First glue and then apply isomorphism $=$ First apply isomorphism and then glue. \\
In other word, if $t$ is a point on the fiber above the $i$th cut, then we must have\\
    \begin{eqnarray*}
        tg_ix^{-1}=tx^{-1}g_i'\\
        \mbox{ or } g_i'=xg_ix^{-1} \mbox{ for } i=1...n\\
    \end{eqnarray*}
And of course we must have $h_ix^{-1}=h'_i$ from the definition of isomorphism of two $G$-cover.\\
The other direction is easier to show. We define the isomorphism on each fiber through the multiplication on the right by $x^{-1}$. This can be easily seen to satisfy all the required property of isomorphism.
\end{proof}

\subsection{Monodromy of a Standard Block}
The next natural question is to ask what is the monodromy of the standard block, $S_n (g_1,g_2...g_n;h_1,h_2...h_n)$. So we have the following lemma:\\
\begin{lemma}\label{monodromyofstandardblock}
    Let $S_n(g_1,....g_n;h_1,...h_n)$ be a $G$-cover of $S_n$. Here we
    assume the orientation of the boundary circles induces by the outward normal vector according to the right hand rule. See section \ref{orientationoftheboundarycircle} on page \pageref{orientationoftheboundarycircle} for the discussion of how the orientation of the surface induces orientation on the boundary. See the picture on figure \ref{sso}.
    \begin{figure}
        \includegraphics[scale=.6]{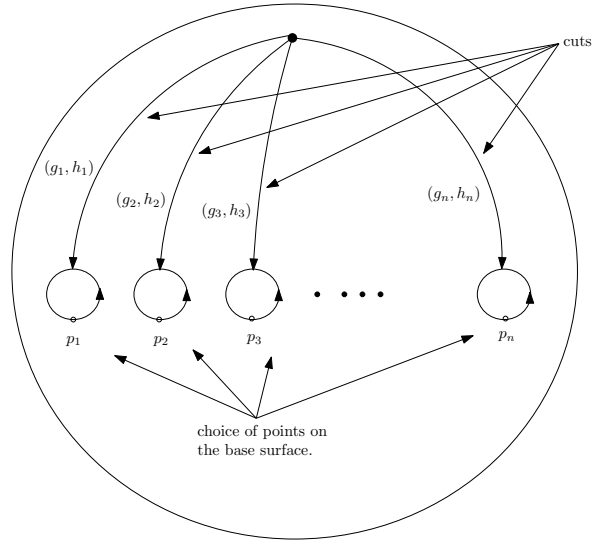}\\
        \caption{observe the orientation of the boundary circles}\label{sso}
    \end{figure}
     Then the monodromy $m_i \in G$ around the $i$th boundary circle is
    given by: $$m_i = h_ig_i^{-1}h_i^{-1}$$

\end{lemma}
 \begin{proof}
    Given a parameterization $\alpha$ : [0,1] $\longrightarrow$ $i$th
    boundary circle of $S_n$, we lift this path starting from
 $(p_i,h_i)$,
    and as we cross the cut labeled by $g_i$ from right to left, we end
 up with $(*,h_ig_i^{-1})$. See the picture on figure \ref{findingmonodromyofs}.
    \begin{figure}
        \includegraphics{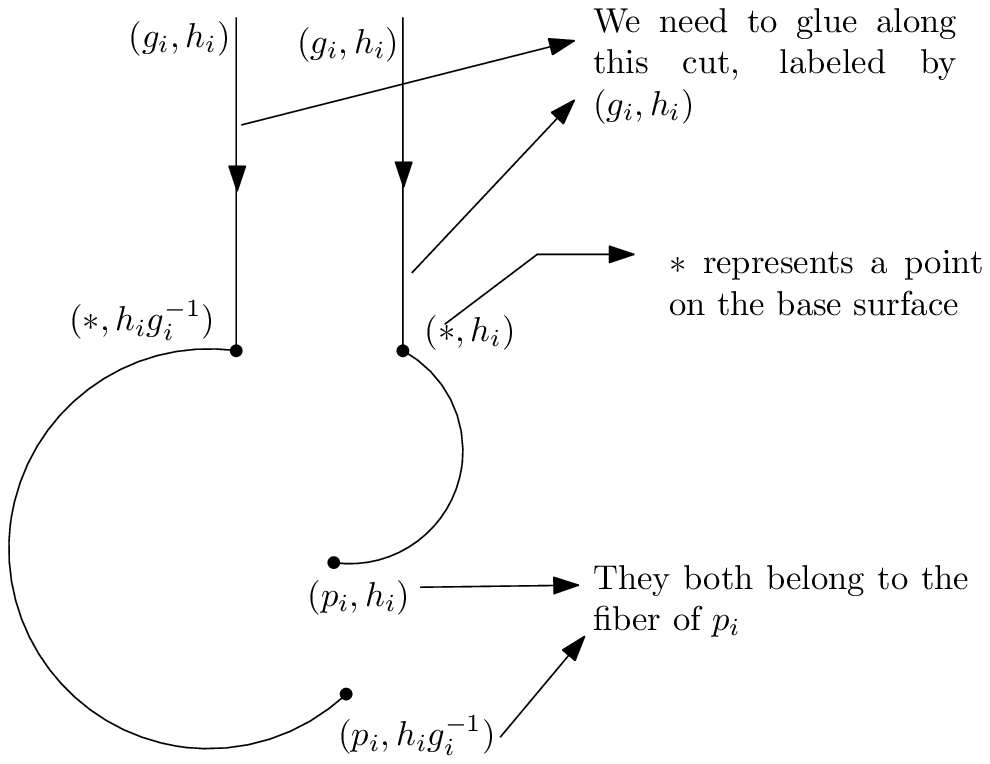}
        \caption{}\label{findingmonodromyofs}

    \end{figure}

    And as we continue all the way, at the end we arrive at the point
    $(p_i,h_ig_i^{-1})$. So from the definition of monodromy, the $i$th
    monodromy $m_i \in G$ is given by\\
    \begin{eqnarray*}
        m_ih_i=h_ig_i^{-1} \mbox{ so }\\
        m_i=h_ig_i^{-1}h_i^{-1}\\
    \end{eqnarray*}
\end{proof}

\subsection{Gluing of two standard blocks}\label{gluenotation}
Let the standard block $S_{n}(g_{1},...,g_{n};h_{1},...,h_{n})$ and $S_{m}(u_{1},...,u_{m};v_{1},...,v_{m})$ are given. We want to know when we can glue these two standard blocks along the $i$ th boundary circle of $S_{n}(g_{1},...,g_{n};h_{1},...,h_{n})$ to the $j$ th boundary circle of $S_{m}(u_{1},...,u_{m};v_{1},...,v_{m})$. Hence we have the following lemma:\\
\begin{lemma}
$i$ th boundary circle of $S_{n}(g_{1},...,g_{n};h_{1},...,h_{n})$ can be glued to the $j$ th boundary circle of $S_{m}(u_{1},...,u_{m};v_{1},...,v_{m})$ iff $h_{i}g_{i}^{-1}h_{i}^{-1} = [v_{j}u_{j}^{-1}v_{j}^{-1}]^{-1}$.
\end{lemma}
\begin{proof}
According to lemma \ref{monodromyofstandardblock} on page \pageref{monodromyofstandardblock}, monodromy of the $i$ th boundary circle of $S_{n}(g_{1},...,g_{n};h_{1},...,h_{n})$ is given by $h_{i}g_{i}^{-1}h_{i}^{-1}$ and the monodromy of the $j$ th boundary circle of $S_{m}(u_{1},...,u_{m};v_{1},...,v_{m})$ is given by $v_{j}u_{j}^{-1}v_{j}^{-1}$. Now according to the lemma \ref{gluinglemma} on page \pageref{gluinglemma}, such a gluing exist if and only if these two monodromy are inverse of each other, that is $h_{i}g_{i}^{-1}h_{i}^{-1} = [v_{j}u_{j}^{-1}v_{j}^{-1}]^{-1}$.
\end{proof}
\begin{remark}
We will use the notation $$S_{n}(g_{1},...,g_{n};h_{1},...,h_{n}) \bigsqcup_{c,h_{i},v_{j}} S_{m}(u_{1},...,u_{m};v_{1},...,v_{m})$$ to indicate that the $i$th boundary of the left standard block is glued along the cut $c$ to the $j$th boundary of the right standard block.
\end{remark}
\subsection{Groupoid}
\begin{definition}
A category is called a groupoid if all of its morphism is invertible.\\
\end{definition}
We list two important examples of groupoid which will be needed later.\\
\begin{example}
 For a fixed $S_n$, consider the category,
$Y_n$, defined in the following way:\\
\ obj($Y_n$)= $G$-covers of $S_n$\\
\ mor($Y_n$)= isomorphisms of $G$-covers which are trivial on $S_{n}$.\\
Since only morphisms are isomorphisms, this is obviously a groupoid.
\end{example}
\begin{example}
 We define the category $T_n$ in the
 following
way:\\
\ obj($T_n$) = \{standard blocks\} = $\{S_{n}(g_1,g_2,....g_n;h_1,h_2...h_n)|g_i \in G \mbox{
 and
} h_i \in G \mbox{ for } i=1...n \mbox{ and } g_1g_2...g_n=1\}$\\
According to lemma in \ref{isostand} on page \pageref{isostand} all the isomorphism between two standard block can be defined in the following way:\\
mor($S_{n}(g_1,g_2,....g_n;h_1,h_2...h_n),
S_{n}(g_1^{'},g_2^{'},....g_n^{'};h_1^{'},h_2^{'}...h_n^{'})$)=$\{x \in
 G |
xg_ix^{-1}=g_i^{'} \mbox{ and } h_ix^{-1}=h_i^{'} \mbox{ for }
 i=1...n\}$.
Since given a morphism $x \in G$ in this category the inverse morphism
 is
$x^{-1}\in G$, every morphism is invertible. So the category $T_n$ is in fact a groupoid.
\end{example}
 We have the following important lemma:\\
\begin{lemma}\label{groupoidlemma}
    The groupoid $Y_n$ is equivalent to the groupoid $T_n$. See example
 1
    and 2 above for the description of $Y_n$ and $T_n$.
\end{lemma}
\begin{proof}
    To prove this lemma, we will use a well known theorem from category theory.
    \begin{theorem}\label{cat}
    Let $A$ and $B$ be two categories, and $F: A \rightarrow B$ is a
    covariant functor, so that the following two conditions hold:\\
    \begin{enumerate}
    \item For any two objects, $X,Y \in \mbox{ Obj }(A)$, the map
        $\mbox{ Mor }_{A} (X,Y) \longrightarrow \mbox{ Mor }_{B}
        (F(X),F(Y))$, induced by $F$, is bijective.
    \item The map $F$ is essentially surjective. That is given any
 object,
        $V \in \mbox{ Obj }(B)$, there exist an object, $U \in \mbox{
 Obj
        }(A)$ so that $F(U) \cong V$. Here the symbol $\cong$ means
        isomorphic.
    \end{enumerate}
    Then $F$ is in fact an equivalence of categories.
\end{theorem}
\begin{proof}
    For a proof of this theorem, see any standard book on category theory, for
    example ``Categories for the working Mathematician'' by Mac Lane, see \cite{C}.
\end{proof}
Now let us come back to the proof of our lemma. We will construct a
 functor
$V: T_n \rightarrow Y_n $ so that $V$ satisfies the two conditions of
Theorem~\ref{cat}. Then this $V$ will define the equivalence between
$T_{n}$ and $Y_{n}$. Given a standard block, $S_{n}(g_1,g_2,....g_n;h_1,h_2...h_n) \mbox{
where }g_i \in G \mbox{ and } h_i \in G \mbox{ for } i=1...n \mbox{ and
 }
g_1g_2...g_n=1$, this standard block, in particular, a $G$ cover of $S_{n}$. That is $$V(S_{n}(g_1,g_2,....g_n;h_1,h_2...h_n)) = \mbox{ the $G$ cover of $S_{n}$ corresponding to the standard block }$$
Similarly, if $\phi_{x}: S_{n}(g_1,g_2,....g_n;h_1,h_2...h_n) \rightarrow S_{n}(g'_1,g'_2,....g'_n;h'_1,h'_2...h'_n)$ is a morphism between two standard blocks, then $V(\phi_{x})$ is the same morphism between the $G$ covers.\\

To show, this functor $V$ satisfies the condition $1$ of the above theorem, we just refer to the lemma \ref{isolemma} on page \pageref{isolemma}. First note that, by this lemma, given any two standard blocks, either there exist an unique morphism, $\phi_{x}$, between them or there is no morphism between them. In either case, the condition $1$ is obviously satisfied.

For condition 2, let $\pi: \tilde{S_{n}} \rightarrow S_{n}$ be a $G$-cover of $S_{n}$. First, we make the same cuts on $S_{n}$ as figure \ref{cutsons11} on page \pageref{cutsons11}. Then the $G$-cover of $S_n \backslash$cuts is trivial. In other word, the $G$-cover of $S_n \backslash$cuts can be identified with \\
\begin{displaymath}
\xymatrix{S_n \backslash \mbox{ cuts } \times G \ar[d] \\
S_n \backslash \mbox{ cuts }}
\end{displaymath}
Now we need to identify the components of $S_n \backslash \mbox{ cuts }
\times G$ with the group $G$. Here we have choice. So we make some choice, it does not matter how we want to do this. Now we basically repeat the same construction when we describe standard blocks, namely, from this trivial cover, to get the original cover we started with, we need to glue the cover along the cuts. See the picture on figure \ref{gluingithcut1}.
\begin{figure}
\includegraphics{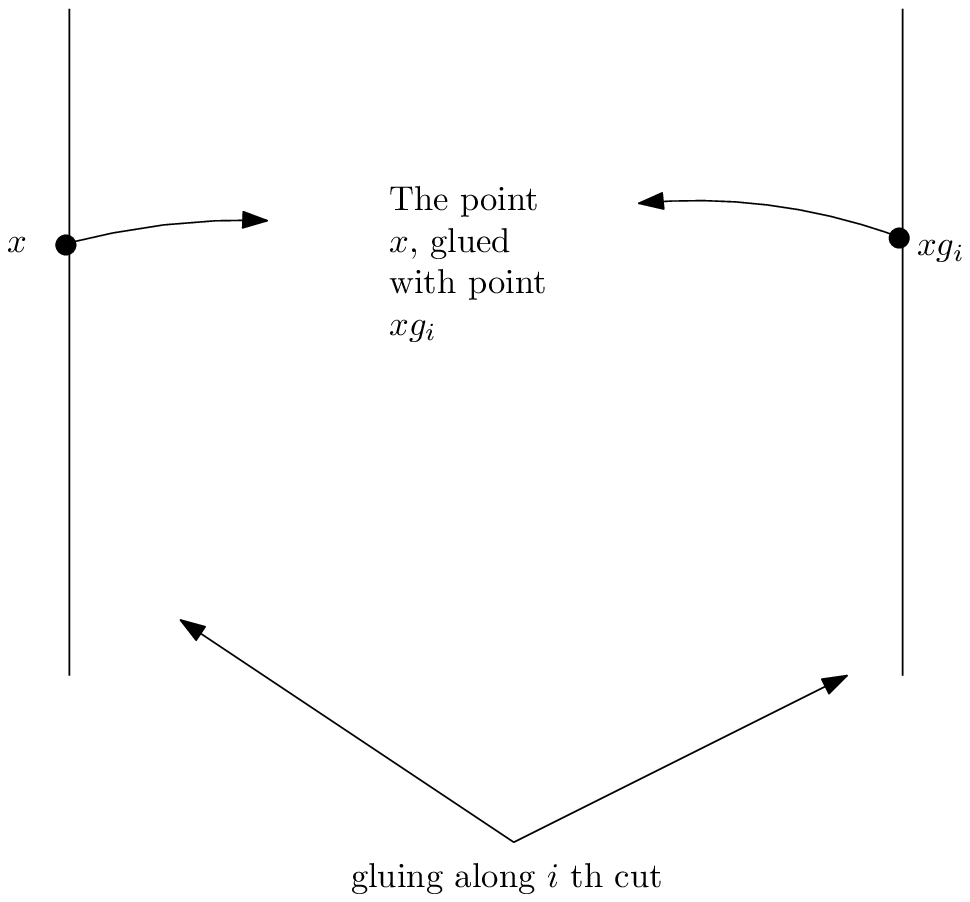}
\caption{}\label{gluingithcut1}
\end{figure}
This gluing must preserve the action of $G$ on the fiber. So to glue along the $i$-th cut, there must exist $g_i \in G$ so that the point $t$ on one side of the cut must glued with the point $tg_i$ from the other side of the cut. Thus for a total of n-cut we get $\{g_1,g_2...g_n \in G\}$. Also to start with, our $G$-cover comes with a point $(p_i,h_i)$ on the fiber above $p_i$. This gives us $\{h_1,h_2...h_n \in G\}$. To show $g_1...g_n=1$ is
easy. Consider a point $(t,x), \mbox{ where } x \in G$, on the $G$-cover. As we moved along each cut, we multiply on the right by $g_i$. After moving along all the cut, we end up with $(\mbox{ some point on of the base surface }, xg_1g_2...g_n)$, but then we must have $x= xg_1g_2...g_n$ or equivalently $g_1...g_n=1$. Thus the $G$ cover $\tilde{S_{n}}$ is isomorphic to the standard block $S_{n}(g_{1},...g_{n};h_{1},...h_{n})$. This shows that the functor $V$ satisfies the second requirement of the above theorem.
\end{proof}
\begin{remark}
In particular, the above theorem shows that every $G$ cover of $S_{n}$ is isomorphic to a standard block.
\end{remark}

\subsection{Review of the parameterization for the case $|G| = 1$}\label{gonecase}
\begin{figure}
\includegraphics{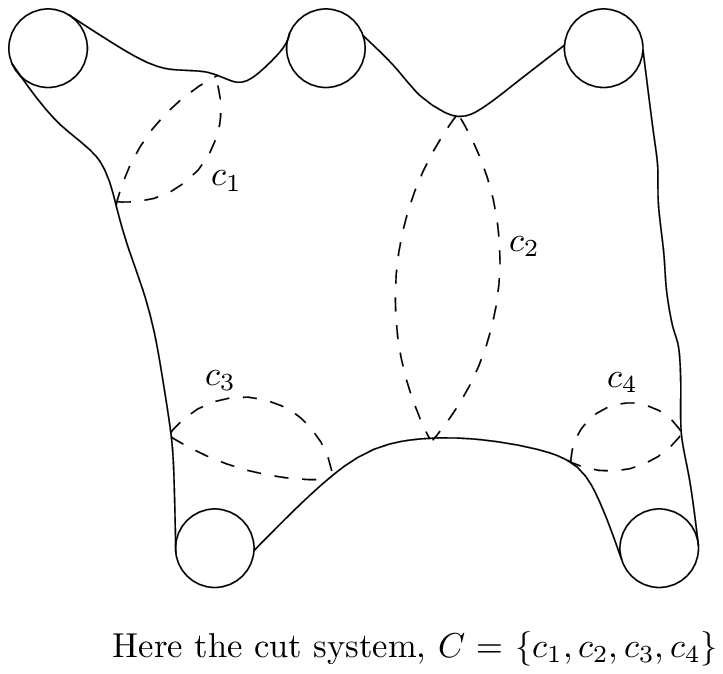}
\caption{An example of a cut system} \label{cutsystem}
\end{figure}

\begin{figure}
\includegraphics{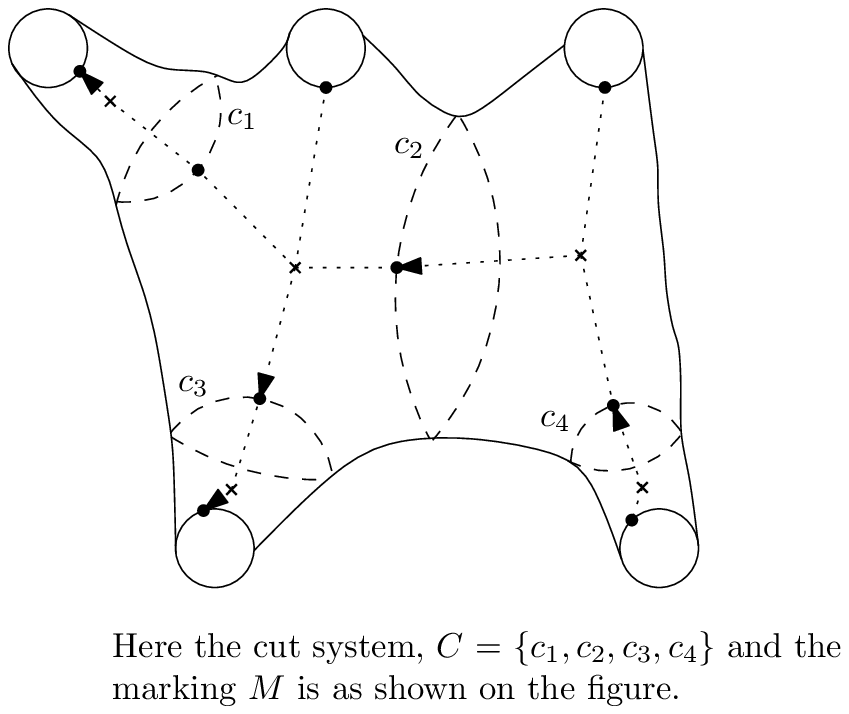}
\caption{An example of a marked surface}\label{markedsurface}
\end{figure}

\begin{figure}
\includegraphics{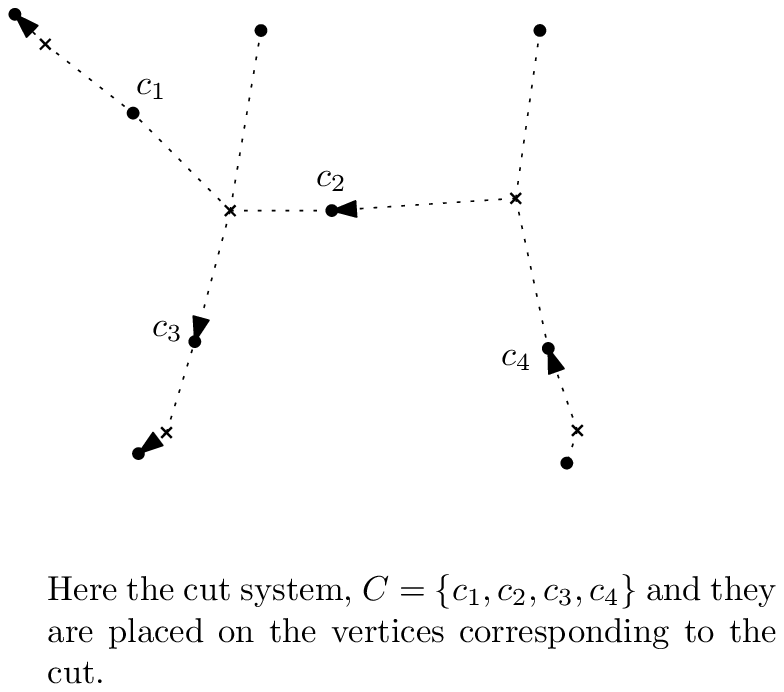}
\caption{Marked graph without the surface}\label{markednosurface}
\end{figure}

For the readers convenience, we will review the concept of parameterization for the trivial case, $|G| =1$, from \cite{BK}. We begin with some definitions.
\begin{definition}
Let $\Sigma$ be an extended surface. A cut system, $C$,  on $\Sigma$ is a finite collection of smooth, simple closed non-intersecting curves on $\Sigma$ such that each connected component of the complement, $\Sigma \backslash C$, is a surface of genus zero. In this paper, we will always assume the surface $\Sigma$ has genus $0$. So in this case, the requirement that each connected component of the complement, $\Sigma \backslash C$, is a surface of genus zero is always satisfied.
\end{definition}
An example of a cut system is given in figure \ref{cutsystem}.
\begin{definition}
Let $\Sigma$ be an extended surface with genus zero. A parameterization without cuts of $\Sigma$ is a homotopy equivalence class of homeomorphisms $\phi : \Sigma \simeq S_{n}$, where $S_{n}$ is the standard sphere with $n$ holes.
\end{definition}
\begin{definition}
Let $\Sigma$ be a an extended surface. A parameterization, $P$, of $\Sigma$ is a collection $(C, \{\phi_{a}\})$, where $C$ is a cut system on $\Sigma$ and $\phi_{a}$ are parameterization without cuts of the connected components $\Sigma_{a}$ of $\Sigma \backslash C$, i.e. homotopy equivalence class of homeomorphisms $\phi_{a}: \Sigma_{a} \simeq S_{n_{a}}$.
\end{definition}
\begin{definition}
Let $S_{n}$ be the standard sphere with $n$ holes. We let $m_{0}$ be the graph on it, shown in figure \ref{standardmarking}, for $n = 5$. This graph has a distinguished edge- one which connected to the boundary component labeled by $1$. This distinguished edge has been marked by an arrow. We call $m_{0}$, the standard marking without cuts on $S_{n}$.\\

Let $\Sigma$ be an extended surface with genus zero. A marking without cuts of $\Sigma$ is a graph, $m$, on $\Sigma$ with one distinguish or marked edge such that $m = \phi^{-1}(m_{0})$ for some homeomorphism $\phi : \Sigma \to S_{n}$. The graphs are considered up to isotopy of $\Sigma$.
\begin{figure}
\includegraphics{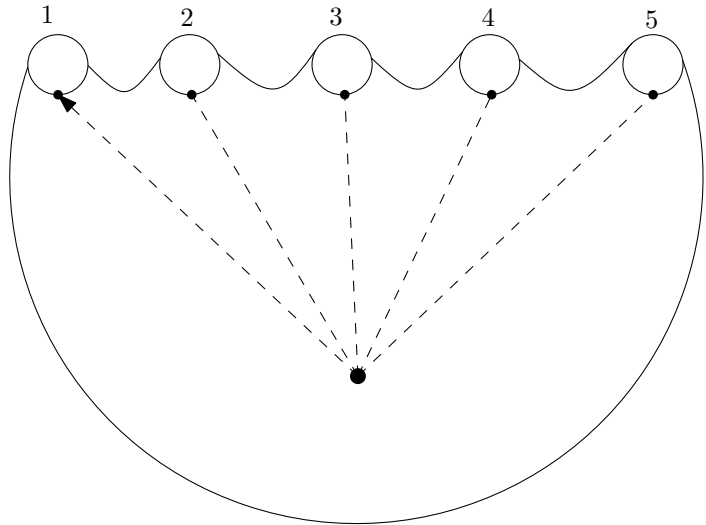}
\caption{Standard marking on $S_{5}$, sphere with $5$ holes}\label{standardmarking}
\end{figure}
\end{definition}
\begin{definition}
Let $\Sigma$ be an extended surface. A marking, $M$ of $\Sigma$ is a pair, $(C,m)$, where $C$ is a cut system on $\Sigma$ and $m$ is a graph on $\Sigma$ with some distinguished edges such that it gives a marking without cuts on each of the connected component of $\Sigma \backslash C$. We will denote the set of all marking of a surface $\Sigma$ modulo isotopy by $M(\Sigma)$. A marked surface is an extended surface, $\Sigma$, together with a marking, $M$ on it.
\end{definition}
An example of a marked surface is shown on figure \ref{markedsurface}.
Sometimes, for the convenience of drawing, we drop the surface from the picture and just draw the marked graph if no confusion arise. So for example, the marked surface on figure \ref{markedsurface} may be just drawn as in figure \ref{markednosurface}.
\subsection{ Parameterization in the general case}
In this paper we will only consider $G$ cover of extended surface for which the base surface is connected and has genus zero. But this definition makes sense in the general case also.
\begin{definition}
Let $\tilde{\Sigma} \to \Sigma$ be a $G$ cover of an extended surface, $\Sigma$. A parameterization of $\tilde{\Sigma}$ is an isomorphism of this $G$-cover with one or gluing of several standard blocks. We have defined what we meant by gluing of $G$-cover (see section \ref{gluingofgcover} on page \pageref{gluingofgcover}). That is if $f$ is a parameterization of the $G$-cover $\ts$ then $f$ is just an isomorphism which looks like :\\
$$f: \ts  \longrightarrow S_{n_{1}}(g_{1}^{1},...g_{n_{1}}^{1} ;
h_{1}^{1},...h_{n_{1}}^{1}) \bigsqcup_{c_{i},h_{p}^{1},h_{q}^{2}} .....\bigsqcup_{c_{j},h_{t}^{k-1},h_{l}^{k}} S_{n_{k}}(g_{1}^{k},...g_{n_{k}}^{k} ;
h_{1}^{k},...h_{n_{k}}^{k})$$
\end{definition}
For the explanation of the above notation see the remark on section \ref{gluenotation} on page \pageref{gluenotation}.
\subsection{Visualizing parameterization}
\label{visualizingparameterization}
\begin{figure}
    \includegraphics{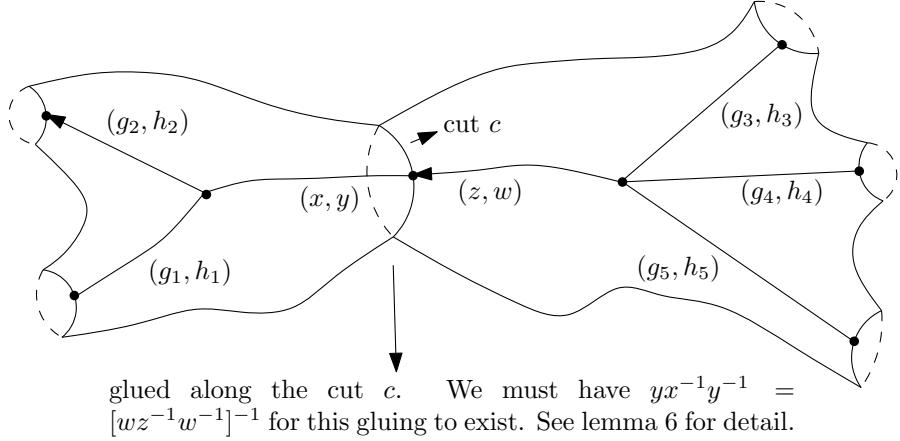}\\
    \caption{Visualization of the parameterization of $G$-cover with the gluing of $S_{3}(g_{1},g_{2},x;h_{1},h_{2},y)$ and $S_{4}(x^{-1},g_{3},g_{4},g_{5};y,h_{3},h_{4},h_{5})$} \label{parameterization1}
\end{figure}

\begin{figure}
    \includegraphics{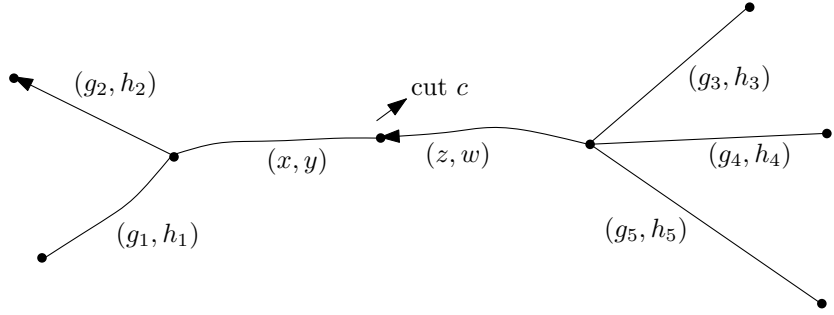}\\
    \caption{We usually do not draw surface and just draw the graph}\label{parameterization2}
\end{figure}

We need some kind of graphical way to represent parameterization of $G$-cover of extended
surface just as in the case of $|G| = 1$ which is described in the subsection \ref{gonecase} on page \pageref{gonecase}. So we need a similar kind of machinery in this general case. The lemma \ref{groupoidlemma} gives us a way to visualize $G$-covers of extended surface. Given a $G$-cover of an extended surface,$\tilde{\Sigma} \to \Sigma$, first we marked the base surface, $\Sigma$. In other word, we identify the base surface, $\Sigma$, with one or gluing of several standard spheres. For details about the marking of an extended surface see the subsection \ref{gonecase} on page \pageref{gonecase}. This marked base surface gives the parameterization of the base surface, $\Sigma$, with one or gluing of several $S_n$ where $n$ may vary. Let $C$ be the cut system of this marking. If we restrict the $G$ cover, $\tilde{\Sigma}$, to each connected component of $\Sigma \backslash C$, then the whole $G$ cover can be realized as a gluing of all this restriction. Each such restriction is isomorphic to one of the standard block, $S_{n}(g_{1},...,g_{n}; h_{1},...,h_{n})$, where $n$ may vary for each restriction. So the whole $G$ cover, $\tilde{\Sigma} \to \Sigma$, can be identify with the gluing of several standard blocks (one for each restriction). We can include all this data into the surface, $\Sigma$ as follows:\\
 We label each edge of our marking graph with a pair $(g_{i},h_{i})$ which comes from the identification of the restriction of $\tilde{\Sigma}$ to the component containing the edge. See the proof of the lemma in section \ref{groupoidlemma} to see how to assign $(g_{i},h_{i})$ to each edge. See a typical picture on figure \ref{parameterization1}.
sometimes, we drop the picture of the surface and just draw the graph
 for
simplicity, if no confusion arise. So for example, we will usually draw
 the
picture on figure \ref{parameterization1} as a simple figure in \ref{parameterization2}

This gives a visual presentation of the parameterization of our
 $G$-cover
of extended surface with one or gluing of several standard blocks.

\section{Moves}
At this point, we want to remind our readers about our main goal of this paper. we will eventually define a
$2$-dimensional CW complex $M(\ts, \Sigma)$, which has the set of
all parameterization of $\ts$ as the set of vertices.  The edges of $M(\ts, \Sigma)$
will be directed; we call them moves.  It is convenient to
look at $M(\ts, \Sigma)$ as a groupoid with objects---all vertices and
morphisms between two vertices---the set of homotopy classes of paths
on the edges of $M(\ts, \Sigma)$ from the first vertex to the second one
(going along an edge in the direction opposite to its orientation is
allowed). We will use group notation writing a path composed of edges
$E_1, E_2, \dots$ as a product $E_1 \dots E_2$, and we will write
$E^{-1}$ if the edge $E$ is traveled in the opposite direction. Then
the $2$-cells are interpreted as relations among the moves: we will
write $E_1\cdots E_k=id$ if the closed loop formed by the  edges
$E_1,\dots, E_k$ is contractible in $M(\ts, \Sigma)$; if we want to specify
the base point for the loop, we will write $E_1\cdots E_k(M)=id(M)$.
We will write $E:M \to M'$ if the edge $E$ goes from $M$ to $M'$.

Our Main Theorems state that the complex $M(\ts, \Sigma)$ is connected and simply-connected. This main theorem will be described in detail after we describe the moves or edges in this section.

\subsection{Standard Morphism} Before describing our moves, we need to
define some standard morphisms between our standard blocks. These are
described in terms of lifting of certain morphisms between base surfaces
$S_n$. These are:\\
$$z: S_n \longrightarrow S_n$$
$$b: S_3 \longrightarrow S_3$$
$$\alpha_{k,l}: S_{k+1} \sqcup S_{l+1} \longrightarrow S_{k+l}$$
For a more elaborate description of these morphism, see the paper \cite{BK}. Here we will just give a quick description.
\begin{description}
\item[1]{$z$ is the rotation of the sphere which cyclically permutes the boundary circles. That is if $m$ is a marking on $S_{n}$ then $z(m)$ will be the same marking (of course up to homotopy) on $S_{n}$ but with a different distinguished edge. See the figure \ref{zmap} where $n = 4$.
    \begin{figure}
    \includegraphics{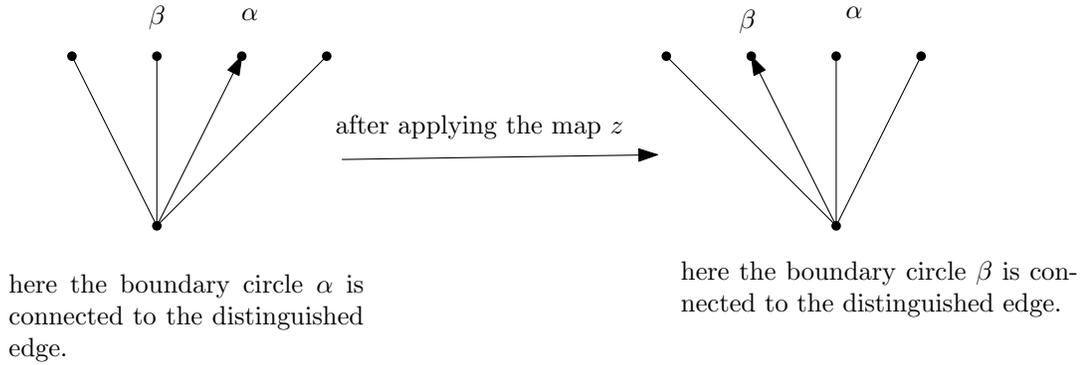}
    \caption{A pictorial description of $z$ map for $n = 4$}\label{zmap}
    \end{figure}}
    \item[2]{Now we will describe briefly the braiding, $b$, see \cite{BK} for detail. If we label the boundary circles of $S_{3}$ by $\alpha, \beta, \gamma$, then sometimes we will denote by $b_{\alpha ,\beta }$ the braiding of the $\alpha$ and $\beta$ component of $S_{3}$. Let $m$ is the graph on $S_{3}$ shown on the left hand side of figure \ref{b}. Then we define the $b_{\alpha, \beta}$ by figure \ref{b}.
        \begin{figure}
        \includegraphics{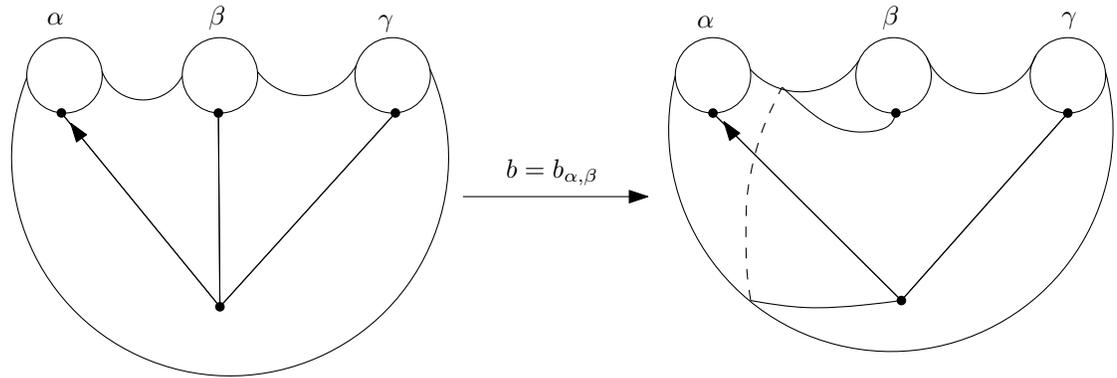}
        \caption{Visual description of $b$ map}\label{b}
        \end{figure}
        }
\item[3]{$\alpha_{k,l}$ is the identification of the result of gluing $S_{k+1}$ and $S_{l+1}$ (along the $(k+1)$ th boundary component of the first one with the $1$ st boundary component of the last one) with, $S_{k+l}$, the standard sphere with
$(k+l)$ hole. For more detail description see \cite{BK}.}
\end{description}
\begin{lemma} \label{liftingstandardmor}
    Each of the above defined morphisms between $S_n$ ($n = 3$ in the case of $b$) can be uniquely lifted
    to the $G$-cover of $S_n$, $S_n(g_1,....g_n;h_1,...h_n)$, described
    by:\\
    $$\tilde{z}: S_{n}(g_1...g_n;h_1...h_n) \longrightarrow
    S_{n}(g_n,g_1...g_{n-1};h_n,h_1...h_{n-1})$$
    $$\tilde{b}: S_{3}(g_{1},g_{2},g_{3};h_{1},h_{2},h_{3}) \longrightarrow S_{3}(g_{1}g_{2}g_{1}^{-1},g_{1},g_{3};h_{2}g_{1}^{-1},h_{1},h_{3})$$
 $$\tilde{\alpha_{k,l}}: S_{k+1}(\mathbf{g},\mathbf{h}) \sqcup
    S_{l+1}(\mathbf{g'},\mathbf{h'})\longrightarrow
    S_{k+l}(\mathbf{g''},\mathbf{h''})$$
where \\
\begin{eqnarray*}
\mathbf{g} & = & (g_1...g_{k+1})\\
\mathbf{h} & = &  (h_1...h_{k+1})\\
\mathbf{g'} & = &  (g_1',...g_{l+1}')\\
\mathbf{h'} & = &  (h_1',...h_{l+1}')\\
\end{eqnarray*}

So that $h_{k+1}g_{k+1}^{-1}h_{k+1}^{-1} =
[h_{1}'g_{1}'^{-1}h_{1}'^{-1}]^{-1}; h_{k+1} = h_1'$\\
And\\

\begin{eqnarray*}
    \mathbf{g''} &=& (g_1...g_k,g_2'...g_{l+1}')\\
    \mathbf{h''} &=& (h_1...h_k,h_2'...h_{l+1}')\\
    \end{eqnarray*}

\end{lemma}

\begin{remark} Notice that we need
    $h_{k+1}g_{k+1}^{-1}h_{k+1}^{-1} =
 [h_{1}'g_{1}'^{-1}h_{1}'^{-1}]^{-1}$
 so that it is possible to glue the $(k+1)$ th
    boundary circle of $S_{k+1}$ with $1$st boundary circle of
 $S_{l+1}$. This is to make sure that the monodromy (see
    \ref{gluingofgcover} on page \pageref{gluingofgcover}) of $(k+1)$
 th
    boundary circle of $S_{k+1}$ is inverse of the monodromy of $1$st
    boundary circle of $S_{l+1}$ (recall that the orientation of two
    boundary circles must be opposite for glue to exist). The second
    equality, $h_{k+1} = h_1'$, is an extra one which is not needed for the gluing to exist but we assume this extra condition whenever we want to apply $\tilde{\alpha_{k,l}}$.
\end{remark}
\begin{proof}
    First consider the $G$-covers,$ S_{n}(g_1...g_n;h_1...h_n)$ and
    $S_{n}(g_n,g_1...g_{n-1};h_n,h_1...h_{n-1})$, and the map $z: S_n
    \longrightarrow S_n$. Our base surface, $S_{n}$ has genus zero and
 the
    map $z$ preserves the monodromy which can be easily checked by hand.
 So
    according to lemma 2 (see \ref{liftinglemme} on page
    \pageref{liftinglemme}), this map $z$ can be uniquely lifted to the
    $G$-cover. Similarly, in all the other cases, all we have to check
 that
    the maps $b: S_3 \longrightarrow S_3$ and
    $\alpha_{k,l}: S_{k+1} \sqcup S_{l+1} \longrightarrow S_{k+l}$
    preserve the monodromy. This can be easily checked by hand, using the lemma 6.
\end{proof}

\subsubsection[$\phi _{x}$ morphism]{The standard morphism, $\phi _{x}, x \in G$}
This is just a reminder to the reader. We have already defined $\phi_{x}$ move before. This is the morphism between standard blocks defined by\\
$\phi_x$ : $S_n(g_1...g_n;h_1...h_n) \longrightarrow
S_n(xg_1x^{-1}...xg_nx^{-1};h_1x^{-1}...h_nx^{-1})$. See lemma \ref{isolemma} on page \pageref{isolemma}.

\subsection{$\mathbf{Z}$, $\mathbf{B}$ ,and $\mathbf{F}$ Move}
For $|G|=1$, these moves are described in the paper \cite{BK}. Reader
 might
want to read this first before continue. We again want to remind our
readers that our base surface will always have genus 0, unless
 otherwise
specified. We will first have three moves similar to the $Z$, $B$ and
 $F$
moves of the paper \cite{BK}. These moves will be called
 $\mathbf{Z}$-move,
$\mathbf{B}$-move and $\mathbf{F}$-move respectively. Note that we
 use
boldface letter to denote these three moves so that it does not get
confused with the $Z$, $B$, and $F$ moves of the paper \cite{BK}.
 Then we
will have two extra moves which do not have any correspondence to the
 paper
\cite{BK}.
So again $Z$, $B$ and $F$ will denote the moves in the case of
 $|G|=1$
and
$\mathbf{Z}$, $\mathbf{B}$ ,and $\mathbf{F}$ will denote the moves in
 the
general case. Each move will take a parameterization to another
parameterization. So let us start describing these moves in more
 detail:\\
\subsubsection{\emph{\textbf{$\mathbf{Z}$ move}}}
Given a $G$-cover $(\tilde{\pi }: \ts  \longrightarrow \Sigma,
\{\tilde{p_a}\})$ of $(\Sigma,\{p_a\})$, and a parameterization of this
$G$-cover with a standard block. This means we have an isomorphism,
 $f$, of
this $G$-cover with one of the standard block say
$S_n(g_1...g_n;h_1...h_n)$. $\mathbf{Z}$ move takes this
 parameterization
$f$ to the parameterization $\tilde{z} \circ f$ where $\tilde{z}$ is
 the
standard morphism defined in the lemma 5. This new parameterization
identify the original $G$-cover with
$S_n(g_n,g_1...g_{n-1};h_n,h_1...h_{n-1})$. Look at the diagram below\\
\begin{displaymath}
\underbrace{\xymatrix{\mbox{Original G-cover} \ar[r]^{f} &
S_n(g_1...g_n;h_1...h_n) \ar[r]^{\tilde{z}} &
S_n(g_n,g_1...g_{n-1};h_n,h_1...h_{n-1})}}_{\tilde{z} \circ f =
\mathbf{Z}(f), \mbox{ $\mathbf{Z}$ move applied to f}}
\end{displaymath}
\subsubsection{\emph{\textbf{$\mathbf{B}$ move}}}
Given a $G$-cover $(\tilde{\pi }: \ts  \longrightarrow \Sigma,
\{\tilde{p_a}\})$ of $(\Sigma,\{p_a\})$, and a parameterization of this
$G$-cover with a standard block, $S_{3}(g_{1},g_{2},g_{3};h_{1},h_{2},h_{3})$. This means we have an isomorphism,$f$,
 of this $G$-cover with $S_{3}(g_{1},g_{2},g_{3};h_{1},h_{2},h_{3})$. $\mathbf{B}$ move takes this
 parameterization $f$ to the parameterization $\tilde{b} \circ f$ where $\tilde{b}$
 is the standard morphism defined in the lemma \ref{liftingstandardmor}. This new parameterization
identify the original $G$-cover with $S_{3}(g_{1}g_{2}g_{1}^{-1},g_{1},g_{3};h_{2}g_{1}^{-1},h_{1},h_{3})$.
\subsubsection{\emph{\textbf{$\mathbf{F}$ move}}}
\begin{figure}
    \includegraphics{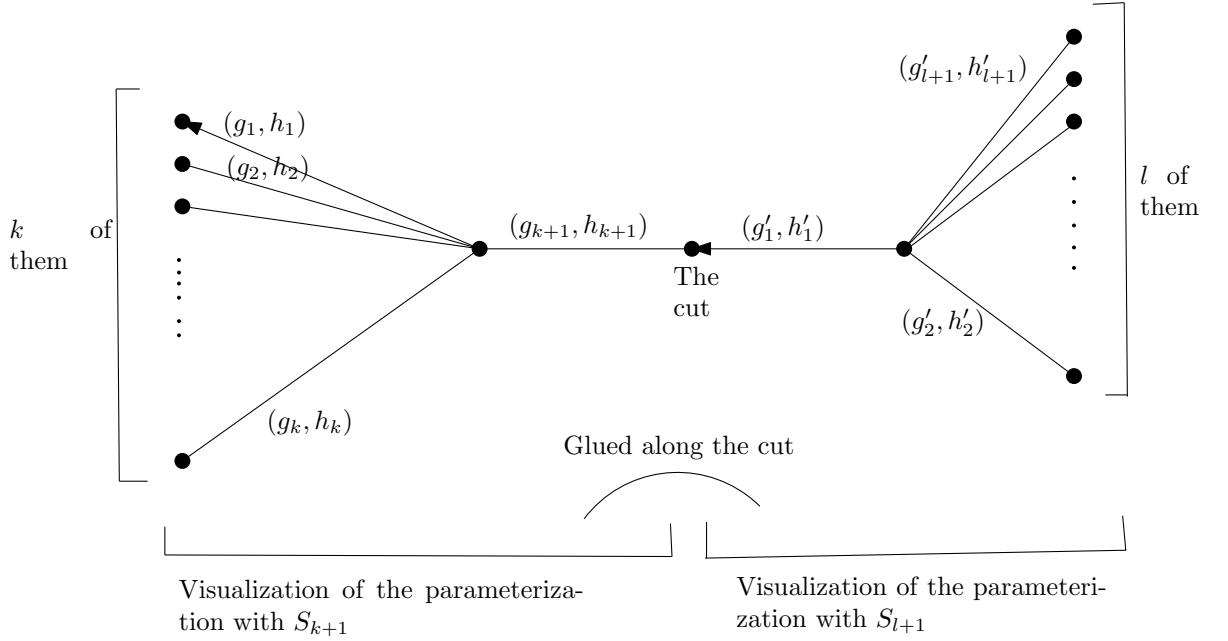}\\
    \caption{Visual description of the parameterization of $f$}\label{fmove}
\end{figure}

\begin{figure}
    \includegraphics{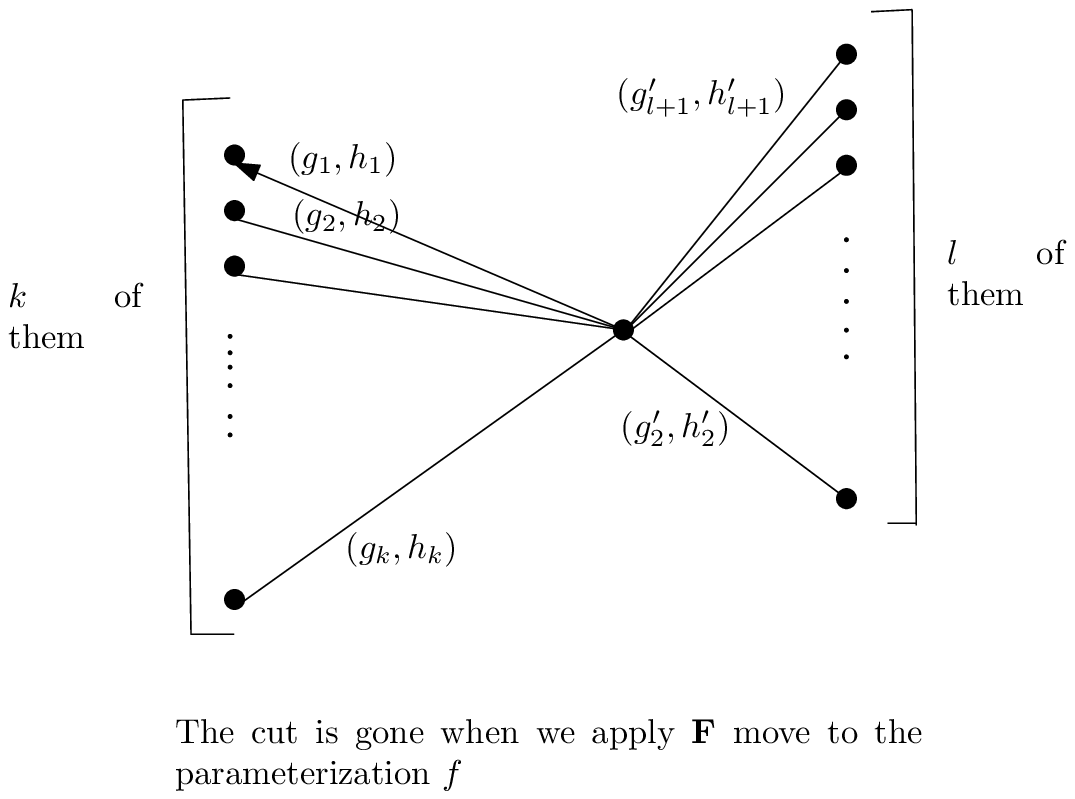}\\
    \caption{Visualization of the parameterization $\mathbf{F}(f)$}\label{fmove2}
\end{figure}

Here the setup is a little different from the above two moves. Fist,
 let
$\Sigma$ be an extended surface of genus $0$ with one cut $\{c\}$. This
cut \{$c$\} divides $\Sigma$ into two pieces. Say, $G$-cover of one
 piece is
parameterized with $S_{k+1}(\mathbf{g},\mathbf{h})$ and the other with
$S_{l+1}(\mathbf{g^{'}},\mathbf{h^{'}})$. That is the parameterization,
$f$, is given by:\\
$$f: \ts  \longrightarrow S_{k+1}(\mathbf{g},\mathbf{h})
 \bigsqcup_{c,h_{k+1}} S_{l+1}(\mathbf{g^{'}},\mathbf{h^{'}})$$
where\\
\begin{eqnarray*}
\mathbf{g}&=&(g_1...g_{k+1})\\
\mathbf{h}& =& (h_1...h_{k+1})\\
\mathbf{g'}& =& (g_1',...g_{l+1}')\\
\mathbf{h'} &= &(h_1',...h_{l+1}')
\end{eqnarray*}
And also $h_{k+1}g_{k+1}^{-1}h_{k+1}^{-1} =
[h_1'g_{1}'^{-1}h_{1}'^{-1}]^{-1} ; h_{k+1} = h_1'$\\
Thus the parameterization of the whole $G$-cover, $f$, identifies the
$G$-cover with the result of gluing these two standard blocks along the
$k+1$ st boundary circle of one side with the $1$st boundary circle of
 the other. See the picture on figure \ref{fmove} for more detail.
Now by applying $\mathbf{F}$ move to $f$, we get a new parameterization
$\tilde{\alpha_{k.l}} \circ f$. Note that, $\tilde{\alpha_{k.l}} \circ
 f$
is a parameterization of G-cover of $\Sigma$ with no cuts with
$S_{k+l}(\mathbf{g''},\mathbf{h''})$ where\\
\begin{eqnarray*}
    \mathbf{g''}& =& (g_1...g_k,g_2'...g_{l+1}')\\
    \mathbf{h''} &=& (h_1...h_k,h_2'...h_{l+1}').
\end{eqnarray*}
See the figure \ref{fmove2}.

\begin{remark}
If two components are glued along a cut where the edges corresponding to this cut are labeled by $(x,y)$ and $(z,w)$, then for this gluing to exist, we must have $$xy^{-1}x^{-1} = [zw^{-1}z^{-1}]^{-1}$$
See subsection \ref{gluenotation} on page \pageref{gluenotation} for more detail explanation. Now to apply $\mathbf{F}$ move to this cut, we also required the extra condition that $y=w$. This will imply that $xz = 1$. The requirement, $y=w$, is not necessary for gluing to exist but we only require this whenever we want to apply $\mathbf{F}$ move to a cut.
\end{remark}

\begin{figure}
\includegraphics{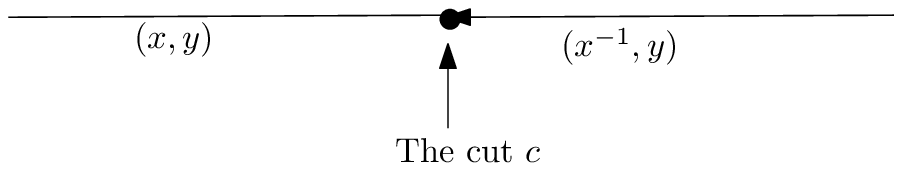}
\caption{}\label{fmove3}
\end{figure}

\textbf{A remark about notation}: Note that by applying a $\mathbf{F}$
move, we are removing a cut. If the cut for one component is labeled by
$(x,y)$ then the cut for the other component must be labeled by
$(x^{-1},y)$ for $\mathbf{F}$ move to apply. See the picture on figure \ref{fmove3}.
This $x \in G$ is completely determined by the label of all the other
 cuts
since they all must multiplied to $1 \in G$. But for the second
 component
$y$, we have a choice. So it is better to encode the cut $c$ and the
 label
$y$ in the notation. From now on, a $\mathbf{F}$ move will always be
denoted by $\mathbf{F_{c,y}}$ to emphasize the cut, $c$, and the choice
 of
point, $y$.\\
\subsection{$\mathbf{P_x}$-move where $x \in G$}\label{pxmove}
Beside $\mathbf{Z}$, $\mathbf{B}$ and $\mathbf{F}$ move, we need one
 more
move which we will denote by $\mathbf{P_x}$ where $x \in G$. This
$\mathbf{P_x}$-move is defined in the following way:\\
Given a $G$-cover $(\tilde{\pi }: \ts  \longrightarrow \Sigma,
\{\tilde{p_a}\})$ of $(\Sigma,\{p_a\})$, and a parameterization of this
$G$-cover with a standard block. This means we have an isomorphism,
 $f$, of
this $G$-cover with one of the standard block say
$S_n(g_1...g_n;h_1...h_n)$. $\mathbf{P_x}$ move takes this
 parameterization
$f$ to the parameterization $\phi_x \circ f$ where $\phi_x$ is the
 morphism
defined by \\
$\phi_x$ : $S_n(g_1...g_n;h_1...h_n) \longrightarrow
S_n(xg_1x^{-1}...xg_nx^{-1};h_1x^{-1}...h_nx^{-1})$. See the lemma \ref{isolemma} for the description of $\phi (x)$. See the diagram below:\\

$$\ts \stackrel{f}{\rightarrow} S_{n}(g_{1}....g_{n};h_{1}...h_{n}) \stackrel{\phi_{x}}{\rightarrow}
S_{n}(xg_1x^{-1}...xg_nx^{-1};h_1x^{-1}...h_nx^{-1})$$
$$\phi_{x} \circ f = \mathbf{P}_{x}(f)$$

\subsection{$\mathbf{T_{z,y}^{c}}$ move}
\begin{figure}
\includegraphics{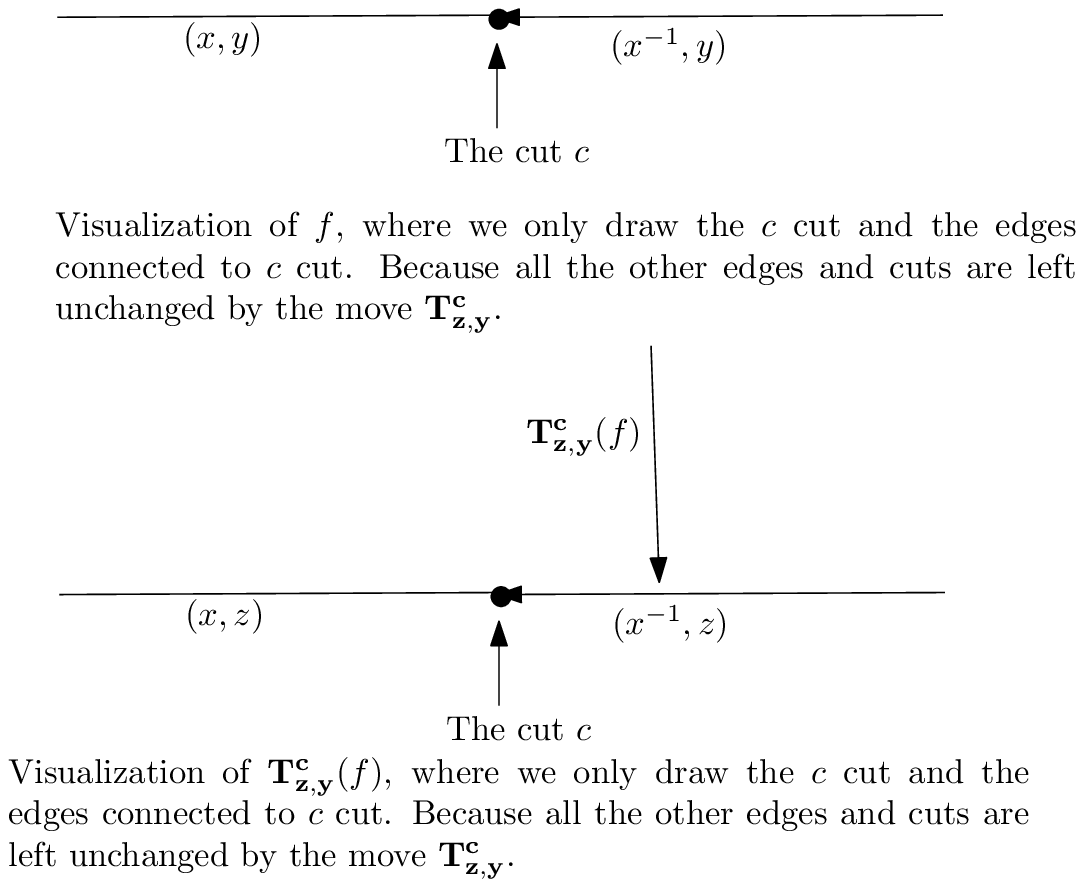}
\caption{The sequence of diagram describing $\mathbf{T}$ move}\label{tmove}
\end{figure}

Say we have the following parameterization, $f$, of our $G$-cover:\\
$$f: \ts  \longrightarrow S_{n}(g_{1}...g_{n-1},x ; h_{1}...h_{n-1},y)
\bigsqcup_{c,y} S_{m}(x^{-1},u_{1}...u_{m-1} ; y,v_{1}...v_{m-1}) $$
$c$ is the cut where we glue the two components. Note that all
 conditions
are satisfied so that the gluing make sense. Here $x$ is determine by
 all
the other element, namely, $x = (g_{1}...g_{n-1})^{-1} =
(u_{1}...u_{m-1})^{-1}$, but $y$ is not determine; we can choose $y$
freely. By definition, $\mathbf{T_{z,y}^{c_{i}}}(f)$ is the following
parameterization:\\
$$\mathbf{T_{z,y}^{c}}(f): \ts  \rightarrow S_{n}(g_{1}..g_{n-1},x ;
h_{1}..h_{n-1},z) \bigsqcup_{c,y} S_{m}(x^{-1},u_{1}..u_{m-1} ;
z,v_{1}..v_{m-1})$$
Where we replace the choice of $``y"$ with the choice of $``z"$.\\
\begin{remark}
Although we introduce this new move,
$\mathbf{T_{z,y}^{c}}$, this can in fact be thought of as composition of two
$\mathbf{F}$ move. Namely $\mathbf{T_{z,y}^{c}} =
\mathbf{F_{c,z}^{-1}F_{c,y}}$. The move, $\mathbf{F_{c,y}}$, will
 remove
the cut, $c$, with the choice of point $y$, while the move
$\mathbf{F_{c,z}^{-1}}$ will replace the cut, $c$, but this time with
 the
choice of point $z$. The reason for introducing such an extra move will
 be
clear later, but introducing this new move $\mathbf{T_{z,y}^{c}}$ and
adding the relation $\mathbf{T_{z,y}^{c}} =
 \mathbf{F_{c,z}^{-1}F_{c,y}}$
will not do any harm, since one can easily sees that the complex $M(\ts
 ,
\Sigma)$ without the $\mathbf{T_{z,y}^{c}}$ move and
 $\mathbf{T_{z,y}^{c}}
= \mathbf{F_{c,z}^{-1}F_{c,y}}$ relation is connected and
 simply-connected
iff the complex $M(\ts , \Sigma)$ with the $\mathbf{T_{z,y}^{c}}$ move
 and
$\mathbf{T_{z,y}^{c}} = \mathbf{F_{c,z}^{-1}F_{c,y}}$ relation is
 connected
and simply-connected.
\end{remark}
See the picture on figure \ref{tmove} for a visual description of $\mathbf{T}$ move.
\subsection{Moves, when more than one standard blocks are glued together}
\begin{figure}
 \includegraphics{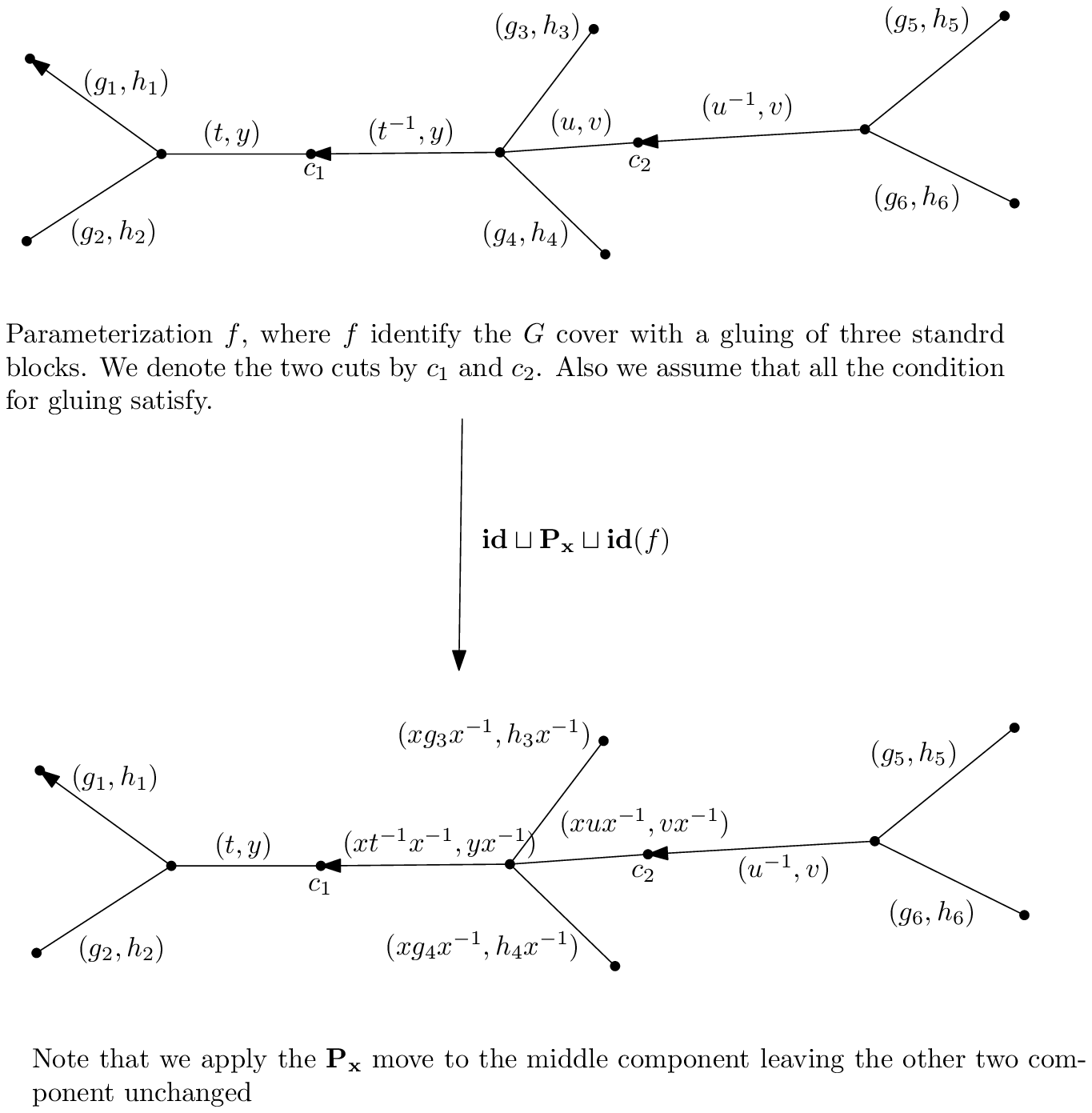}
 \caption{An example of a move when having more than one component.} \label{pmove}
 \end{figure}

Let $\ts$ be a $G$ cover of $\Sigma$. Also let $f$ be a parameterization of this $G$ cover with gluing of several standard blocks. In other word,
 $$f: \ts \rightarrow S_{n_1}(\mathbf{g_1,h_1}) \bigsqcup_{glued} S_{n_2}(\mathbf{g_2,h_2}) \bigsqcup_{glued}.....\bigsqcup_{glued} S_{n_k}(\mathbf{g_k,h_k})$$
 This parameterization of $f$, by the restriction on each component, can be realized as a gluing of $k$ parameterization. That is $f = f_{1} \sqcup f_{2} \sqcup .....\sqcup f_{k}$.
 Let $\mathbf{E}$ be one of the above five moves, that is $\mathbf{E} \in \mathbf{Z,B,F,P_{x},T}$. Then by the move $id \sqcup id \sqcup ....\sqcup \mathbf{E} \sqcup ...\sqcup id$($f$), where $\mathbf{E}$ appears in the $i$th component, we mean that we only apply $\mathbf{E}$ move to the $f_{i}$ parameterization and identity to all others. An example of such a move is shown on figure \ref{pmove}. For simplicity of the picture, we take three components.
\subsection{$\mathbf{B_{i}}$ move}\label{bimove1}
\begin{figure}
\includegraphics{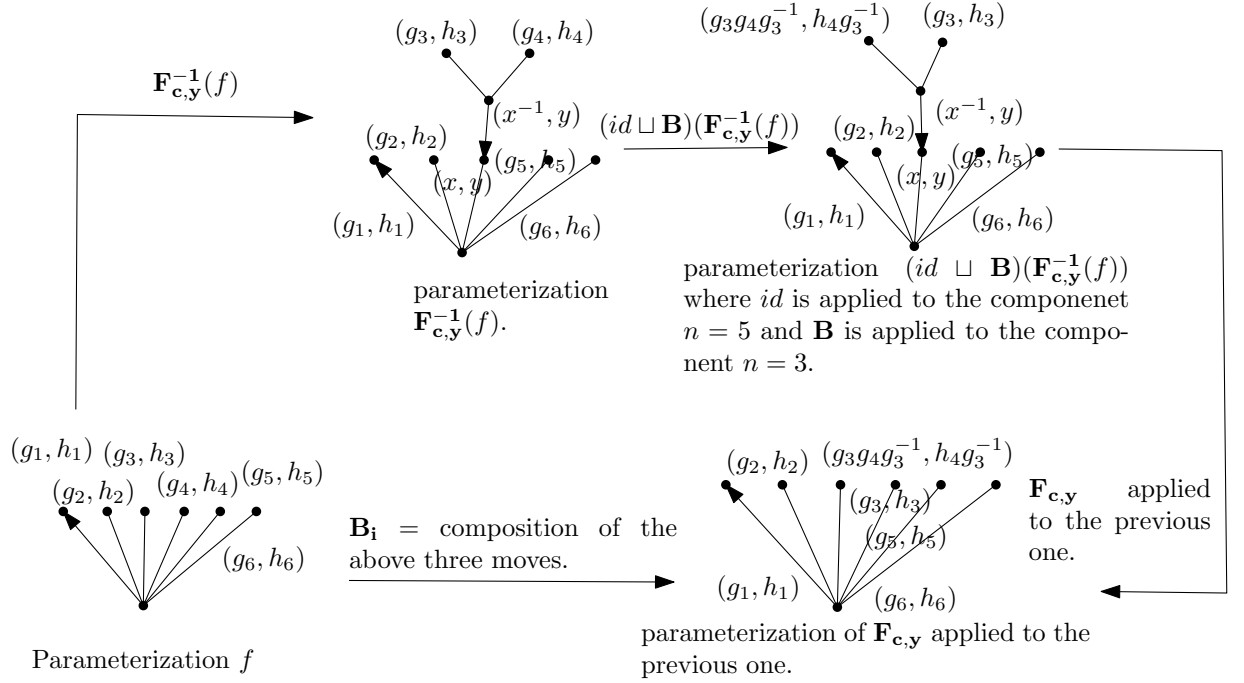}
\caption{description of $\mathbf{B_{i}}$ move as a composition of previously defined move}\label{bimove}
\end{figure}

The $\mathbf{B}$ move is very restricted; it is only defined for $n = 3$. So we need a little flexibility. This $\mathbf{B_{i}}$ move will be a braiding of $i$ and $i+1$ boundary circles. This is in fact not a new move but a composition of previously defined moves. This will be used many times from now on when we will describe the other moves and relation. So it is important to do this right now. Let $f$ be a parameterization of a $G$ cover, $\tilde{\Sigma}$, with a standard block say $S_{n}(g_{1},...,g_{n};h_{1},...,h_{n})$. Then $\mathbf{B_{i}}$ is the composition of the moves described in the figure \ref{bimove}. Here we only draw pictures for $n=6$ for convenience but readers are clear of what should be done for other $n$.
If $$\ts \stackrel{f}{\rightarrow} S_{n}(g_{1},...,g_{n};h_{1},...,h_{n})$$
Then $$\ts \stackrel{\mathbf{B_{i}}(f)}{\rightarrow} S_{n}(g_{1},..,g_{i-1},g_{i}g_{i+1}g_{i}^{-1},g_{i},...,g_{n};h_{1},..,h_{i-1},h_{i+1}g_{i}^{-1},h_{i},..,h_{n})$$
\section{The Complex}
\subsection{Definition of the complex $M(\ts ,\Sigma)$}
We are given a fixed $G$-cover, $\Pi: \ts  \longrightarrow \Sigma$,
 where
$\Sigma$ will always denote an extended surface which is compact,
orientable, closed surface of genus zero unless otherwise specified. We
will first define the $1$-skeleton of the complex, $M(\ts ,\Sigma)$,
 that
is the vertex and the edges.\\
\{The set of vertex of $M(\ts ,\Sigma)$\} = \{The set of parameterization of the $G$-cover $\Pi: \ts  \longrightarrow
\Sigma$\}\\
A vertex, $\alpha$, is connected to vertex $\beta$, directed from
 $\alpha$
to $\beta$, if $\beta$ can be obtained from $\alpha$ by applying one of
 the
following moves:\\
$\mathbf{Z,Z^{-1},B_i,B_{i}^{-1},F,F^{-1}, P_x,P_x^{-1},
\mathbf{T_{z,y}^{c}},(\mathbf{T_{z,y}^{c}})^{-1}}$ in one of the components.
Having defined the vertex and edges of the complex $M(\ts ,\Sigma)$ ,
 now
we need to define the relations or 2-cells of $M(\ts ,\Sigma)$. We
 define
the relations in the next section.
\subsection{The Statement Of the Main Theorem}
\begin{maintheorem}
The complex, $M(\ts ,\Sigma)$, with the above defined edges and
 relations
(which will be described in the next section), is connected and simply
connected.
\end{maintheorem}

\section{Relations} Now it is time to define all the relations between
these moves. These relation will be described in the following
subsections:\\
\subsection{Obvious Relation}
\begin{center}
$\mathbf{EE^{-1}} = \mathbf{E^{-1}E} = 1$
\end{center}
Where $\mathbf{E}$ is one of the five moves, that is $\mathbf{E \in
\{Z,B_{i},F,P_{x},T_{z,y}^{c}\}}$
\subsection{$\mathbf{P_{x}}$ Relation}
\begin{figure}
\includegraphics{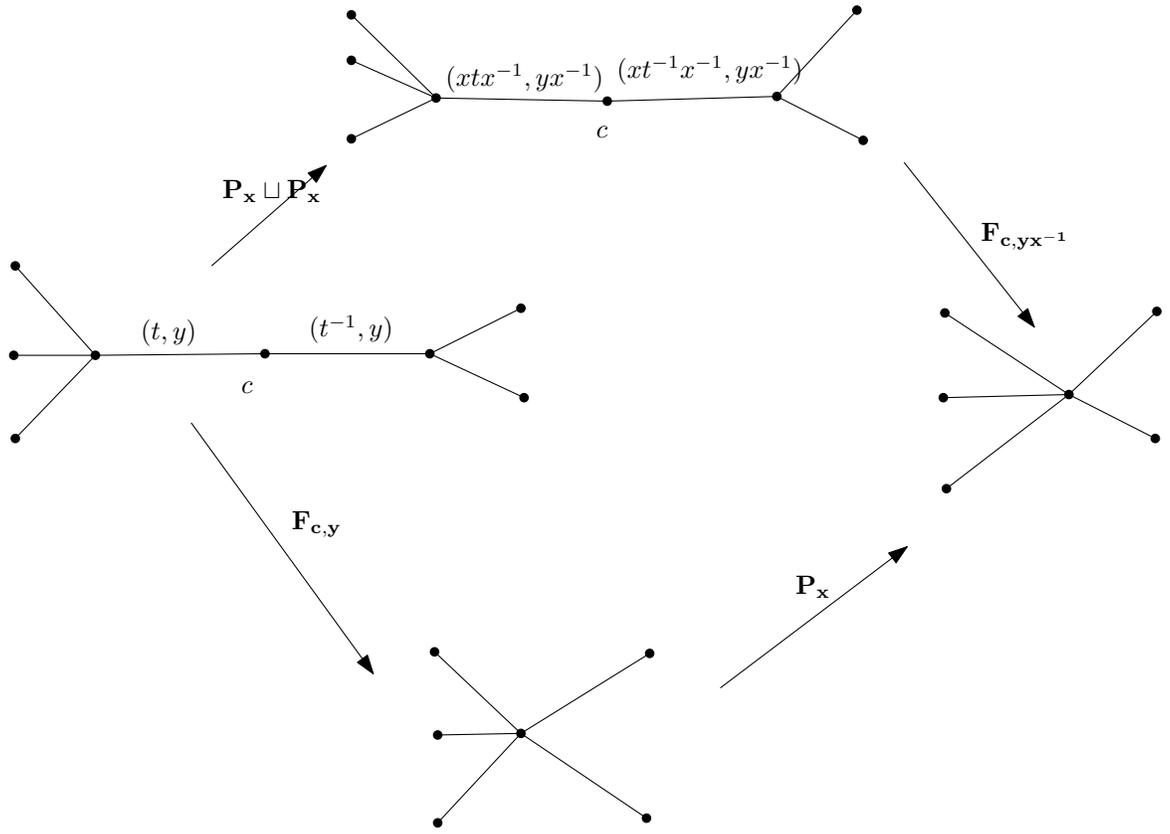}
\caption{Visual description of the third $\mathbf{P_{x}}$ relation}\label{indforpx}
\end{figure}

$\mathbf{P_x}$ move commutes with all the other moves. More precisely
 for
all $x,y \in G$ and all cuts, $c$, we have the following relations:\\
\begin{enumerate}
\item {$\mathbf{P_xZ} = \mathbf{ZP_x}$}
\item{$\mathbf{P_xB_i} = \mathbf{B_iP_x}$}
\item{$\mathbf{P_xF_{c,y}} = \mathbf{F_{c,yx^{-1}}(P_x \sqcup P_{x})}$. See figure \ref{indforpx}.}
\item{$\mathbf{P_xP_y} = \mathbf{P_{xy}}$}
\end{enumerate}
\begin{remark}
Note the change of indices for $\mathbf{F}$. Why we need to change the indices is clear from figure \ref{indforpx}. Here we only label the edges corresponding to the cut $c$.
\end{remark}
\begin{remark}
The last relation shows that whatever relation satisfied by the elements of $G$, the same relation hold for $\mathbf{P_x}$. So in particular this implies $\mathbf{(P_x)^{-1}} = \mathbf{P_{x^{-1}}}$.
\end{remark}
\subsection{$\mathbf{P-F}$ relation}
\begin{figure}
\includegraphics{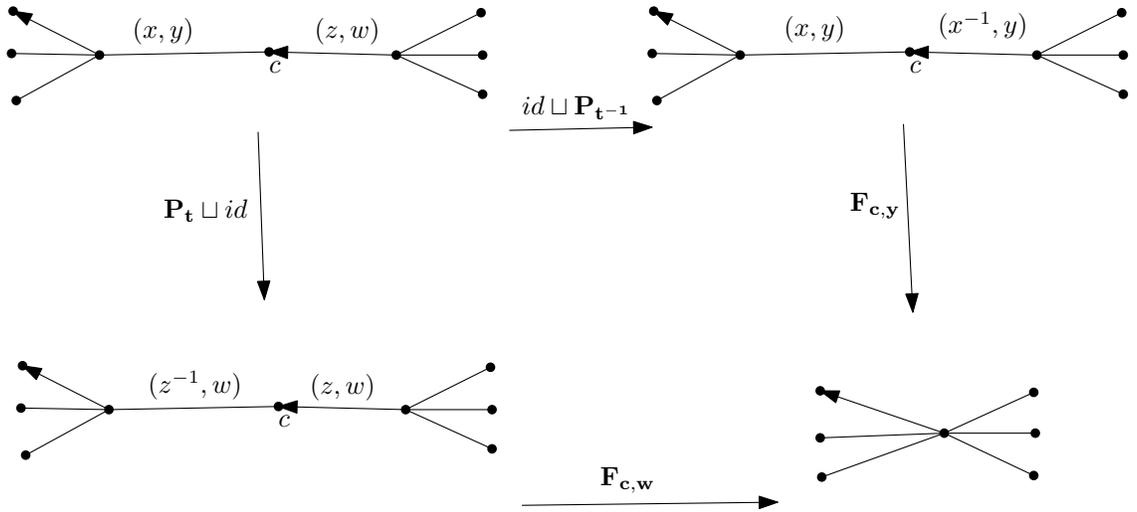}
\caption{Visual description of $\mathbf{P-F}$ relation}\label{pfrelation}
\end{figure}

Let $f$ be a parameterization having two component glued to a cut $c$ where one edge of the cut labeled by $(x,y)$ and the other edge of the cut labeled by $(z,w)$. For this glue to exist, we must have $yx^{-1}y^{-1} = [wz^{-1}w^{-1}]^{-1}$, see the subsection \ref{gluenotation} on page \pageref{gluenotation}. Let $t = w^{-1}y$. Then one can easily check that $txt^{-1} = z^{-1}$, $yt^{-1} = w$, $t^{-1}zt = x^{-1}$ and $wt = y$. Then $\mathbf{P-F}$ relation is the following relation://
$$\mathbf{F_{c,w}(P_{t} \sqcup id)} = \mathbf{F_{c,y}(id \sqcup P_{t^{-1}})}$$
See the figure \ref{pfrelation} for a visual description of this relation.
\subsection{$\mathbf{Z}$ Relation}
\begin{figure}
\includegraphics{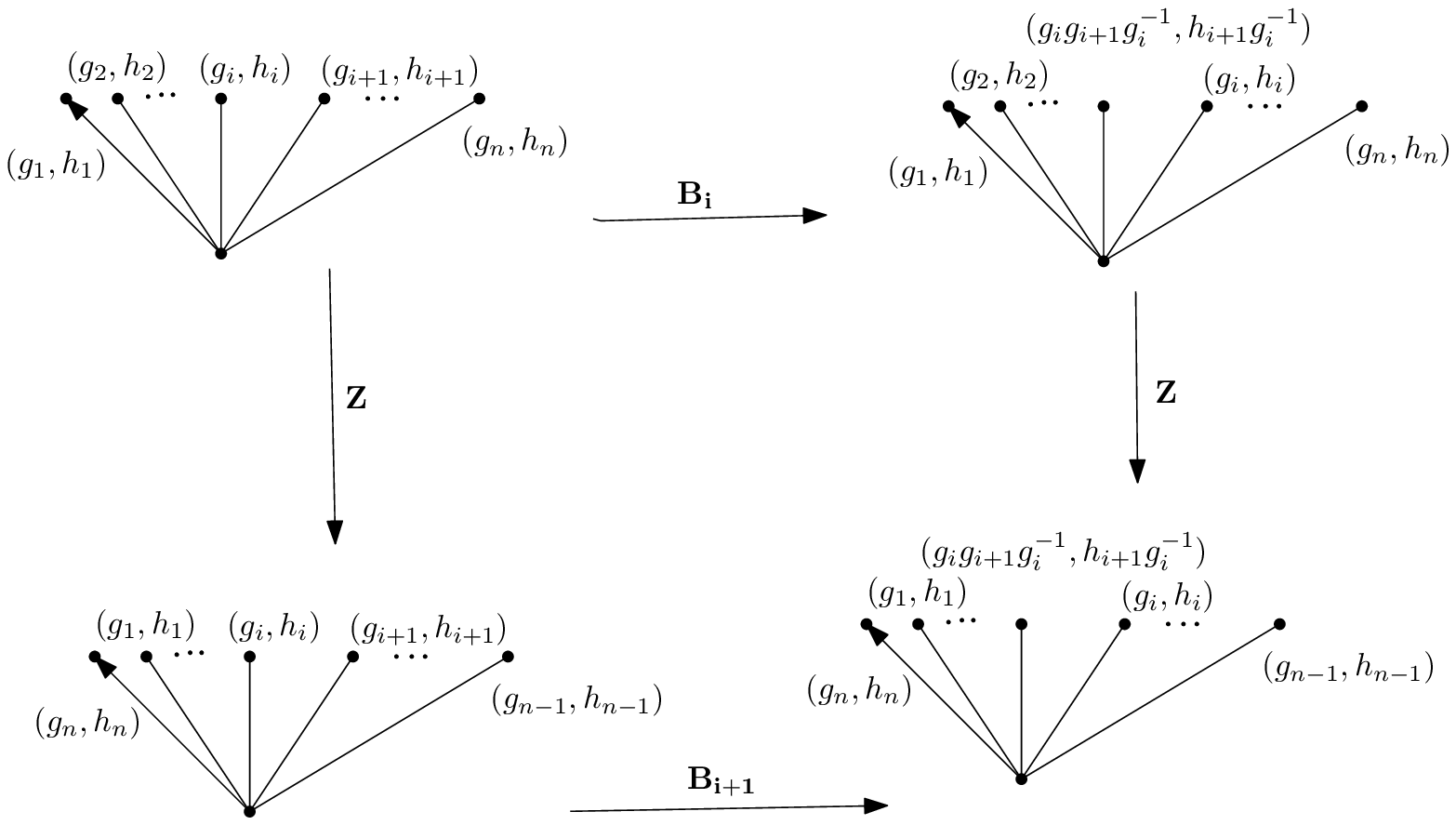}
\caption{Notice the change of indices for $\mathbf{B}$.} \label{bzrelation}
\end{figure}

$\mathbf{Z}$ commutes with all the other moves. More precisely, let a parameterization, $f$, identifies a $G$ cover,$\ts$, with the standard block $S_{n}(g_{1},...,g_{n};h_{1},...,h_{n})$. Then we have the following relation:\\
$$\mathbf{ZB_i} = \mathbf{B_{i+1}Z}$$
See figure \ref{bzrelation} for a visual description of this relation.\\

\begin{figure}
\includegraphics{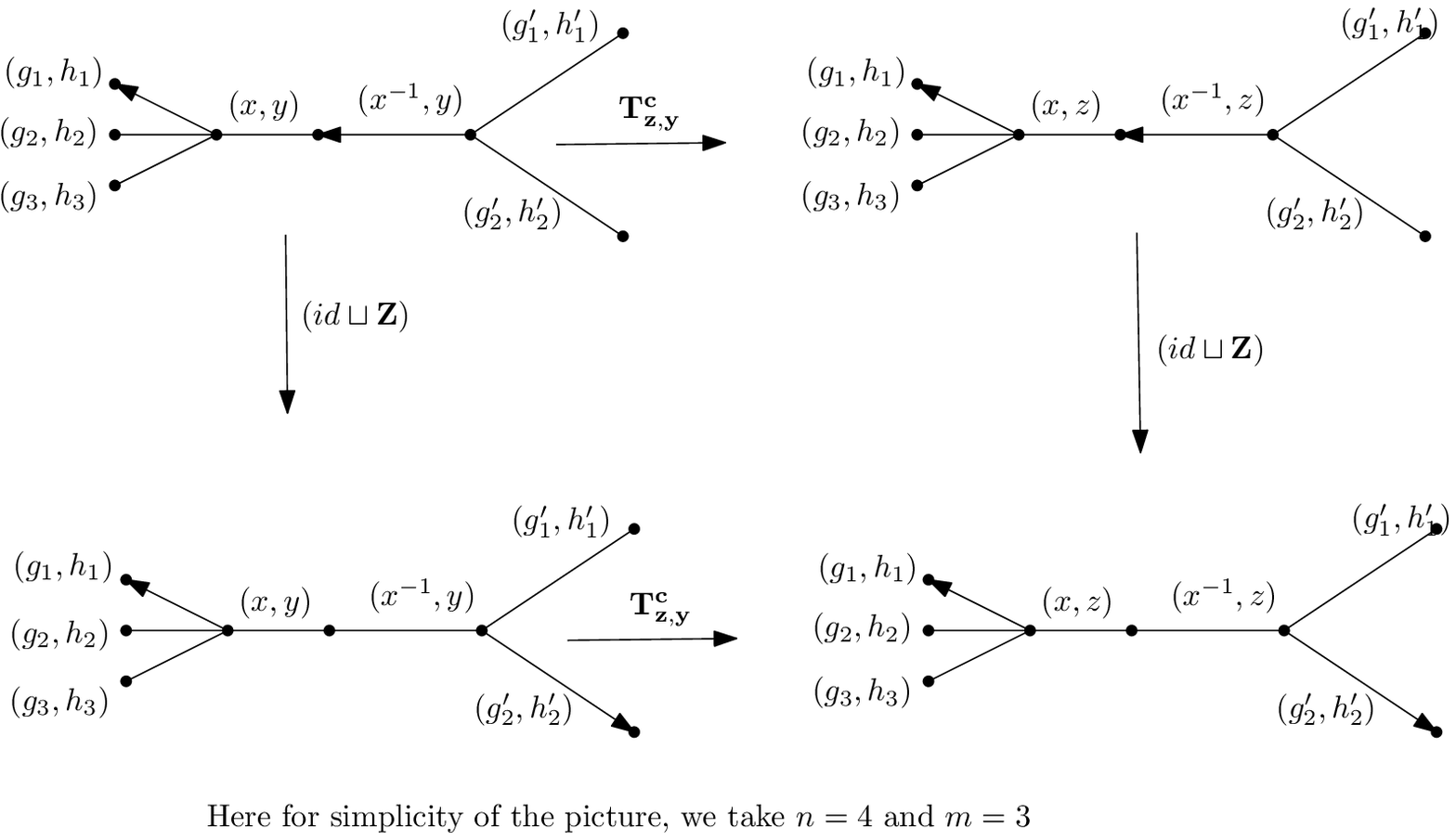}
\caption{The commutativity of $\mathbf{Z}$ and $\mathbf{T}$ moves}\label{ztrelation}
\end{figure}

Also let a parameterization, $f$, identifies a $G$ cover, $\ts$, with gluing of two standard blocks along a cut $c$.
 $$f: \ts \rightarrow S_{n}(g_{1},...,g_{n-1},x;h_{1},...,h_{n-1},y) \bigsqcup_{c,y} S_{m}(x^{-1},g_{1}',...,g_{m-1}';y,h_{1}',...,h_{m-1}')$$
 Then we have the following relation:\\
 $$(id \sqcup \mathbf{Z})\mathbf{T^{c}_{z,y}} = \mathbf{T^{c}_{z,y}}(id \sqcup \mathbf{Z})$$
See figure \ref{ztrelation} for a visual description of this relation.

\subsection{$\mathbf{B}$ Relation}
\begin{figure}
\includegraphics{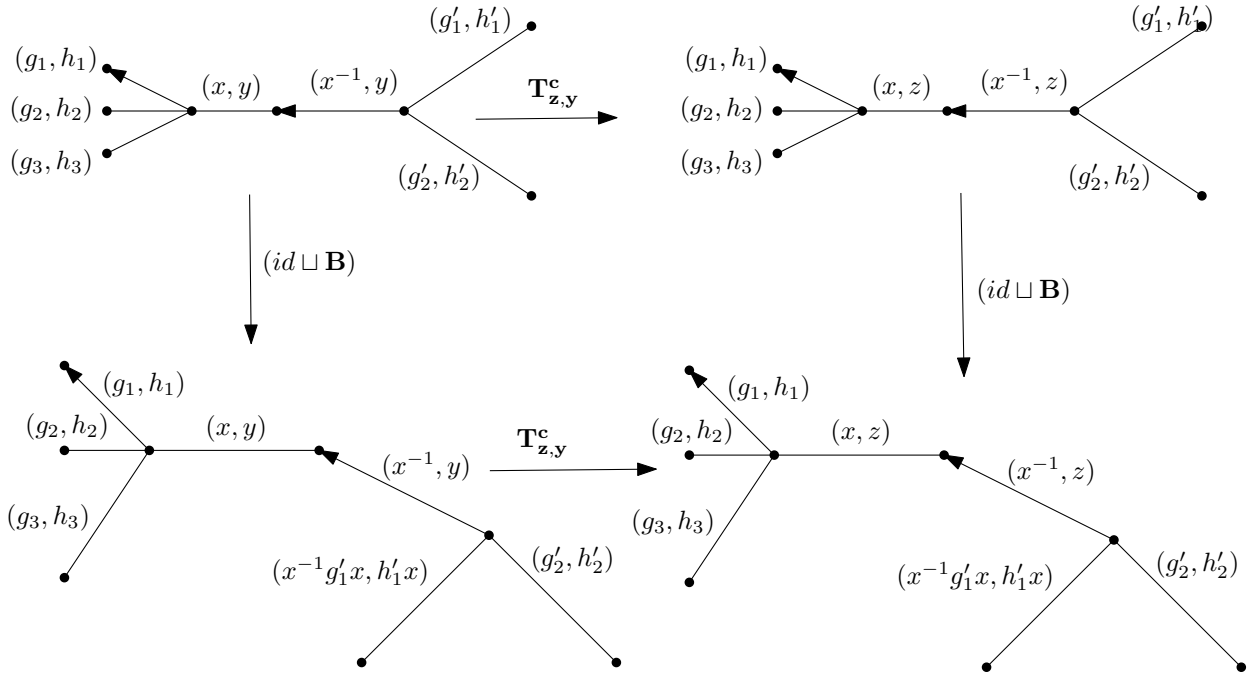}
\caption{The commutativity of $\mathbf{B}$ and $\mathbf{T}$ moves}\label{btrelation}
\end{figure}

From above relations, we already know that $\mathbf{B}$ move commutes
 with
$\mathbf{Z}$ and $\mathbf{P_{x}}$ move. But we have more--$\mathbf{B}$
 also
commutes with the $\mathbf{T}$ move. More precisely, let a parameterization, $f$, identifies a $G$ cover,$\ts$, with gluing of two standard blocks where one of the standard blocks has three boundary circles.
$$f: \ts \rightarrow S_{n}(g_{1},...,g_{n-1},x;h_{1},...,h_{n-1},y) \bigsqcup_{c,y} S_{3}(x^{-1},g_{1}',g_{2}';y,h_{1}',h_{2}')$$
Here $c$ denotes the cut where these two standard blocks are glued together. Then we have the following relation:\\
$$(id \sqcup \mathbf{B})\mathbf{T^{c}_{z,y}} = \mathbf{T^{c}_{z,y}}(id \sqcup \mathbf{B})$$
See figure \ref{btrelation} where for simplicity we take $n=4$.

\subsection{$\mathbf{T_{z,y}^{c}}$ relation}
This is simply the relation $\mathbf{T_{z,y}^{c}} =
\mathbf{F_{c,z}^{-1}F_{c,y}}$ which is expected.
\subsection{Rotation axiom}
\begin{figure}
\includegraphics{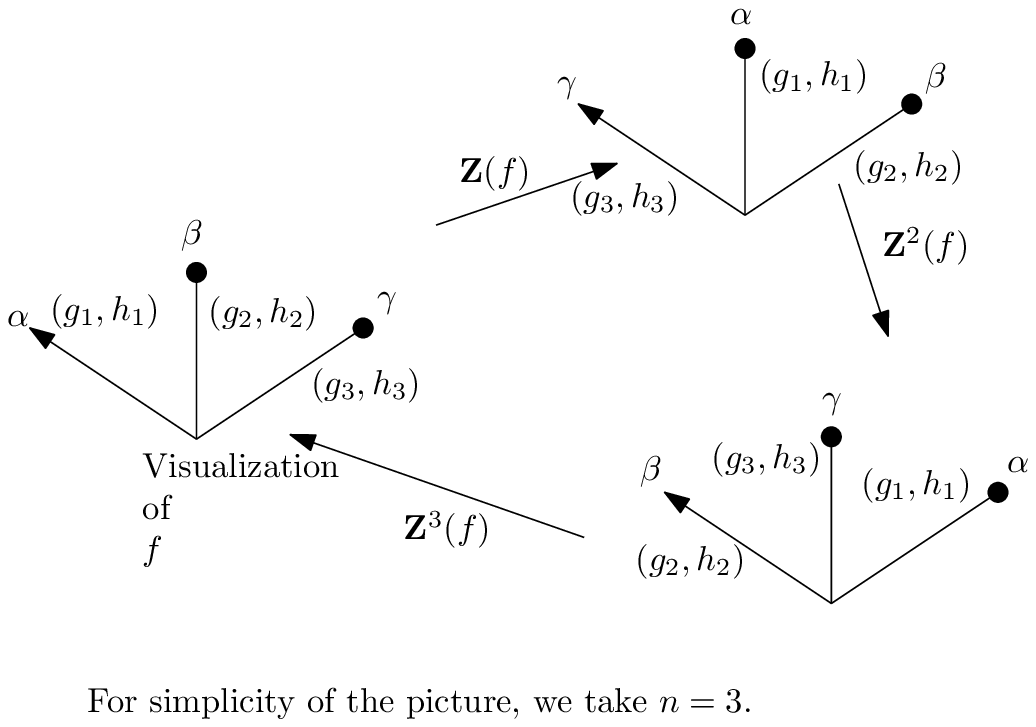}\\
\caption{The 2-cell of $\mathbf{Z}$ relation, for $n=3$}\label{zrelation}
\end{figure}

Let $f: \ts  \longrightarrow S_n(g_1...g_n;h_1...h_n)$ be a
parameterization of the $G$-cover, $\ts $. Then the rotation axiom says
that if we apply $\mathbf{Z}$ move $n$ times to the parameterization,
 $f$,
we get back the same parameterization, $f$. That is $\mathbf{Z^n}f = f$
 or
$\mathbf{Z^n} = 1$.
See the diagram on figure \ref{zrelation}.
\subsection{Commutativity of Disjoint Union}
\begin{figure}
\includegraphics{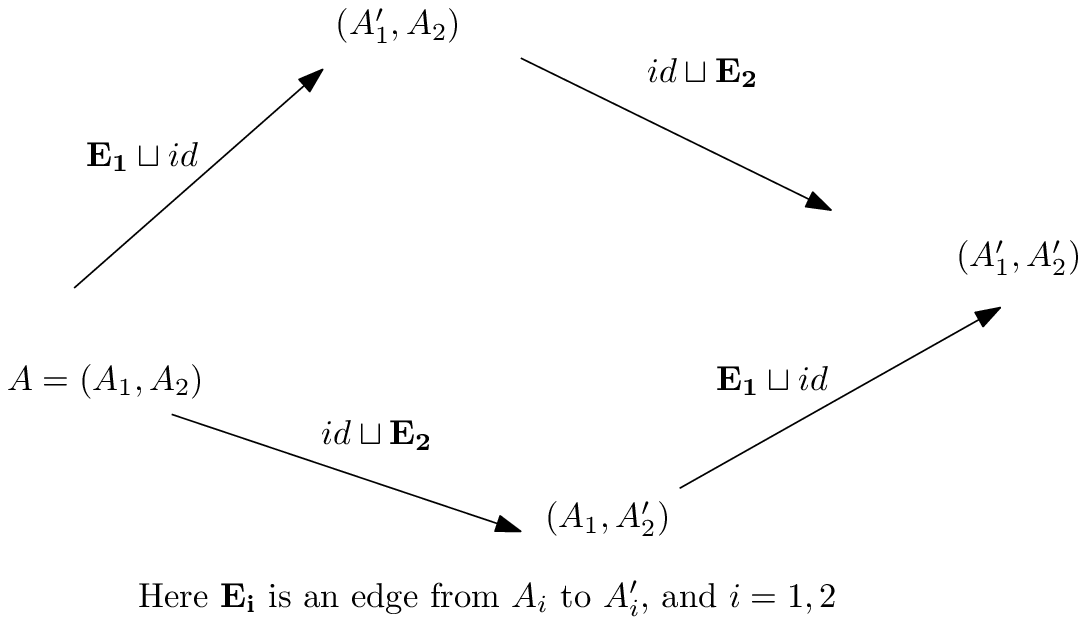}\\
\caption{The 2-cell of commutativity of disjoint union}\label{commofdisunionrelation}
\end{figure}

Let $\Sigma = \Sigma_1 \bigsqcup \Sigma_2 $. Then by restriction, the
$G$-cover $\ts  \longrightarrow \Sigma$ can be broken down into two
 pieces,
say,\\
$\tilde{\Sigma_1} \longrightarrow \Sigma_1$ and $\tilde{\Sigma_2}
\longrightarrow \Sigma_2$.\\
Thus a parameterization of $\ts  \longrightarrow \Sigma$ can be
 realized as
disjoint union of two parameterization, one for $\tilde{\Sigma_1}
\longrightarrow \Sigma_1$, and the other for $\tilde{\Sigma_2}
\longrightarrow \Sigma_2$. We can reformulate this statement in the
language of complex. In this case, it means that if $A$ is a vertex of
 the
complex $M(\ts ,\Sigma)$, then $A$ can be written as $A = ( A_1, A_2 )$
where $A_1$ is a vertex of $M(\tilde{\Sigma_1},\Sigma_1)$ and $A_2$ is
 a
vertex of $M(\tilde{\Sigma_2},\Sigma_2)$. Let $\mathbf{E_i}$ be an edge
 of
$M(\tilde{\Sigma_i},\Sigma_i)$ directed from $A_i$ to $A'_i$, here $i =
1,2$. Then the commutativity of disjoint union is the following
 relation:\\
\begin{center}
$(id \bigsqcup \mathbf{E_2})(\mathbf{E_1} \bigsqcup id) =
 (\mathbf{E_1}
\bigsqcup id)(id \bigsqcup \mathbf{E_2})$
\end{center}
See the diagram on figure \ref{commofdisunionrelation} for a visual presentation.
\subsection{Symmetry of $\mathbf{F}$ Move}
Let $f$ be the following parameterization of the $G$-cover, $\ts $,\\
$$f: \ts  \longrightarrow S_{k+1}(g_{1},...g_{k+1} ; h_{1},...h_{k+1})
\bigsqcup_{c,h_{k+1}} S_{l+1}(g'_{1},...g'_{l+1} ;
h'_{1},...h'_{l+1})$$ \\
Here the symbol, $$\bigsqcup_{c,h_{k+1}}$$ means that we glued along the cut $c$ where the edges of the cut are labeled by $(t,h_{k+1})$ and $(t^{-1},h_{k+1})$. This $t$ is completely determined by all the other labels. More precisely we have the following:\\
The boundary component labeled by $g_{k+1}$ of the standard block,
$S_{k+1}(g_{1},...g_{k+1} ; h_{1},...h_{k+1})$ is glued with the boundary component
labeled by $g'_{1}$ of the standard block $S_{l+1}(g'_{1},...g'_{l+1} ;
h'_{1},...h'_{l+1})$. of course we must have $h_{k+1} = h'_{1}$ and
$g_{k+1} = (g'_{1})^{-1}$ for gluing to exist.
From this parameterization, $f$, we can get an another
 parameterization,
$j$, where
$$j: \ts  \longrightarrow S_{l+1}(g'_{2},...,g'_{l+1},g'_{1}; h'_{2},...,h'_{l+1},h'_{1}) \bigsqcup_{c,h'_{1}} S_{k+1}(g_{k+1},g_{1},...,g_{k} ; h_{k+1},h_{1},...,h_{k})$$
This is done by first interchange the component $S_{k+1}(g_{1},...g_{k+1};h_1,...h_{k+1})$ and $S_{l+1}(g'_{1},...,g'_{l+1} ; h'_{1},...,h'_{l+1})$ and then we apply $\tilde{z}^{-1}$ to $S_{l+1}(g'_{1},...,g'_{l+1} ; h'_{1},...,h'_{l+1})$ and $\tilde{z}$ to $S_{k+1}(g_1,...,g_{k+1} ; h_1,...,h_{k+1})$ so that the boundary component labeled by $g'_{1}$ becomes the last component of $S_{l+1}(g'_{1},...g'_{l+1} ; h'_{1},...h'_{l+1})$ and the boundary component labeled by $g_{k+1}$ becomes the first component of $S_{k+1}(g_1,...g_{k+1} ; h_1,...h_{k+1})$. See Lemma 7 for the description of $\tilde{z}$\\
Now the symmetry of $\mathbf{F}$ move is the following relation:\\
\begin{center}
$\mathbf{Z^{l}F_{c,h_{k+1}}}(f) = \mathbf{F_{c,h_{1}^{'}}}(j)$
\end{center}
Or in short, we can just say\\
\begin{center}
$\mathbf{Z^{l}F_{c,h_{k+1}}} = \mathbf{F_{c,h_{1}^{'}}}$
\end{center}
Here $c$ is the cut where we glue the $G$-cover. For the notation,
$F_{c,h_{k+1}}$, see the remark on sec 6.2.3. In word, this relation
 just
says that, we can first interchange the component, apply $\mathbf{Z}$
 and
$\mathbf{Z^{-1}}$ move to the components and then apply $\mathbf{F}$
 move
or we can just first apply the $\mathbf{F}$ move and then apply some
appropriate power of $\mathbf{Z}$ move.
\subsection{Associativity of Cuts}
Let $f$ be the following parameterization of the $G$-cover $\ts $\\
$$f: \ts  \rightarrow S_{k+1}(g_1,..g_{k+1} ; h_1,..h_{k+1})
\bigsqcup_{c_{1},h_{k+1}} S_{l+1}(g'_{1},..g'_{l+1} ;
 h'_{1},..h'_{l+1}) \bigsqcup_{c_{2},h'_{l+1}} S_{n+1}(g''_{1},..g''_{n+1} ;
h''_{1},..h''_{n+1})$$
For the description of the notation,
$$\bigsqcup_{c_{1},h_{k+1}} \mbox{ and } \bigsqcup_{c_{2},h'_{l+1}}$$
see the relation``Symmetry of $\mathbf{F}$ move". Here $c_{1}$ denote
 the
first cut and $c_{2}$ denote the second cut of the gluing. Of course, we
assume that all the conditions are satisfied for gluing to exist. More
specifically, we assume the following conditions:\\
$h_{k+1} = h'_{1}$ , $h'_{l+1} = h''_{1}$ , $g_{k+1}g'_{1} = 1$ ,
$g'_{l+1}g''_{1} = 1$\\
Now the "associativity of cuts" is the following relation:\\
\begin{center}
$\mathbf{F_{c_{2},h'_{l+1}}F_{c_{1},h_{k+1}}}(f) =
\mathbf{F_{c_{1},h_{k+1}}F_{c_{2},h'_{l+1}}}(f)$
\end{center}
or in short, we can just say\\
\begin{center}
$\mathbf{F_{c_{2},h'_{l+1}}F_{c_{1},h_{k+1}}} =
\mathbf{F_{c_{1},h_{k+1}}F_{c_{2},h'_{l+1}}}$
\end{center}
In word, this just says that, given two distinct cuts, it does not
 matter
in which order we apply the $\mathbf{F}$ move.
\subsection{Cylinder Axiom}
\begin{figure}
\includegraphics{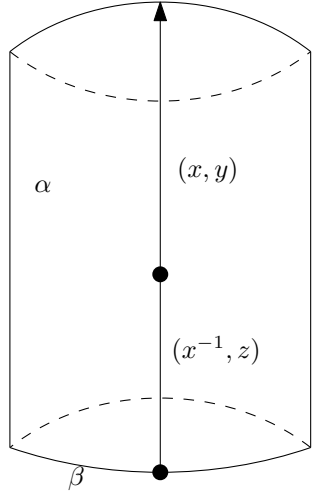}\\
\caption{Standard marking on $S_2$}\label{cylinderaxiom}
\end{figure}

We consider the standard cylinder ($S_2$) with the standard marking.
 See the picture on figure \ref{cylinderaxiom}.
Let $f: \tilde{\Omega} \longrightarrow S_2(x,x^{-1} ; y,z)$ be a
parameterization of a $G$-cover, $\tilde{\Omega} \rightarrow \Omega$.
 Also
let $j$ be the following parameterization of $G$-cover, $\ts
  \rightarrow
\Sigma$:\\
$$j: \ts  \longrightarrow S_{k+1}(g_{1},...,g_{k},x^{-1} ;
h_{1},...,h_{k},y)$$
Note that, here we choose all the elements of $G$ in such a way that we
 can
glue the G-cover $\tilde{\Omega}$ and $\ts $. So after this gluing, we
 get
a new parameterization of the $G$-cover $$\ts  \bigsqcup_{glued}
\tilde{\Omega}$$ Let us denote this parameterization by $j \sqcup f$.
 That
is \\
$$j \sqcup f : \ts  \bigsqcup_{glued} \tilde{\Omega} \rightarrow
S_{k+1}(g_{1},...,g_{k},x^{-1} ; h_{1},...,h_{k},y)
 \bigsqcup_{x^{-1},x}
S_2(x,x^{-1} ; y,z)$$
Finally let $\mathbf{E}$ be a move which can be applied to the
parameterization $j$. In the language of complex this just means that
$\mathbf{E}$ is one of the edges of the complex $M(\ts , \Sigma)$
 starting
from the vertex $j$.
Here of course, $\mathbf{E}$ must be one of the following edges\\
$\mathbf{Z,Z^{-1},B_{i},B{i}^{-1},F,F^{-1},T_{z,y}^{c},
 (T_{z,y}^{c})^{-1},
P_{x} \mbox{ or } P_{x}^{-1}}.$
Then the cylinder axiom is the following relation or 2-cell:\\
\begin{center}
$\mathbf{EF}(j \sqcup f ) = \mathbf{F( E \sqcup id )} (j \sqcup f)$
\end{center}
Or in short, we can just say\\
\begin{center}
$\mathbf{EF} = \mathbf{F( E \sqcup id )}$
\end{center}
\subsection{Braiding Axiom}
Before we define the braiding axiom, we need the following two
definition:\\
\subsubsection{Generalized $\mathbf{F}$ move}
\begin{figure}
\includegraphics{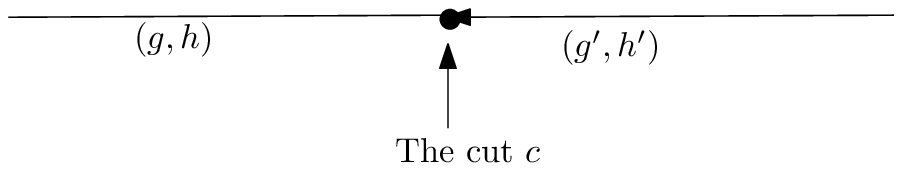}
\caption{}\label{gfmove}
\end{figure}

We know that before we can apply $\mathbf{F}$ move to a
 parameterization,
there are some assumption we must satisfy. More precisely, if $(g,h)$
 is
the label for the last boundary circle of one component, $\Sigma_{1}$,
 and
$(g',h')$ is the label for the first boundary circle of the other
component, $\Sigma_{2}$, then to apply $\mathbf{F}$ move we must
 satisfy:\\
$gg' = 1$ and $h = h'$\\
See the picture on figure \ref{gfmove}.
If to start with the boundary circle associated with $(g,h)$ is not the
last one and the boundary circle associated with $(g',h')$ is not the
 first
one, we can apply appropriate $\mathbf{Z}$ move to $\Sigma_{1}$ and
$\Sigma_{2}$ to get what we want. In other word, the generalized
$\mathbf{F}$ move is the composition of the following moves:\\

$$\overrightarrow{\mathbf{F_{c,h}}} = \mathbf{F_{c,h}(Z^{a} \sqcup
Z^{b})}$$

Note that we used an arrow above the generalized $\mathbf{F}$ move to
distinguish it from the usual $\mathbf{F}$ move. Sometimes if there is
 no
confusion, we will just use the notation $\mathbf{F_{c,h}}$ to denote a
generalized $\mathbf{F}$ move. Here $``a"$ and $``b"$ is chosen
appropriately so that $(g,h)$-boundary circle of $\Sigma_{1}$ becomes
 the
last one and $(g',h')$-boundary circle of $\Sigma_{2}$ becomes the
 first
one.
\begin{remark}
In the case $a = b = 0$, the generalized
$\overrightarrow{\mathbf{F_{c,h}}}$ move becomes the usual
$\mathbf{F_{c,h}}$ move. So we can think the usual $\mathbf{F_{c,h}}$
 move
as a special case of generalized $\overrightarrow{\mathbf{F_{c,h}}}$
 move.
 \end{remark}
\subsubsection{generalized Braiding move}
\begin{figure}
\includegraphics{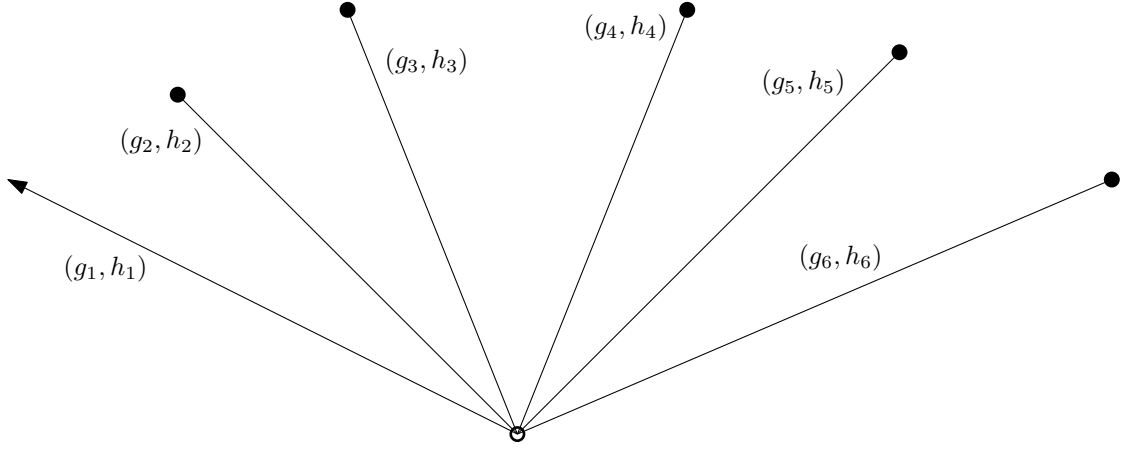}\\
\caption{Visualization of $f$}\label{gbmpic}
\end{figure}

\begin{figure}
\includegraphics{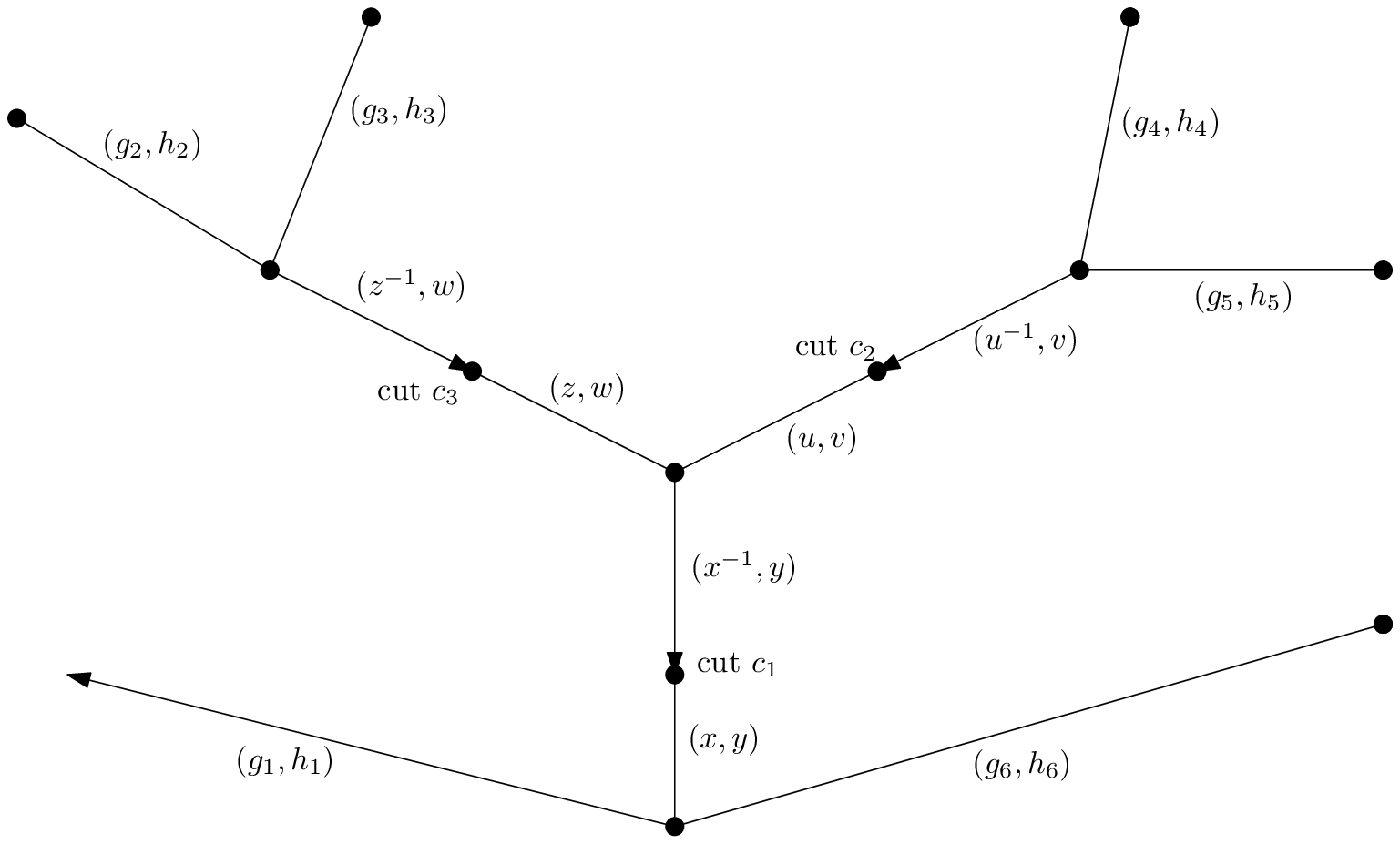}\\
\caption{Visualization of $q =
\mathbf{F_{c_{3},w}^{-1}F_{c_{2},v}^{-1}F_{c_{1},y}^{-1}}(f)$} \label{gbmpic2}
\end{figure}

\begin{figure}
\includegraphics{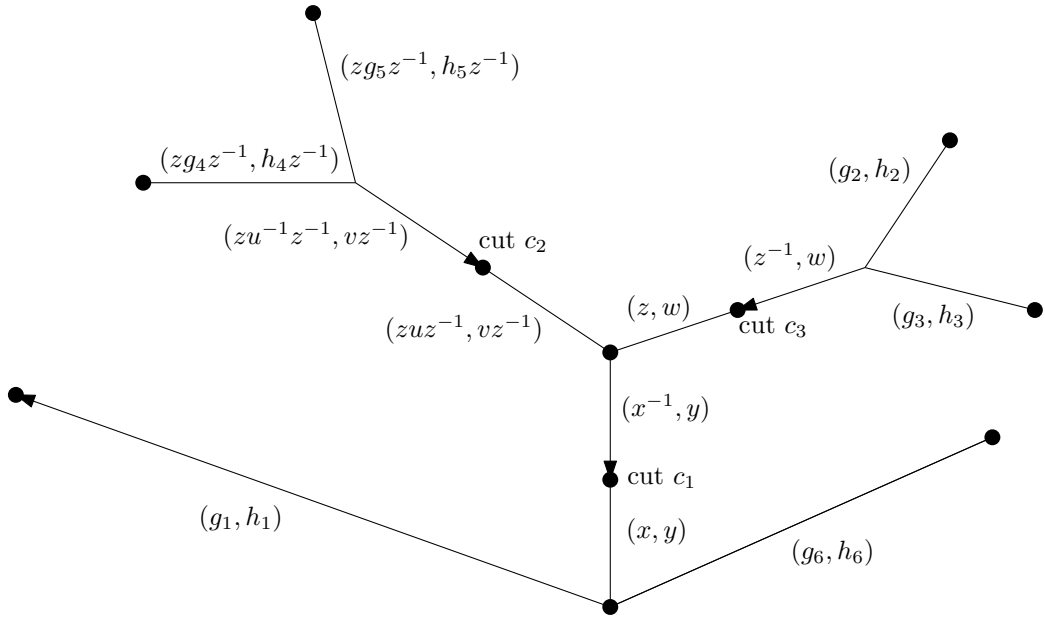}\\
\caption{Visualization of $r$}\label{gbmpic3}
\end{figure}

\begin{figure}
\includegraphics{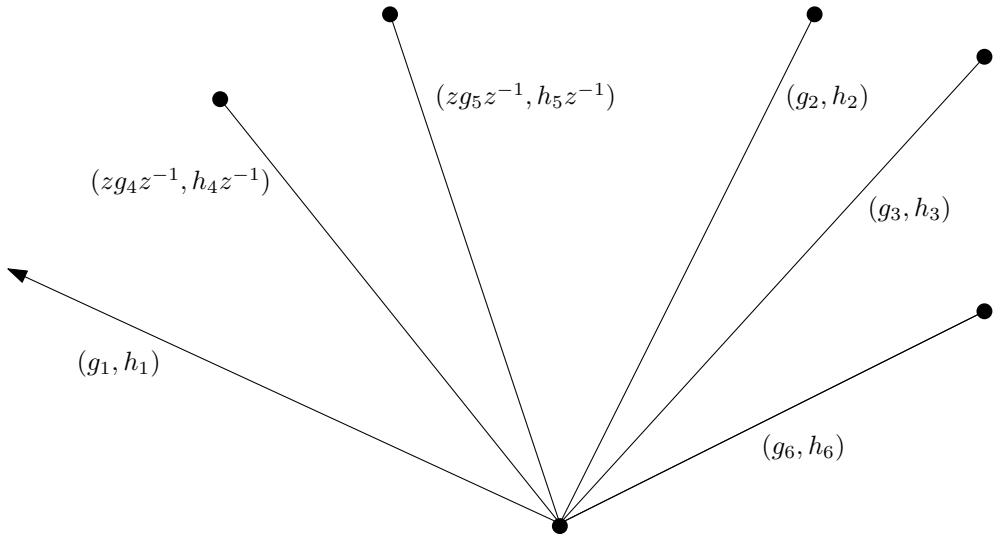}\\
\caption{Visualization of
$\mathbf{F_{c_{3},vz^{-1}}F_{c_{2},w}F_{c_{1},y}}(r)$}\label{gbmpic4}
\end{figure}

Let $f$ be the following parameterization of the $G$-cover, $\ts $, \\
$$f: \ts  \rightarrow S_{(k+l+n+m)}(g_{1},...,g_{k+l+n+m} ;
h_{1},...h_{k+l+n+m})$$
We also let \\
$I_{1} =$ \{the set of boundary circle associated with
$(g_{1},...g_{k})\}$\\
$I_{2} =$ \{the set of boundary circle associated with
$(g_{k+1},...g_{k+l})\}$\\
$I_{3} =$ \{the set of boundary circle associated with
$(g_{k+l+1},...g_{k+l+n})\}$\\
$I_{4} =$ \{the set of boundary circle associated with
$(g_{k+l+n+1},...g_{k+l+n+m})\}$\\
We want to define the generalized $\mathbf{B}$ move denoted by
$\mathbf{B_{I_{2},I_{3}}}$. The best way to explain this is through an
example. For simplicity, we take $|I_{1}| = 1, |I_{2}| = 2, |I_{3}| =
 2,
|I_{4}| = 1 $, but readers easily see that it works in general. Now
 without
further delay, we define the generalized $\mathbf{B_{I_{2},I_{3}}}$ as
 the
composition of the following moves:\\
First of all we are given:\\
$f: \ts  \rightarrow
S_{6}(\underbrace{g_{1}}_{I_{1}},\underbrace{g_{2},g_{3}}_{I_{2}},
\underbrace{g_{4},g_{5}}_{I{3}},\underbrace{g_{6}}_{I_{4}} ;
h_{1},h_{2},h_{3},h_{4},h_{5},h_{6})$.
See figure \ref{gbmpic} for the visualization of $f$.
In what follows, all the $\mathbf{F}$ and the $\mathbf{F^{-1}}$ move
 will
be generalized $\overrightarrow{\mathbf{F}}$ and
$(\overrightarrow{\mathbf{F_{c,h}}})^{-1}$ move, but we do not put
 arrow
sign on top to make the picture simple.
we apply $\mathbf{F_{c_{3},w}^{-1}F_{c_{2},v}^{-1}F_{c_{1},y}^{-1}}$ to
 the
parameterization $f$ to get the following parameterization of $\ts $\\
$$\mathbf{F_{c_{3},w}^{-1}F_{c_{2},v}^{-1}F_{c_{1},y}^{-1}}(f): \ts
\longrightarrow S_{3}(g_{1},x,g_{6} ; h_{1},y,h_{6}) \bigsqcup_{(c_{1},y)}
S_{3}(x^{-1},z,u;y,w,v) $$ $$\bigsqcup_{(c_{2},v)} S_{3}(u^{-1},g_4,g_5 ;
v,h_{4},h_{5}) \bigsqcup_{(c_{3},w)} S_3(z^{-1},g_{2},g_{3} ;
 w,h_{2},h_{3})$$
For simplicity of writing things down let us denote\\
$q = \mathbf{F_{c_{3},w}^{-1}F_{c_{2},v}^{-1}F_{c_{1},y}^{-1}}(f)$. See figure \ref{gbmpic2}.
Now to the parameterization $q$, we apply the move $id \sqcup ($the
 usual
$\mathbf{B}$ move on the edge labeled by $(u,v)$ and $(z,w)$)$ \sqcup
\mathbf{P_{z}} \sqcup id$\\
Recall that $\mathbf{P_{z}}$ means the $\mathbf{P}$ move for the
 element $z
\in G$. See section 6.3 for detail.\\
let us denote by $r$ the parameterization \{$id \sqcup ($the usual
$\mathbf{B}$ move on the edge labeled by $(u,v)$ and $(z,w)$)$ \sqcup
\mathbf{P_{z}} \sqcup id$\}($q$)\\
Here $$r: \ts  \longrightarrow S_{3}(g_{1},x,g_{6} ; h_{1},y,h_{6})
\bigsqcup_{(c_{1},y)} S_{3}(x^{-1},zuz^{-1},z;y,vz^{-1},w)$$ $$
\bigsqcup_{c_{2},vz^{-1}}
S_{3}(zu^{-1}z^{-1},zg_{4}z^{-1},zg_5z^{-1};vz^{-1},h_{4}z^{-1},h_{5}z^{-1}
) \bigsqcup_{c_{3},w} S_{3}(z^{-1},g_{2},g_{3}; w,h_{2},h_{3})$$\\
See the diagram on figure \ref{gbmpic3}.
Finally, to get rid off all the cuts, $c_{1},c_{2}, c_{3}$, that have
 been
artificially created, we apply the move,
$\mathbf{F_{c_{3},vz^{-1}}F_{c_{2},w}F_{c_{1},y}}$, to $r$. Thus
$$\mathbf{F_{c_{3},vz^{-1}}F_{c_{2},w}F_{c_{1},y}}(r): \ts
  \longrightarrow
S_6(g_{1},zg_{4}z^{-1},zg_{5}z^{-1},g_{2},g_{3},g_{6} ;
h_1,h_{4}z^{-1},h_{5}z^{-1},h_{2},h_{3},h_{6})$$
See the figure \ref{gbmpic4}.
This is the end of the generalized $\mathbf{B_{I_{2},I_{3}}}$ move. If
 we
put everything together, we have the following composition for
$\mathbf{B_{I_{2},I_{3}}}$ move:\\
$\mathbf{B_{I_{2},I_{3}}}(f) =
\mathbf{F_{c_{3},vz^{-1}}F_{c_{2},w}F_{c_{1},y}}\{id \sqcup ($the usual
$\mathbf{B}$ move on the edge labeled by $(u,v)$ and $(z,w)$)$ \sqcup
\mathbf{P_{z}} \sqcup
id\}\mathbf{F_{c_{3},w}^{-1}F_{c_{2},v}^{-1}F_{c_{1},y}^{-1}}(f)$\\
\begin{remark}
Note that we did not use any parenthesis to
denote the move $\mathbf{F_{c_{3},vz^{-1}}F_{c_{2},w}F_{c_{1},y}}$ or
$\mathbf{F_{c_{3},w}^{-1}F_{c_{2},v}^{-1}F_{c_{1},y}^{-1}}$ since by
 the
"associativity of cuts", it does not matter on which order we add or
 remove
cuts.
\end{remark}
\begin{remark}
We might wonder, what is $z$ that appear on
 the
generalized $\mathbf{B}$ move above? We can in fact calculate the value
 of
$z$ quite easily.\\
$xg_{6}g_{1} = 1 \Rightarrow x = g_{1}^{-1}g_{6}^{-1}$\\
$ux^{-1}g_{2}g_{3} = 1 \Rightarrow u =
g_{3}^{-1}g_{2}^{-1}g_{1}^{-1}g_{6}^{-1}$\\
$zux^{-1} = 1 \Rightarrow z = xu^{-1} =
g_{1}^{-1}g_{6}^{-1}g_{6}g_{1}g_{2}g_{3} = g_{2}g_{3}$\\
$\mbox{ So } z = g_{2}g_{3}$
In general $z$ = multiplication of all $g$'s in $\mathbf{I_{2}}$.\\
\end{remark}
\begin{remark}
The generalized braiding move agree with the
usual braiding move in the case when $|\mathbf{I_{2}}| =
 |\mathbf{I_{3}}| =
1$.
\end{remark}
\subsubsection{ Braiding Axiom}
\begin{figure}
\includegraphics{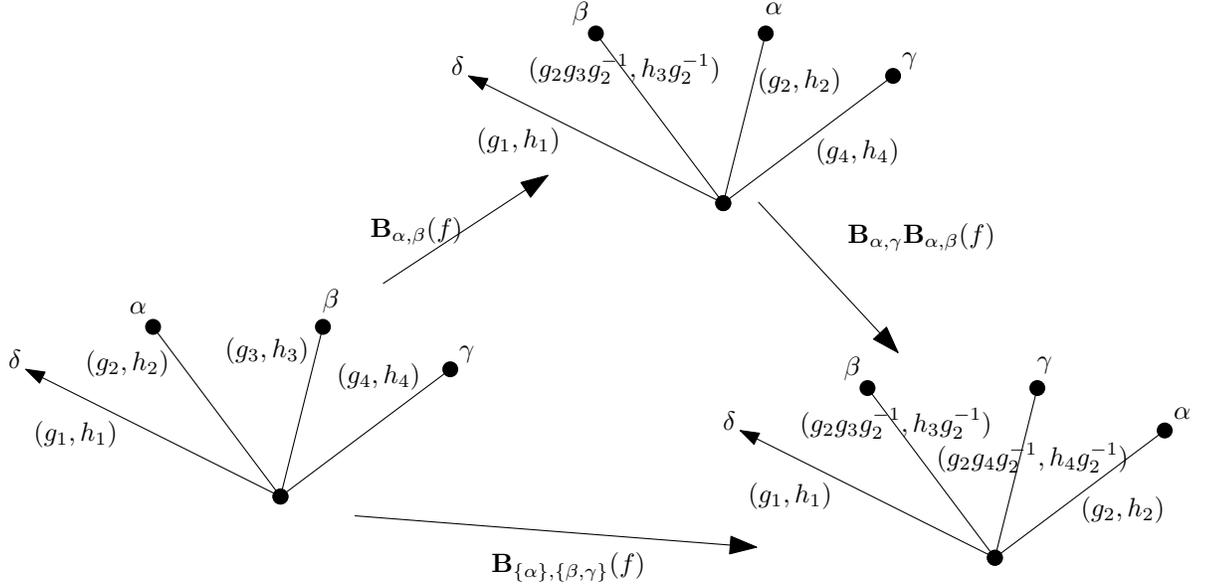}\\
\caption{2-cell or relation of braiding axiom}\label{finalbraidingaxiom}
\end{figure}

Now it is time to describe the braiding axiom. We start with the
 following
parameterization of our G-cover:\\
$f: \ts  \longrightarrow S_4(g_{1},g_{2},g_{3},g_{4} ;
h_{1},h_{2},h_{3},h_{4})$ where we name the boundary circle $\delta,
\alpha, \beta, \gamma$, in the increasing order; that is $\delta$ is
 the
boundary circle associated with $(g_{1},h_{1})$ and $\gamma$ is the
boundary circle associated with $(g_{4},h_{4})$. \\
Then the Braiding axiom is the following two relation:\\
\begin{center}
$\mathbf{B_{\alpha,\gamma}B_{\alpha,\beta}}(f) =
\mathbf{B_{\{\alpha\},\{\beta,\gamma\}}}(f)$
\end{center}
and\\
\begin{center}
$\mathbf{B_{\alpha,\gamma}B_{\beta,\gamma}}(f) =
\mathbf{B_{\{\alpha\,\beta\},\{\gamma\}}}(f)$
\end{center}
Here $\mathbf{B_{\alpha,\gamma}}$ etc denote the usual $\mathbf{B}$
 move
and $\mathbf{B_{\{\alpha\},\{\beta,\gamma\}}}$ etc denote the
 generalized
$\mathbf{B}$ move. \\
We will describe step by step move for the first relation of Braiding
 Axiom
since the description of the other relation is similar. \\
$$\mathbf{B_{\alpha,\beta}}(f): \ts  \longrightarrow
S_{6}(g_{1},g_{2}g_{3}g_{2}^{-1},g_{2},g_{4} ;
h_{1},h_{3}g_2^{-1},h_{2},h_{4})$$
$$\mathbf{B_{\alpha,\gamma}}\mathbf{B_{\alpha,\beta}}(f): \ts
\longrightarrow
 S_{6}(g_{1},g_{2}g_{3}g_{2}^{-1},g_{2}g_{4}g_{2}^{-1},g_{2}
; h_{1},h_{3}g_{2}^{-1},h_{4}g_{2}^{-1},h_{2})$$\\
Look at the definition of $\mathbf{B}$ move above. On the other hand
$$\mathbf{B_{\{\alpha\},\{\beta,\gamma\}}}(f): \ts  \longrightarrow
S_{6}(g_{1},g_{2}g_{3}g_{2}^{-1},g_{2}g_{4}g_{2}^{-1},g_{2} ;
h_{1},h_{3}g_{2}^{-1},h_{4}g_{2}^{-1},h_{2})$$\\
See the diagram on figure \ref{finalbraidingaxiom}  for a visual presentation.
\subsection{Dehn Twist Axiom}
\begin{figure}
\includegraphics{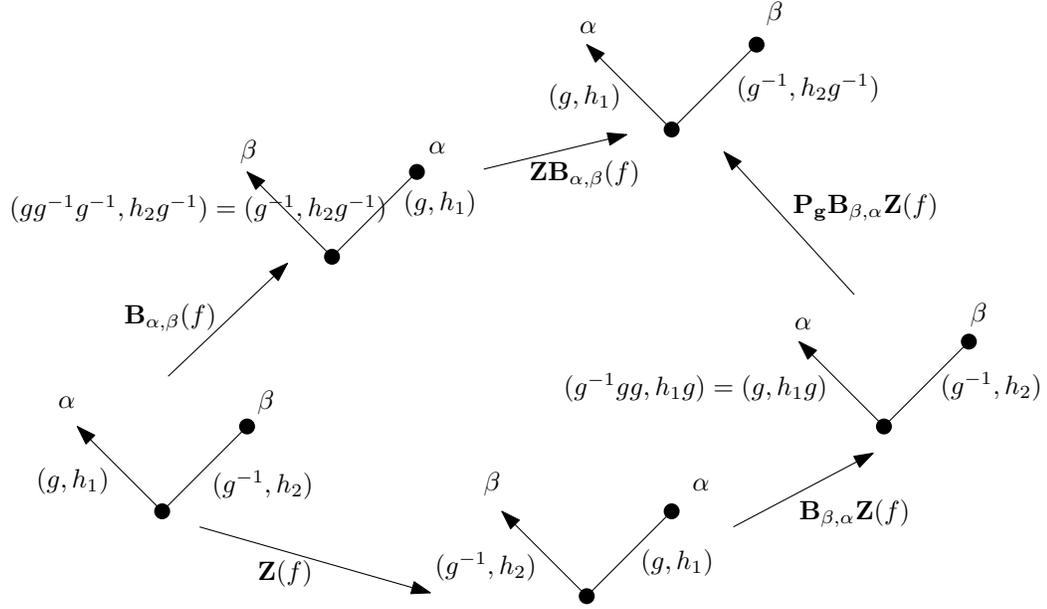}\\
\caption{2-cell or relation of Dehn twist axiom}\label{dehntwistaxiom}
\end{figure}

Say, we are given a Standard Cylinder ($S_{2}$) with the standard graph
(see the picture below) where we denote by $\alpha$, the first boundary
component and by $\beta$, the second boundary component. Also let $f$
 be a
parameterization of our G-cover, $\ts $, given by\\
$f: \ts  \longrightarrow S_{2}(g,g^{-1} ; h_{1},h_{2})$\\
Then the dehn twist axiom is the following relation:\\
\begin{center}
$\mathbf{ZB_{\alpha,\beta}}(f) = \mathbf{P_{g}B_{\beta,\alpha}Z}(f)$
\end{center}
Here $\mathbf{B_{\alpha,\beta}}$ denote the usual $\mathbf{B}$ move and
$\mathbf{P_{g}}$ denote the $\mathbf{P}$ move for the element $g \in
 G$.
For a detail description of these moves, see section 6.
For a better understanding, we breakdown this relation piece by piece:
$$\mathbf{B_{\alpha,\beta}}(f): \ts  \longrightarrow
 S_{2}(gg^{-1}g^{-1},g
; h_{2}g^{-1},h_{1}) = S_{2}(g^{-1},g ; h_{2}g^{-1},h_{1})$$
$$\mathbf{ZB_{\alpha,\beta}}(f): \ts  \longrightarrow S_{2}(g,g^{-1} ;
h_{1},h_{2}g^{-1})$$
On the other hand\\
$$\mathbf{Z}(f): \ts  \longrightarrow S_{2}(g^{-1},g ; h_{2},h_{1})$$
$$\mathbf{B_{\beta,\alpha}Z}(f): \ts  \longrightarrow
 S_{2}(g^{-1}gg,g^{-1}
; h_{1}g,h_{2}) = S_{2}(g,g^{-1} ; h_{1}g,h_{2})$$
$$\mathbf{P_{g}B_{\beta,\alpha}Z}(f): \ts  \longrightarrow
S_{2}(ggg^{-1},gg^{-1}g^{-1} ; h_{1}gg^{-1},h_{2}g^{-1}) = S_2(g,g^{-1}
 ;
h_{1},h_{2}g^{-1})$$
See the diagram on figure \ref{dehntwistaxiom}  for a visual description of the Dehn Twist
 Axiom.
\begin{remark}
Note that, here this dehn twist axiom differs
 from
the dehn twist axiom of the paper \cite{BK} (the case when $|G| = 1$),
 because
of this extra $\mathbf{P_{g}}$ move appearing in the relation.
\end{remark}
We conclude this section with the following lemma\\
\begin{lemma}
All the relations describe above make sense. In other word, all the
relation describe closed loop in our complex, M($\ts ,\Sigma$).
\end{lemma}
\begin{proof}
This is basically done when we explained the relation above. We can
 review
each relation one by one and infact see that they are closed loop in
 our
complex. The detail are left to the reader.
\end{proof}
\section{Proving the Main Result}
\subsection{A general Theorem about Complexes}
\begin{theorem}
We will use this result to prove our main theorem. Let $A$ and $B$ be
 two
2-dim complex. Let $\Pi : B^{[1]} \longrightarrow A^{[1]}$ be a map of
their 1-skeleton(vertices and edges), which is surjective both on
 vertices
and edges. Also suppose the following condition is satisfied:\\
\begin{enumerate}
\item {$A$ is connected and simply-connected}
\item{For every vertex $a \in A$, $\Pi^{-1}(a)$ is connected and
simply-connected. That is every loop, completely lie in $\Pi^{-1}(a)$,
 is
contractable in $B$. }
\item{Let $$b_{1}^{'} \stackrel{f^{'}}{\rightarrow} b_{2}^{'}$$ and
$$b_{1}^{''}\stackrel{f^{''}}{\rightarrow} b_{2}^{''}$$ be two lifting
 of
$$a_{1} \stackrel{f}{\rightarrow} a_{2}$$ Then there is a path $e_{1}$,
starting from $b_{1}^{''}$ and end with $b_{1}^{'}$, and completely lie
 in
$\Pi^{-1}(a_{1})$ and a path $e_{2}$, starting from $b_{2}^{''}$ and
 end
with $b_{2}^{'}$, and completely lie in $\Pi^{-1}(a_{2})$, so that the
following relation hold:\\
\begin{center}
$e_{2}f^{''} = f^{'}e_{1}$
\end{center}
In other word, the following diagram commute:\\
\begin{center}
$b_{1}^{''} \stackrel{f^{''}}{\rightarrow} b_{2}^{''}$\\
$e_{1} \downarrow \hspace{1 cm} \downarrow e_{2}$\\
$b_{1}^{'} \stackrel{f^{'}}{\rightarrow} b_{2}^{'}$\\
\end{center}}
\item{Every loop in $A$ can be lifted to a contractible loop in $B$.}
\end{enumerate}
Then the complex $B$ is connected and simply-connected.\\
\end{theorem}
\begin{proof}
Not hard and we will leave the proof to the reader.
\end{proof}
\subsection{Proof Of the Main Theorem}
Recall that the main theorem says that the complex $M(\ts , \Sigma)$ is
connected and simply-connected. To use the previous theorem, we let \\
$A = M(\Sigma)$:= (as defined in the paper \cite{BK}.)\\
\begin{itemize}
\item {Vertex of $M(\Sigma)$ = Marking with cuts}
\item{Edges of $M(\Sigma)$ = $Z, B, F$ move described in the paper
\cite{BK}}
\item{Relation = Described in the paper \cite{BK}}
\end{itemize}
We also let $B = M(\ts , \Sigma)$ := see section 7 for the definition
 of
this complex.
\subsubsection{Description of $\Pi : M(\ts , \Sigma)^{[1]} \rightarrow
M(\Sigma)^{[1]}$}
We first describe the map of 1-skeleton. If $v$ is a vertex in $M(\ts ,
\Sigma)$, then $v$ is given by marked graph, with $(g_{i},h_{i})$ in
 each
cut. $\Pi$ just take this vertes $v$ to the marked graph, forgetting
 about
$(g_{i},h_{i})$. This marked graph without $(g_{i},h_{i})$ is a vertex
 in
$M(\Sigma)$\\
If $\mathbf{E}$ is an edge in $M(\ts , \Sigma)$, where $\mathbf{E}$ is
 one
of $\mathbf{Z,B_{i}}$ or $\mathbf{F}$ move, then $\Pi(\mathbf{E})$ is
 the
corresponding $Z,B_{i}$ and $F$ move in $M(\Sigma)$.
Since a $\mathbf{P_g}, g\in G$ and $\mathbf{T_{z,y}^{c}}$ move does not
change the underling marked graph, $\Pi(\mathbf{P_{g}}) =
\Pi(\mathbf{T_{z,y}^{c}})$ = empty edge. This conclude the description
 of
$\Pi$. Note that $\Pi$ is obviously surjective both on the vertex and
edges.
\begin{lemma}
The above map of complex, $\Pi : M(\ts , \Sigma) \rightarrow M(\Sigma)$
make sense. In other word, $\Pi$ does define a map of complex.
\end{lemma}
\begin{proof}
obvious.
\end{proof}
\subsubsection{Description of a typical fiber}
Let $a \in$ vertex ($M(\Sigma)$). We want to describe $\Pi^{-1}(a)$.
 Hence
the following lemma:\\
\begin{lemma}
Let $b$ and $b^{'}$ $\in \Pi^{-1}(a)$. Then $b$ and $b^{'}$ is
 connected by
a sequence of moves of the form $\mathbf{F_{c_{i},z}^{-1}F_{c_{i},y}}
 :=
\mathbf{T_{z,y}^{c_{i}}}$ and $\mathbf{P_{x}}$. See section 6 for the
description of these moves and the notation. Here $c_{i}$ is the i-th
 cut
and $z,y,x \in G$. Also conversely, only move that takes a vertex of a
fiber to the same fiber, is the move $\mathbf{T_{z,y}^{c_{i}}}$ and
$\mathbf{P_{x}}$.\\
\end{lemma}
\begin{proof}
Let the vertex $a \in M(\Sigma)$ be parameterized by $k$ cuts. That
 is
$$\Sigma = \Sigma_{1} \bigsqcup_{ \mbox{ glued }} \Sigma_{2}
 \bigsqcup_{
\mbox{ glued }} ......\bigsqcup_{ \mbox{ glued }} \Sigma_{k}$$
Then the corresponding G-cover can also be break down as a gluing of
 $k$ G-cover. That is $$\ts  = \ts _{1} \bigsqcup_{ \mbox{ glued }} \ts _{2}
\bigsqcup_{ \mbox{ glued }} ......\bigsqcup_{ \mbox{ glued }} \ts
 _{k}$$
Then the parameterization, $b$ and $b^{'}$ can be realized as a gluing
 of
$k$ parameterization; each one coming from the parameterization of
 $\ts
_{i} \longrightarrow \Sigma_{i}$ which we will denote by $b_{i}$ and
$b_{i}^{'}$. So in short hand notation we can write $b = b_{1} \sqcup
 b_{2}
\sqcup...\sqcup b_{k}$ and $b^{'} = b_{1}^{'} \sqcup b_{2}
 \sqcup...\sqcup
b_{k}^{'}$ . Here $i=1...k$. Let
$$b: \ts  \longrightarrow S_{n_{1}}(g_{1}^{1},...g_{n_{1}}^{1} ;
h_{1}^{1},...h_{n_{1}}^{1}) \bigsqcup_{c_{1},h_{n_{1}}^{1}} .....\bigsqcup_{c_{k},h_{n_{k-1}}^{k-1}} S_{n_{k}}(g_{1}^{k},...g_{n_{k}}^{k} ;
h_{1}^{k},...h_{n_{k}}^{k})$$ and let
$$ b^{'}: \ts  \longrightarrow S_{n_{1}}(p_{1}^{1},...p_{n_{1}}^{1} ;
q_{1}^{1},...q_{n_{1}}^{1}) \bigsqcup_{c_{1},q_{n_{1}}^{1}} .....\bigsqcup_{c_{k},q_{n_{k-1}}^{k-1}} S_{n_{k}}(p_{1}^{k},...p_{n_{k}}^{k} ;
q_{1}^{k},...q_{n_{k}}^{k})$$
consider the parameterization $$b_{1}: \ts _{1} \longrightarrow
S_{n_{1}}(g_{1}^{1},...g_{n_{1}}^{1} ; h_{1}^{1},...h_{n_{1}}^{1})$$
 and
$$b_{1}^{'}: \ts _{1} \longrightarrow
 S_{n_{1}}(p_{1}^{1},...p_{n_{1}}^{1}
; q_{1}^{1},...q_{n_{1}}^{1})$$ Then
$$b_{1}{'}b_{1}^{-1} : S_{n_{1}}(g_{1}^{1},...g_{n_{1}}^{1} ;
h_{1}^{1},...h_{n_{1}}^{1}) \longrightarrow
S_{n_{1}}(p_{1}^{1},...p_{n_{1}}^{1} ; q_{1}^{1},...q_{n_{1}}^{1})$$ is
 an
isomorphism of G-cover. So from lemma 3, we know that there exist an $x_{1}
 \in
G$ so that $S_{n_{1}}(x_{1}g_{1}^{1}x_{1}^{-1},...x_{1}g_{n_{1}}^{1}x_{1}^{-1} ;
h_{1}^{1}x_{1}^{-1},...h_{n_{1}}^{1}x_{1}^{-1}) =
S_{n_{1}}(p_{1}^{1},...p_{n_{1}}^{1} ; q_{1}^{1},...q_{n_{1}}^{1})$ \\
Similarly we can find $x_{2},x_{3},...x_{k}$ so that $x_{i}g_{j}^{i}x_{i}^{-1} = p_{j}^{i}$ and $h^{i}_{j}x_{i}^{-1} = q_{j}^{i}$ for $i = 2...k$.
Now let us apply the $\mathbf{P_{x_{1}} \sqcup .... \sqcup P_{x_{k}}}$ move to $b$, Then
$$\mathbf{P_{x_{1}} \sqcup .... \sqcup P_{x_{k}}}(b): \ts  \longrightarrow
S_{n_{1}}(x_{1}g_{1}^{1}x_{1}^{-1},...x_{1}g_{n_{1}}^{1}x_{1}^{-1} ;
h_{1}^{1}x_{1}^{-1},...h_{n_{1}}^{1}x_{1}^{-1}) \bigsqcup_{c_{1},h_{n_{1}}^{1}x_{1}^{-1}}
.....\bigsqcup_{c_{k},h_{n_{k-1}}^{k-1}x_{k-1}^{-1}}$$ $$
S_{n_{k}}(x_{k}g_{1}^{k}x_{k}^{-1},...x_{k}g_{n_{k}}^{k}x_{k}^{-1} ;
h_{1}^{k}x_{k}^{-1},...h_{n_{k}}^{k}x_{k}^{-1})$$\\
By the choice of $x_{1},x_{2},...,x_{k}$, we must have, $x_{i}g_{i}^{j}x_{i}^{-1} = p_{i}^{j}$ and $h^{i}_{j}x_{i}^{-1} = q_{j}^{i}$ as $j = 1...k$ and $i = 1...n_{j}$.\\
In other word, $\mathbf{P_{x_{1}} \sqcup .... \sqcup P_{x_{k}}}(b) = b'$.\\

Now conversely, let $b \in \Pi ^{-1}(a)$ and $\mathbf{E}$ is a move in
 \mm,
so that $\mathbf{E}(b) \in \Pi ^{-1}(a)$.
At first site, of course, all possibilities for $\mathbf{E}$ are
$\mathbf{Z,Z^{-1}B,B^{-1}F,F^{-1}}$ $\mathbf{P,P^{-1},T \mbox{ and }
T^{-1}}$. But the move $\mathbf{Z,Z^{-1}B,B^{-1}F,F^{-1}}$ will take
 the
vertex, $b$, outside the fiber $\Pi ^{-1}(a)$. So the only possible
 value
for $\mathbf{E}$ is $\mathbf{P_{x}}$ and \T move and their inverses.
 But
$\mathbf{P_{x}^{-1}} = \mathbf{P_{x^{-1}}}$ and
$(\mathbf{T_{z,y}^{c}})^{-1} = (\mathbf{F_{c,z}^{-1}F_{c,y}})^{-1} =
\mathbf{F_{c,y}^{-1}F_{c,z}} = \mathbf{T_{y,z}^{c}}$. So the inverse of
$\mathbf{P,T}$ move are another $\mathbf{P,T}$ move.
\end{proof}
\emph{Remark}: This lemma shows in particular that the fiber
 $\Pi^{-1}(a)$
is connected.
\subsubsection{The fiber $\Pi^{-1}(a)$ is simply-connected}
\begin{lemma}
The fiber $\Pi^{-1}(a)$ is simply-connected. That is given any loop,
 where
each vertex of the loop belong to $\Pi^{-1}(a)$, is contractable using
 the
2-cell or relations of $M(\ts ,\Sigma)$.
\end{lemma}
\begin{proof}
We divide the proof in three cases. \\
\underline{\emph{case 1: no cuts and no boundary circle}}\\
Here the base surface, $\Sigma$, is isomorphic to $S_{0}$, the sphere,
 and
the G-cover of $\Sigma$ is trivial that is $\Sigma \times G$. This case
must be treated separately. Although there is only one G-cover of
 $\Sigma$
up to isomorphism, there are a total of $|G|$ many parameterization of
 this
G-cover. To see this, recall that a parameterization, $f$, of our
 G-cover
$\Sigma \times G$ is just an isomorphism from the G-cover $\Sigma
 \times G$
to the G-cover $S_{0} \times G$. of course this $f$ must maps a
 component
of $\Sigma \times G$ isomorphically to a component of $S_{0} \times G$.
 Let
say that $$f(\Sigma \times 1, \mbox{ where } 1 \in G) = S_{0} \times h,
\mbox{ where } h \in G$$
Then this information will determine $f$ completely since f must
 preserve
the action of $G$ on the fiber. More precisely we have:
$$f(\Sigma \times x, \mbox{ where } x \in G) = S_{0} \times xh$$
So a parameterization is completely determine by an element of $h \in
 G$.
So in this way, we can identify the set of parameterization of our
 G-cover
to the group $G$. That is, in this case, vertex ($M(\ts ,\Sigma)$) =
 $G$\\
What about the moves and relation? None of the moves,
$\mathbf{Z,B_{i},F_{c,y},P_{x} T_{z,y}^{c}}$, make sense in this case
 since
we do not have any cuts. But nevertheless, it is possible to define
 moves.
Recall that in general, given a parameterization, $f$, of our G-cover,
applying a move to $f$ means we compose $f$ with some standard
 automorphism
of our "Standard Block", G-cover of $S_{n}$. In this case, the standard
block is just$ S_{0} \times G$. So what are all the automorphism of the
$S_{0} \times G$? Again by the same argument as above we can identify
 the
set of automorphism of $S_{0} \times G$ with $G$. So moves of $M(\ts
,\Sigma)$ = $G$. More precisely we connect $x$ to $yx$ by an edge
 directed
from $x$ to $yx$.
only relation here is precisely the relation satisfied by the group.
Trivially this complex is connected and simply-connected (any closed
 loop
starting from vertex $x$ has the form $g_1 g_2 ...g_k x = x$ but then
 $g_1
g_2 ...g_k = 1$ which is a relation of our complex). \\
\underline{\emph{case 2: no cuts but at least one boundary circle}}\\
In this case, the base surface, $\Sigma$, is still simply-connected
(remember that we always assume $\Sigma$ has genus $0$). So the G-cover
 of
this is again $\Sigma \times G$. So all the argument above goes through
 and
we see that vertex ($M(\ts ,\Sigma)$) = $G$ and the edges of $M(\ts
,\Sigma)$ = $G$. And the complex $M(\ts ,\Sigma)$ is connected and
simply-connected.\\
\underline{\emph{case 3: at least one cut}}\\
This is the general situation. First recall the $\mathbf{P_{x}}$
relation:\\
\begin{enumerate}
\item {$\mathbf{P_xZ} = \mathbf{ZP_x}$}
\item{$\mathbf{P_xB_i} = \mathbf{B_iP_x}$}
\item{$\mathbf{P_xF_{c,y}} = \mathbf{F_{c,yx^{-1}}(P_{x} \sqcup P_{x})}$}
\item{$\mathbf{P_xP_y} = \mathbf{P_{xy}}$}
\end{enumerate}
and also the relation $\mathbf{T_{z,y}^{c}} =
 \mathbf{F_{c,z}^{-1}F_{c,y}}$.\\
Now given a loop completely lie inside the fiber $\Pi^{-1}(a)$, we know
from lemma 10 that this loop consist entirely of moves of the form
$\mathbf{P_{x}}$ and $\mathbf{F_{c_{i},z}^{-1}F_{c_{i},y}} =
\mathbf{T_{z,y}^{c_{i}}}$. First we translate all the move of the form
$\mathbf{P_{x}}$ to the right and combine all the $\mathbf{P}$ move
together to create a single $\mathbf{P}$ move. This can be done by the
above three $\mathbf{P_{x}}$-relation. So we may assume that the loop
 looks
like :\\
$$\mathbf{T_{z_{1,1},y_{1,1}}^{c_{1}}}
\mathbf{T_{z_{1,2},y_{1,2}}^{c_{1}}}....
\mathbf{T_{z_{1,n_{1}},y_{1,n_{1}}}^{c_{1}}}......
\mathbf{T_{z_{k,1},y_{k,1}}^{c_{k}}}
\mathbf{T_{z_{k,2},y_{k,2}}^{c_{k}}}....
\mathbf{T_{z_{k,n_{k}},y_{k,n_{k}}}^{c_{k}}} \mathbf{P_{x} \sqcup P_{x} \sqcup .... \sqcup P_{x}} (b)$$
Here we choose $b$ to be the starting vertex for our loop. Note that the indices for the $\mathbf{P}$ move are all same, which we denote by $x$. This must be the case for $\mathbf{T}$ move to make sense. Note also that that we
gather all $\mathbf{T}$ move associated to a given cut $c$. This can be
done since first of all, the associativity of cuts says that
$\mathbf{F_{c_{j},y}F_{c_{i},z}} = \mathbf{F_{c_{i},z}F_{c_{j},y}}$
whenever $i$ is different from $j$. From this it follows that
$\mathbf{T_{z_{1},y_{1}}^{c_{i}}} \mathbf{T_{z_{2},y_{2}}^{c_{j}}} =
\mathbf{T_{z_{2},y_{2}}^{c_{j}}} \mathbf{T_{z_{1},y_{1}}^{c_{i}}}$
 whenever
$i$ is different from $j$. Thus we can bring all $\mathbf{T}$ move
associated to a given cut in one place. Now consider the $i$th cut,
$c_{i}$. What can we say about $\mathbf{T_{z_{i,1},y_{i,1}}^{c_{i}}}
\mathbf{T_{z_{i,2},y_{i,2}}^{c_{i}}}.....\mathbf{T_{z_{i,n_{i}},y_{i,n_{i}}
}^{c_{i}}}$? We use the relation $\mathbf{T_{z,y}^{c}} =
\mathbf{F_{c,z}^{-1}F_{c,y}}$ to conclude that
$\mathbf{T_{z_{i,1},y_{i,1}}^{c_{i}}}
\mathbf{T_{z_{i,2},y_{i,2}}^{c_{i}}}.....\mathbf{T_{z_{i,n_{i}},y_{i,n_{i}}
}^{c_{i}}}$ = $\mathbf{F_{c_{i},z_{i,1}}^{-1}F_{c_{i},y_{i,n_{i}}}}$.
 All
the middle part will be $\mathbf{FF^{-1} = 1}$. So now our loop has the
form\\
$$\mathbf{F_{c_{1},z_{1,1}}^{-1}F_{c_{1},y_{1,n_{1}}}}
\mathbf{F_{c_{2},z_{2,1}}^{-1}F_{c_{2},y_{2,n_{2}}}}....\mathbf{F_{c_{k},z_
{k,1}}^{-1}F_{c_{k},y_{k,n_{k}}}} \mathbf{P_{x} \sqcup P_{x} \sqcup .... \sqcup P_{x}}(b)$$
or
$$\mathbf{T_{z_{1,1},y_{1,n_{1}}}^{c_{1}}}
\mathbf{T_{z_{2,1},y_{2,n_{2}}}^{c_{2}}}....\mathbf{T_{z_{k,1},y_{k,n_{k}}}
^{c_{k}}} \mathbf{P_{x} \sqcup P_{x} \sqcup .... \sqcup P_{x}}(b)$$
Now consider the $i$th cut, $c_{i}$. Say this cut is labeled by
$(u_{i},y_{i,n_{i}})$ for one component and by
 $(u_{i}^{-1},y_{i,n_{i}})$
for the other component. When we apply $\mathbf{P_{x} \sqcup P_{x} \sqcup .... \sqcup P_{x}}$ move to $b$,
 this
cut will be relabeled by $(xu_{i}x^{-1},y_{i,n_{i}}x^{-1})$ and
$(xu_{i}^{-1}x^{-1},y_{i,n_{i}}x^{-1})$ respectively. See how we change
 the
label for $\mathbf{P_{x}}$ move in section[7]. Then we apply
$\mathbf{T_{z_{i,1},y_{i,n_{i}}}^{c_{i}}}$ move to get the label
$(xu_{i}x^{-1},z_{i,1})$ and $(xu_{i}^{-1}x^{-1},z_{i,1})$. See the
sequence of picture on figure \ref{fibersimplyconnected} for a visual description.
\begin{figure}
\includegraphics{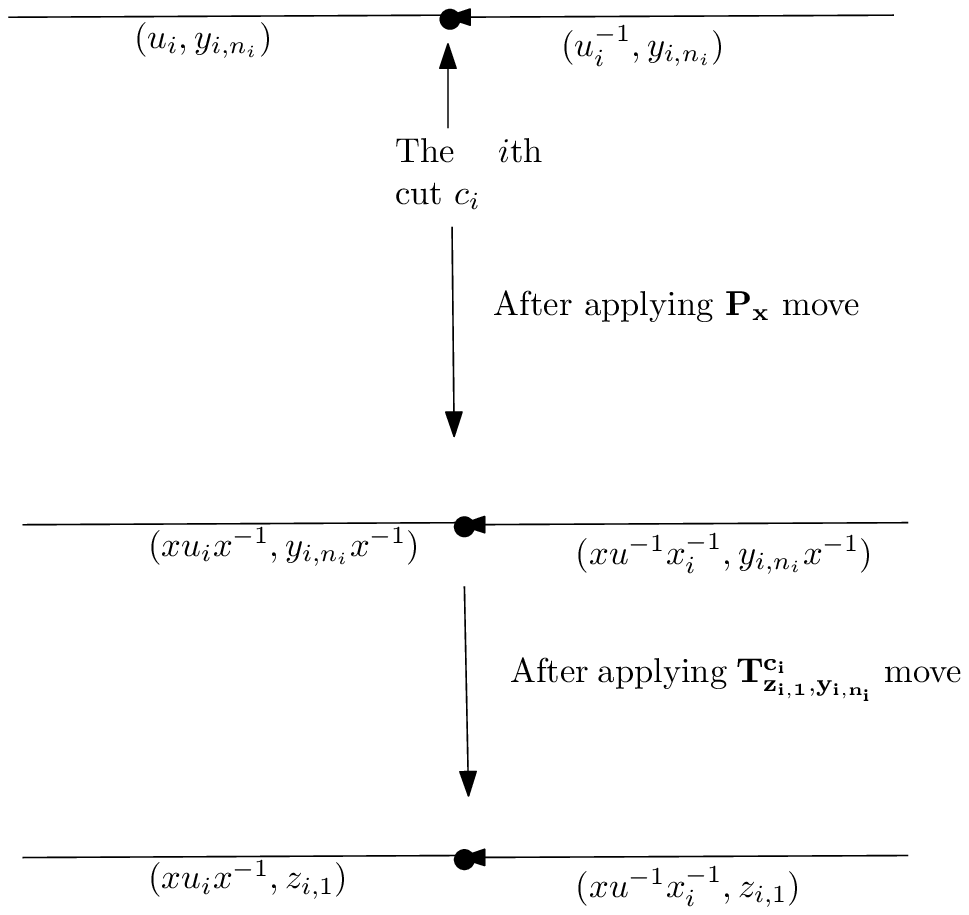}
\caption{}\label{fibersimplyconnected}
\end{figure}
But since it is a closed loop, we must have
$$xu_{i}x^{-1} = u_{i} \mbox{ and } z_{i,1} = y_{i,n_{i}}$$
as $i=1...k$
Now we go back to the beginning and rewrite our original loop but this
 time
we will move $\mathbf{P_{x} \sqcup P_{x} \sqcup .... \sqcup P_{x}}$ to all the way left. Recall the relation
$\mathbf{P_xF_{c,y}} = \mathbf{F_{c,yx^{-1}}(P_{x} \sqcup P_{x})}$. This will imply the
following two relation:\\
$$\mathbf{P_xF_{c,yx}} = \mathbf{F_{c,y}(P_x \sqcup P_x)} \mbox{ and }
\mathbf{P_xF_{c,yx}^{-1}} = \mathbf{F_{c,y}^{-1}(P_x \sqcup P_x)} $$
These above two relation together with $\mathbf{T_{z,y}^{c}} =
\mathbf{F_{c,z}^{-1}F_{c,y}}$ will imply the following:\\
$$\mathbf{T_{z,y}^{c}(P_{x} \sqcup P_{x})} = \mathbf{(P_{x} \sqcup P_{x})T_{zx,yx}^{c}}$$
See the figure \ref{ptrelation} for a visual description of the above relation.
\begin{figure}
\includegraphics{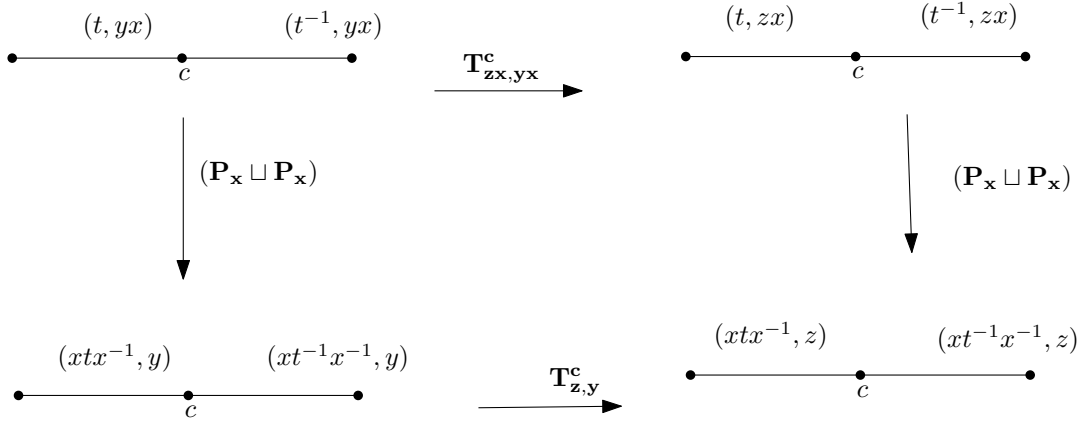}
\caption{Visual description of what happens when we interchange $\mathbf{T}$ and $\mathbf{P}$ moves} \label{ptrelation}
\end{figure}
So when we move $\mathbf{P_{x} \sqcup P_{x}... \sqcup P_{x}})$ to the left, our loop will look
 like:\\
$$\mathbf{P_{x} \sqcup P_{x} \sqcup .... \sqcup P_{x}} \mathbf{F_{c_{1},z_{1,1}x}^{-1}F_{c_{1},y_{1,n_{1}}x}}
\mathbf{F_{c_{2},z_{2,1}x}^{-1}F_{c_{2},y_{2,n_{2}}x}}....\mathbf{F_{c_{k},
z_{k,1}x}^{-1}F_{c_{k},y_{k,n_{k}}x}} (b)$$
or
$$\mathbf{P_{x} \sqcup P_{x} \sqcup .... \sqcup P_{x}} \mathbf{T_{z_{1,1}x,y_{1,n_{1}}x}^{c_{1}}}
\mathbf{T_{z_{2,1}x,y_{2,n_{2}}x}^{c_{2}}}....\mathbf{T_{z_{k,1}x,y_{k,n_{k
}}x}^{c_{k}}}(b)$$
Similarly as before, we consider the $i$th cut, $c_{i}$, and label it
 with
$(u_{i},y_{i,n_{i}}x)$ for one component and
 $(u_{i}^{-1},y_{i,n_{i}}x)$
for the other component. Now we apply
$\mathbf{T_{z_{i,1}x,y_{i,n_{i}}x}^{c_{i}}}$ move first to get the
 label
$(u_{i},z_{i,1}x)$ and $(u_{i}^{-1},z_{i,1}x)$ respectively. Now we
 need to
apply $\mathbf{P_{x} \sqcup P_{x} \sqcup .... \sqcup P_{x}}$ move, and this will give us the label
$(xu_{i}x^{-1},z_{i,1})$ and $(xu_{i}^{-1}x^{-1},z_{i,1})$. Since this
 is a
closed loop, in particular this will imply $z_{i,1} =y_{i,n_{i}}x$ as
$i=1..k$. But we already know from above that $z_{i,1} =y_{i,n_{i}}$.
 So
this means $x = 1$. That is $\mathbf{P_{x}}$ = $\mathbf{P_{1}}$ =
 identity
or empty edge. Now our loops look like
$$\mathbf{F_{c_{1},z_{1,1}}^{-1}F_{c_{1},y_{1,n_{1}}}}
\mathbf{F_{c_{2},z_{2,1}}^{-1}F_{c_{2},y_{2,n_{2}}}}....\mathbf{F_{c_{k},z_
{k,1}}^{-1}F_{c_{k},y_{k,n_{k}}}} (b)$$
But since we already found out that $z_{i,1} =y_{i,n_{i}}$ as $i=1..k$.
 So
each
$\mathbf{F_{c_{i},z_{i,1}}^{-1}F_{c_{i},y_{i,n_{i}}}} = 1$.
So our loop is contractable.
\end{proof}
\subsubsection{The complex $M(\Sigma)$}
\begin{lemma}
The complex $M(\Sigma)$ is connected and simply-connected. Here the
 edges
of the complex are $Z,B,F$, define on the \cite{BK}.
Also the relation is defined on the same paper.
\end{lemma}
\begin{proof}
This is exactly the \cite{BK} is all about. So this paper is heavily
 depend
on this paper. We will not repeat the proof here. Interested readers
 are
referred to the paper \cite{BK}.\\
\end{proof}
\subsubsection{Proving part 3 of sec 9.1}
Let $$b'_{1} \stackrel{f'}{\rightarrow} b'_{2}$$ and
$$b''_{1}\stackrel{f''}{\rightarrow} b''_{2}$$ be two lifting of
 $$a_{1}
\stackrel{f}{\rightarrow} a_{2}$$ Then we need to show that there is a
path $e_{1}$, starting from $b''_{1}$ and end with $b'_{1}$, and
 completely
lie in $\Pi^{-1}(a_{1})$ and a path $e_{2}$, starting from $b''_{2}$
 and
end with $b'_{2}$, and completely lie in $\Pi^{-1}(a_{2})$, so that the
following relation hold:\\
\begin{center}
$e_{2}f''(b''_{1}) = f'e_{1}(b''_{1}) $
\end{center}
We know from lemma 8 that the vertex $b''_{1}$ and $b'_{1}$ is
 connected by
a sequence of move of the form $\mathbf{T_{z,y}^{c}}$ and
 $\mathbf{P_{x}}$
and their inverses. We choose such a path from $b''_{1}$ to $b'_{1}$.
 Let
us say the path is \\
$$(\mathbf{P_{x_{1}} \sqcup P_{x_{2}} \sqcup....\sqcup P_{x_{k}}}) \mathbf{T_{z_{1},y_{1}}^{c_{1}}}
\mathbf{T_{z_{2},y_{2}}^{c_{2}}}....\mathbf{T_{z_{k},y_{k}}^{c_{k}}}(b''_{1
})$$
For later reference, we will denote $(\mathbf{P_{x_{1}} \sqcup P_{x_{2}} \sqcup....\sqcup P_{x_{k}}})
\mathbf{T_{z_{1},y_{1}}^{c_{1}}}
\mathbf{T_{z_{2},y_{2}}^{c_{2}}}....\mathbf{T_{z_{k},y_{k}}^{c_{k}}} =
 L$.
So this path from $b''_{1}$ to $b'_{1}$ is just $L( b''_{1})$. Note
 that
this path $L( b''_{1})$ completely lie inside the fiber
 $\Pi^{-1}(a_{1})$.
So we get our $e_{1}$ which is just $L( b''_{1})$. Now depending on
 what
kind of path the $f$ is, we will have a different construction for the
 path
$e_{2}$.\\

\underline{case 1}: $f = Z$. In this case, we choose $e_{2} =
 L(b''_{2})$
more precisely, the path from $b''_{2}$ to $b'_{2}$ is $$e_{2} =
 L(b''_{2})
= (\mathbf{P_{x_{1}} \sqcup P_{x_{2}} \sqcup....\sqcup P_{x_{k}}}) \mathbf{T_{z_{1},y_{1}}^{c_{1}}}
\mathbf{T_{z_{2},y_{2}}^{c_{2}}}....\mathbf{T_{z_{k},y_{k}}^{c_{k}}}
(b''_{2})$$
Since in this case, both $f'$ and $f''$ is the $\mathbf{Z}$ move, the
commutativity that we want to show is
\begin{center}
$L\mathbf{Z}(b''_{1}) = \mathbf{Z}L(b''_{1}) $
\end{center}
Recall $(\mathbf{P_{x_{1}} \sqcup P_{x_{2}} \sqcup....\sqcup P_{x_{k}}}) \mathbf{T_{z_{1},y_{1}}^{c_{1}}}
\mathbf{T_{z_{2},y_{2}}^{c_{2}}}....\mathbf{T_{z_{k},y_{k}}^{c_{k}}} =
 L$.
But the $\mathbf{Z}$ move commute with all the moves, so in particular,
this means $L\mathbf{Z} = \mathbf{Z}L$ and we are done.\\

\underline{case 2}: $f = B$. In this case, again we choose $e_{2} =
L(b''_{2})$ more precisely, the path from $b''_{2}$ to $b'_{2}$ is
 $$e_{2}
= L(b''_{2}) = (\mathbf{P_{x_{1}} \sqcup P_{x_{2}} \sqcup....\sqcup P_{x_{k}}}) \mathbf{T_{z_{1},y_{1}}^{c_{1}}}
\mathbf{T_{z_{2},y_{2}}^{c_{2}}}....\mathbf{T_{z_{k},y_{k}}^{c_{k}}}
(b''_{2})$$
Since in this case, both $f'$ and $f''$ is the $\mathbf{B}$ move, the
commutativity that we want to show is
\begin{center}
$L\mathbf{B}(b''_{1}) = \mathbf{B}L(b''_{1}) $
\end{center}
here we suppress the indices for $\mathbf{B}$ since it is not
 important.
Again recall
$$(\mathbf{P_{x_{1}} \sqcup P_{x_{2}} \sqcup....\sqcup P_{x_{k}}}) \mathbf{T_{z_{1},y_{1}}^{c_{1}}}
\mathbf{T_{z_{2},y_{2}}^{c_{2}}}....\mathbf{T_{z_{k},y_{k}}^{c_{k}}} =
 L$$
But the $\mathbf{B}$ move also commute with the $\mathbf{T}$ and
$\mathbf{P}$ moves (see the $\mathbf{B}$ relation). So in particular,
 this
means $L\mathbf{B} = \mathbf{B}L$ and we are done.\\

\underline{case 3}: $f = F$. This situation is little bit different
 from
the above two. Here the $F$ move for $M(\Sigma)$ and the $\mathbf{F}$
 move
for $M(\ts ,\Sigma)$, both will remove a cut, say the $i$-th cut. Then
 we
choose\\
\begin{displaymath}
e_{2} = (\mathbf{P_{x_{1}} \sqcup P_{x_{2}} \sqcup...\sqcup \hat{P_{x_{i}}} \sqcup...\sqcup P_{x_{k}}}) \underbrace{\mathbf{T_{z_{1},y_{1}}^{c_{1}}}
\mathbf{T_{z_{2},y_{2}}^{c_{2}}}....\mathbf{T_{z_{k},y_{k}}^{c_{k}}}(b''_{2
})}_{\mbox{ but we do not include } \mathbf{T_{z_{i},y_{i}}^{c_{i}}}}
\end{displaymath}
Here $\hat{\mathbf{P_{x_{i}}}}$ means we do not include $\mathbf{P_{x_{i}}}$.
then the commutativity that we want show is:\\
\begin{displaymath}
\mathbf{F}(\mathbf{P_{x_{1}} \sqcup... P_{x_{k}}})\mathbf{T_{z_{1},y_{1}}^{c_{1}}}
\mathbf{T_{z_{2},y_{2}}^{c_{2}}}....\mathbf{T_{z_{k},y_{k}}^{c_{k}}} =
(\mathbf{P_{x_{1}} \sqcup...\sqcup \hat{P_{x_{i}}} \sqcup...\sqcup P_{x_{k}}}) \underbrace{\mathbf{T_{z_{1},y_{1}}^{c_{1}}}
\mathbf{T_{z_{2},y_{2}}^{c_{2}}}....\mathbf{T_{z_{k},y_{k}}^{c_{k}}}\mathbf
{F}}_{\mbox{ but we do not include } \mathbf{T_{z_{i},y_{i}}^{c_{i}}}}
\end{displaymath}
This is very easy to show. We will give the argument anyway for
completeness. We consider each cut one at a time. First consider the
 $i$ th
cut, the cut removed by $\mathbf{F}$. The left hand side of the above
equation will apply a bunch of $\mathbf{T}$ and a $\mathbf{P_{x}}$ move
 to
this cut first. But it does not matter, since at the end the
 $\mathbf{F}$
move will remove this cut. And the right hand side of the above euation
will apply $\mathbf{F}$ move first and remove this $i$-th cut. Now
 consider
the $j$-th cut where $i$ and $j$ are different. In this case, the
$\mathbf{F}$ move does not have any effect (does not change the label)
 on
the $j$-th cut. So both side of the above equation give rise to the
 same
label on this $j$-th cut (we can just forget about the $\mathbf{F}$
 move
from both side since it does not have any effect). This finishes the
 proof.
\subsubsection{Proving part 4 of sec 9.1}
\begin{figure}
\includegraphics{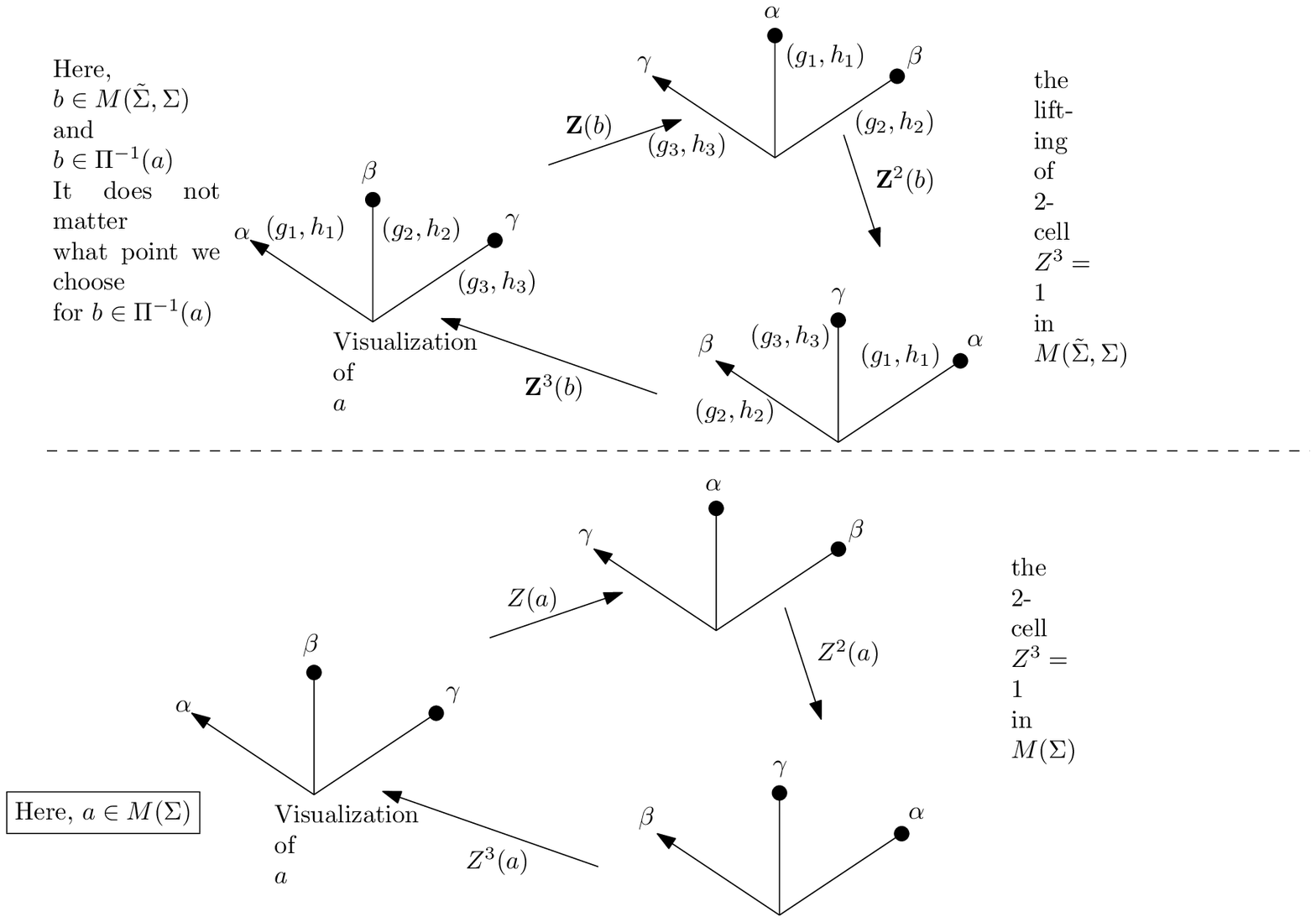}
\caption{Here we take $n = 3$ for simplicity of drawing diagram}\label{liftingofzcell}
\end{figure}

\begin{figure}
\includegraphics{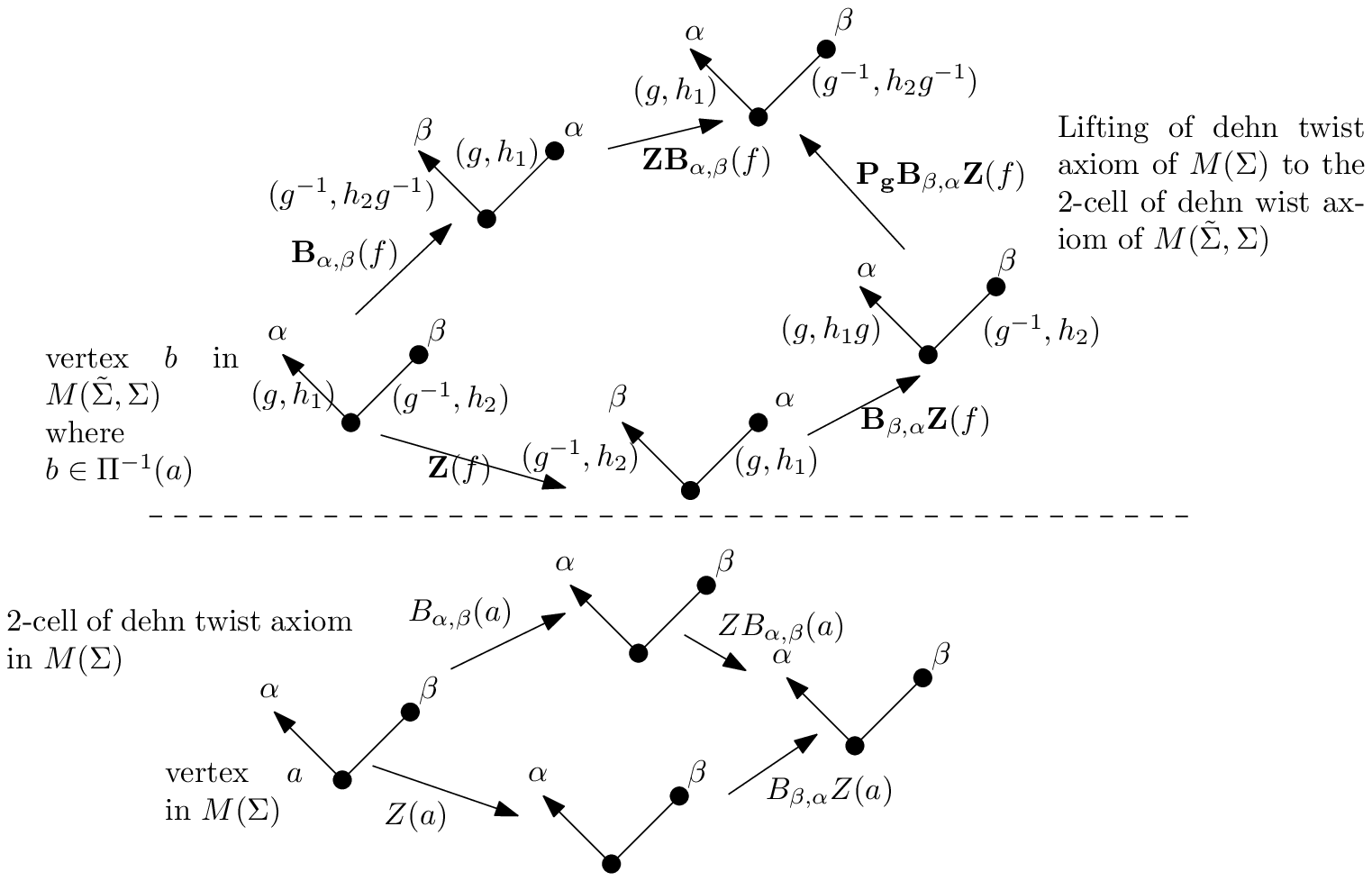}
\caption{}\label{liftingofdehntwist}
\end{figure}

We need to show that Every loop in $M(\Sigma)$ can be lifted to a
contractable loop in $M(\ts ,\Sigma)$. First note that, it is enough to
prove this for 2-cell or relations in $M(\Sigma)$. Because, then any
 loop
in $M(\Sigma)$ can be break down to several 2-cell in $M(\Sigma)$,
 since
$M(\Sigma)$ is simply-connected. We then lift each 2-cell to a
 contractable
loop in $M(\ts ,\Sigma)$. Then finally, we use part 3 of section 9.1,
 to
arrive at our answer. \\
Now we will show that every 2-cell in $M(\Sigma)$ can be lifted to a
contractable loop in $M(\ts ,\Sigma)$. We first observe that all the
 2-cell
in $M(\Sigma)$, consist of $Z,B$ and $F$ move. And each 2-cell has a
corresponding exact 2-cell in $M(\ts ,\Sigma)$. For example, consider
rotation axiom, the 2-cell in $M(\Sigma)$ is $Z^n = 1$ and the
corresponding 2-cell in $M(\ts ,\Sigma)$ is $\mathbf{Z}^{n} = 1$.
 Similarly
for the braiding axiom, the 2-cell of $M(\Sigma)$ is
$B_{\alpha,\gamma}B_{\alpha,\beta} = B_{\{\alpha\},\{\beta,\gamma\}}$
 and
$B_{\alpha,\gamma}B_{\beta,\gamma} = B_{\{\alpha\,\beta\},\{\gamma\}}$
and the corresponding 2-cell of $M(\ts ,\Sigma)$ is
$\mathbf{B_{\alpha,\gamma}B_{\alpha,\beta}} =
\mathbf{B_{\{\alpha\},\{\beta,\gamma\}}}$ and
$\mathbf{B_{\alpha,\gamma}B_{\beta,\gamma}} =
\mathbf{B_{\{\alpha\,\beta\},\{\gamma\}}}$. The only exception to this
 rule
is the dehn twist axiom. Dehn twist axiom for $M(\Sigma)$ is
$ZB_{\alpha,\beta} = B_{\beta,\alpha}Z$ and the corresponding dehn
 twist
axiom for $M(\ts ,\Sigma)$ is $\mathbf{ZB_{\alpha,\beta}} =
\mathbf{P_{g}B_{\beta,\alpha}Z}$, we have this extra $\mathbf{P_{g}}$
 move
appearing in the case of $M(\ts ,\Sigma)$. But this will not be a
 problem.
more precisely, we do the following:\\
If $L$ is a 2-cell in $M(\Sigma)$, say starting at the vertex $a$, then
 we
first pick any point on the fiber of $a$. Say $b \in \Pi^{-1}(a)$, it
 does
not matter which point on the fiber we choose. Then we apply the
corresponding 2-cell move of $M(\ts ,\Sigma)$ to $b$ and this will give
 a
contractable loop in $M(\ts ,\Sigma)$. For example if we have the
 2-cell
$Z^{n} (a)$ in $M(\Sigma)$ then a lifting of this 2-cell which is
contractable is going to be $\mathbf{{Z}^{n}}(b)$ in $M(\ts ,\Sigma)$.
 Look at the diagram on figure \ref{liftingofzcell}.
For the dehn twist, it is not a big difference. If $ZB_{\alpha,\beta} =
B_{\beta,\alpha}Z(a)$ is the 2-cell in $M(\Sigma)$ then
$\mathbf{ZB_{\alpha,\beta}} = \mathbf{P_{g}B_{\beta,\alpha}Z}(b)$ will
 be
the corresponding lifting in $M(\ts ,\Sigma)$ which is of course
contractible for the simple reason that $\mathbf{ZB_{\alpha,\beta}} =
\mathbf{P_{g}B_{\beta,\alpha}Z}$ is a 2-cell in $M(\ts ,\Sigma)$. See
 the diagram on figure \ref{liftingofdehntwist}.
So we satisfy all the condition of sec 9.1. This finishes the proof of
our main theorem.

\bibliographystyle{amsalpha}

\end{document}